\documentclass[11pt]{article} 
\usepackage[papersize={8.5in,11in}, left=1.2in, right=1.2in, top=1.1in, bottom=1.1in]{geometry} 

\usepackage{setspace}
\usepackage{indentfirst} 
\usepackage{xcolor} 
\usepackage{amsmath}
\usepackage{amssymb}
\usepackage{mathtools}
\usepackage{centernot} 
\usepackage{amsthm}
\usepackage{dsfont} 
\usepackage{authblk} 
\usepackage{graphicx} 
\usepackage{algorithmic,algorithm} 
\usepackage[hidelinks]{hyperref}
\usepackage{caption} 
\usepackage{array}
\newcolumntype{C}[1]{>{\centering\arraybackslash}m{#1}} 
\usepackage{blindtext}

\usepackage{doi}
\usepackage[numbers,sort]{natbib}

\setcitestyle{authoryear,open={(},close={)}} 
\setcitestyle{square,numbers} 
\usepackage{sectsty} 
\sectionfont{\large\bfseries}

\newcommand{\A}{\mathcal{A}}
\newcommand{\B}{\mathcal{B}}
\newcommand{\D}{\mathcal{D}}

\newcommand{\N}{\mathcal{N}}

\newcommand{\X}{\mathcal{X}}

\newcommand{\blambda}{\boldsymbol{\lambda}}
\newcommand{\bzeta}{\boldsymbol{\zeta}}

\newcommand{\bxi}{\boldsymbol{\xi}}

\newcommand{\argmin}{\mathop{\mathrm{arg\,min}{}}}
\newcommand{\argmax}{\mathop{\mathrm{arg\,max}{}}}

\newcommand{\bzero}{{\mathbf{0}}}
\newcommand{\bg}{{\mathbf g}}
\newcommand{\ba}{{\mathbf a}}

\newcommand{\bh}{{\mathbf h}}

\newcommand{\bx}{{\mathbf x}}
\newcommand{\bz}{{\mathbf z}}
\newcommand{\by}{{\mathbf y}}
\newcommand{\bu}{{\mathbf u}}
\newcommand{\bv}{{\mathbf v}}
\newcommand{\bb}{{\mathbf b}}

\newcommand{\bw}{{\mathbf w}}
\newcommand{\bA}{{\mathbf A}}
\newcommand{\bB}{{\mathbf B}}

\newcommand{\bO}{{\mathbf O}}

\newcommand{\vertiii}[1]{{\left\vert\kern-0.20ex\left\vert\kern-0.20ex\left\vert #1
		\right\vert\kern-0.20ex\right\vert\kern-0.20ex\right\vert}}

\theoremstyle{plain}
\newtheorem{theorem}{Theorem}[section]
\newtheorem{proposition}[theorem]{Proposition}
\newtheorem{lemma}[theorem]{Lemma}

\theoremstyle{definition}
\newtheorem{definition}[theorem]{Definition}
\newtheorem{assumption}[theorem]{Assumption}

\theoremstyle{remark}

\newcommand*{\affaddr}[1]{#1} 
\newcommand*{\affmark}[1][*]{\textsuperscript{#1}}
\newcommand*{\email}[1]{\texttt{#1}}

\title{Inexact Moreau Envelope Lagrangian Method for Non-Convex Constrained Optimization under Local Error Bound Conditions on Constraint Functions}
\author{%
Yankun Huang\affmark[1],~~Qihang Lin\affmark[2],~~Yangyang Xu\affmark[3]\\
\affaddr{\affmark[1]Department of Information System, Arizona State University, Tempe, AZ 85281}\\
\affaddr{\affmark[2]Department of Business Analytics, University of Iowa, Iowa City, IA 52242}\\
\affaddr{\affmark[3]Department of Mathematical Sciences, Rensselaer Polytechnic Institute, Troy, NY 12180}\\
\affmark[1]\email{yankun.huang@asu.edu},~~\affmark[2]\email{qihang-lin@uiowa.edu},~~\affmark[3]\email{xuy21@rpi.edu}\\
}
\date{}

\begin{document}

\maketitle
\begin{abstract}
In this paper, we investigate how structural properties of the constraint system impact the oracle complexity of smooth non-convex optimization problems with convex inequality constraints over a simple polytope. In particular, we show that, under a local error bound condition with exponent $d\in[1,2]$ on constraint functions, an inexact Moreau envelope Lagrangian method can attain an $\epsilon$-Karush--Kuhn--Tucker point with $\tilde O(\epsilon^{-2d})$ gradient oracle complexity. When $d=1$, this result matches the best-known complexity in literature up to logarithmic factors. Importantly, the assumed error bound condition with any $d\in[1,2]$ is strictly weaker than the local linear independence constraint qualification that is required to achieve the best-known complexity. Our results clarify the interplay between error bound conditions of constraints and algorithmic complexity, and extend complexity guarantees to a broader class of constrained non-convex problems.
\end{abstract}

\section{Introduction}

In this paper, we study first-order methods (FOMs) for nonlinear constrained optimization problems formulated as 
\begin{align}
\label{eq:gco_general}
\min_{\bx\in\X}f(\bx)~~\text{s.t.}~~ \bg(\bx)=(g_1(\bx),\dots,g_m(\bx))^\top\leq \bzero,
\end{align}
where $\X\subset\mathbb{R}^n$ is a bounded polytope that allows a computationally easy projection operator, $f$ is smooth but not necessarily convex, and $g_i$ for $i=1,\dots,m$ are smooth and convex. Since \eqref{eq:gco_general} is non-convex, computing its optimal or $\epsilon$-optimal solution 
is in general intractable. Therefore, the objective of this study is to compute an $\epsilon$-Karush--Kuhn--Tucker ($\epsilon$-KKT) point of \eqref{eq:gco_general} (see Definition~\ref{dfn:epsilon_KKT}) using an FOM.

The FOM studied in this paper is the inexact Moreau envelope Lagrangian (iMELa) method (see Algorithm~\ref{alg:imela}). This method 
is related to the proximal augmented Lagrangian method (ALM), a technique with a long history~\cite{hestenes1969multiplier,powell1969method,rockafellar1973multiplier,rockafellar1973dual,rockafellar1974augmented,rockafellar1976augmented} that remains an active area of research~\cite{lin2022complexity,kong2019complexity,kong2020efficient,melo2020iteration,kong2023iteration-SIOPT,kong2023iteration-MathOR,zhang2020proximal,zhang2022global,zhang2022iteration,pu2024smoothed,zeng2022moreau}. Following the literature, we measure the computational efficiency of the iMELa method by its first-order \emph{oracle complexity}, which is the number of gradient evaluations of $f$ or $g_i$ for $i=1,\dots,m$ required to reach an $\epsilon$-KKT point. 

Under a Slater's condition and a local \emph{error bound condition} (EBC) on the constraint set (Assumption~\ref{assume:local_EBC_exponent_d}) with exponent $d\in[1,2]$, we show that the iMELa method can find an $\epsilon$-KKT point for \eqref{eq:gco_general} with complexity $\tilde{O}(\epsilon^{-2d})$.\footnote{Here and in the rest of the paper, any logarithmic factor of $\epsilon$ is omitted in notation $\tilde{O}(\cdot)$.} When $d=1$, this complexity matches the optimal complexity $O(\epsilon^{-2})$ of an unconstrained problem~\cite{carmon2020lower,carmon2021lower} up to a logarithmic factor and thus is nearly optimal for a constrained problem. For constrained problems, the optimal complexity of $O(\epsilon^{-2})$ has been established in~\citet{zhang2020proximal,zhang2022global,zhang2022iteration,pu2024smoothed}
under different sets of assumptions. The assumptions used in this paper are not directly comparable to theirs:
some are weaker in certain aspects, while others are stronger. In other words, the complexity result in this work applies to the cases which are not covered by the existing methods. Moreover, our complexity $\tilde{O}(\epsilon^{-2d})$ with $d\in(1,2]$ is more general and new in the literature. A detailed comparison with ~\citet{zhang2020proximal,zhang2022global,zhang2022iteration,pu2024smoothed} and other existing methods is provided in Section~\ref{sec:relatedworks}. 

Beyond the complexity results, this work contributes to a deeper understanding
of the interaction between the convergence behavior of an FOM 
and EBCs of constraint functions.
Error bound conditions have a long history in optimization, due to their connections to constraint qualifications and their impact on the convergence of optimization algorithms (see, e.g.,
\citet{lewis1998error,luo1994error,pang1997error,li1997abadie}).
While EBCs on objective functions have been extensively
studied and exploited to accelerate the convergence of FOMs,
the role of EBCs on constraint functions has received
relatively limited attention. This work helps bridge this gap by showing how such conditions are utilized in the convergence analysis of an FOM for nonconvex constrained optimization.

\section{Related Works}
\label{sec:relatedworks}
In the recent literature~\cite{boob2023stochastic, ma2020quadratically,jia2025first,huang2023oracle,liu2025single}, subgradient-based methods have been developed for non-smooth continuous optimization problems with convex or weakly convex inequality constraints and achieved the oracle complexity of $O(\epsilon^{-4})$ for finding a nearly $\epsilon$-KKT point. Although those methods can be applied to \eqref{eq:gco_general}, their analysis does not utilize the smoothness of the functions and thus yields a complexity higher than the one in this work when the constraints are convex. 

The iMELa method we study is closely related to the classical  augmented Lagrangian method (ALM), which has a long history of study~\cite{hestenes1969multiplier,powell1969method,rockafellar1973multiplier,rockafellar1973dual,rockafellar1974augmented,rockafellar1976augmented} and is still one of the most effective approaches for constrained optimization. Below, we review the recent studies on ALM and its variants, including the penalty methods, for the following general problem
\begin{equation}
\label{eq:gco_general_new}
\begin{aligned}
    \min_{\bx}&~ f(\bx)+r(\bx) \\
    \text{s.t.}&~ \bg(\bx)=(g_1(\bx),\dots,g_m(\bx))^\top\leq \bzero, \\  
    &~\bh(\bx)=(h_1(\bx),\dots,h_l(\bx))^\top= \bzero, 
\end{aligned}
\end{equation}
where $f$, $g_i$, $i=1,\dots,m$, and $h_j$, $j=1,\dots,l$, are smooth and have Lipschitz continuous gradients, and $r$ is a proper convex lower semi-continuous function that allows a computationally easy proximal mapping.  When \eqref{eq:gco_general_new} is convex, the oracle complexity of the ALM for finding an $\epsilon$-optimal solution is well studied. See, e.g.~\citet{lan2016iteration,xu2021iteration,xu2021first,he2010acceleration,liu2019nonergodic} for the results under different settings. We next focus on the works that are applicable to \eqref{eq:gco_general_new} or its special cases when the problem is non-convex.  

In several studies on ALM (e.g.~\citet{hong2016decomposing,hajinezhad2019perturbed,zeng2022moreau}), the efficiency of an algorithm is characterized by its \emph{iteration complexity}, defined as the number of main iterations needed for computing an $\epsilon$-KKT point. Since a (proximal) augmented Lagrangian subproblem must be solved in each main iteration of ALM—typically using a separate FOM—the oracle complexity is generally higher than the iteration complexity. 

\textbf{Methods for linear constraints.} Suppose $g_i$'s are not present and $h_j$'s are linear in \eqref{eq:gco_general_new}. When $r\equiv 0$, \citet{hong2016decomposing} introduced a proximal primal-dual algorithm (prox-PDA) that finds an $\epsilon$-KKT point of \eqref{eq:gco_general_new} with an iteration complexity of $O(\epsilon^{-2})$. When $r$ in \eqref{eq:gco_general_new} is the characteristic function of a box or a bounded polytope, it has been shown in~\citet{zhang2020proximal,zhang2022global} that a smoothed proximal ALM (SP-ALM), which generalizes the classical proximal ALM~\cite{rockafellar1976augmented}, can achieve an $\epsilon$-KKT point with an oracle complexity of $O(\epsilon^{-2})$. When $r$ is the characteristic function of a compact set defined by convex inequalities, it has been shown that SP-ALM also achieves an $O(\epsilon^{-2})$ oracle complexity under CRCQ~\cite{zhang2022iteration}. 

Suppose no additional structural assumption is made on $r$. \citet{hajinezhad2019perturbed} developed a variant of prox-PDA with an iteration complexity of $O(\epsilon^{-4})$. \citet{zeng2022moreau} proposed the Moreau envelope ALM (MEAL), a type of proximal ALM, which achieves an iteration complexity of $o(\epsilon^{-2})$ when $f$ satisfies an implicit Lipschitz subgradient property and $O(\epsilon^{-2})$ when $f$ satisfies an implicit bounded subgradient property. The convergence property of MEAL is also established  when the augmented Lagrangian function satisfies the Kurdyka--{\L}ojasiewicz property~\cite{zeng2022moreau}. \citet{kong2019complexity,kong2020efficient} proposed a quadratic penalty accelerated inexact proximal point (QP-AIPP) method, which finds an $\epsilon$-KKT point with an oracle complexity of $\tilde O(\epsilon^{-3})$. Under additional mildly strong assumptions, a reduced oracle complexity of $\tilde O(\epsilon^{-2.5})$ is achieved by an inexact ALM (iALM)~\cite{li2021rate}, an inexact proximal accelerated augmented Lagrangian (IPAAL) method~\cite{melo2020iteration}, and an inner accelerated inexact proximal augmented Lagrangian (IAIPAL) method~\cite{kong2023iteration-SIOPT}. 


\textbf{Methods for nonlinear convex constraints.} Suppose $g_i$'s are convex and $h_j$'s are linear in \eqref{eq:gco_general_new}. \citet{lin2022complexity} proposed an inexact proximal point penalty (iPPP) method that achieves an oracle complexity of $\tilde{O}(\epsilon^{-2.5})$ under a Slater's condition. Under the similar assumptions, \citet{li2021augmented} extended this approach to a hybrid method combining ALM and the penalty method with the same complexity. \citet{dahal2023damped} provided a damped proximal ALM that achieves the same order of oracle complexity. The iALM by \citet{li2021rate} can also achieve the same complexity but requires an additional regularity assumption (see \eqref{eq:PLconstraint} below). The IAIPAL introduced in~\citet{kong2023iteration-SIOPT} was also extended by \citet{kong2023iteration-MathOR} to handle nonlinear convex constraints and obtains an $\tilde{O}(\epsilon^{-3})$ oracle complexity. 

When $r$ is the characteristic function of a compact set defined by convex inequalities, \citet{pu2024smoothed} developed a smoothed proximal Lagrangian method (SP-LM), which finds an $\epsilon$-KKT point with an oracle complexity of $O(\epsilon^{-2})$. In contrast, our method assumes $r$ is the characteristic function of a compact polytope and requires a complexity of $\tilde O(\epsilon^{-2})$. Besides the assumption on $r$, there are two key differences between their method and ours in the updating schemes. First, their method projects the dual variables onto a compact, artificially constructed set,  whereas our method only needs to project the dual variables onto $\mathbb{R}^m_+$. Second, their method only uses a single loop but ours employs a double-loop structure. More importantly, their method relies on a Slater-like condition and a local LICQ condition at all KKT points, whereas we assume a Slater's condition and a local error bound condition on each active subset of the constraint set (see Assumption~\ref{assume:local_EBC_exponent_d}). We will show (see Proposition~\ref{thm:LICQ_implies_local_EBC}) that, if the local LICQ condition in~\citet{pu2024smoothed} holds over all KKT points, our local error bound assumption also holds. However, the converse is not true. Thus, our work extends the complexity analysis of SP-LM to the case not covered in \citet{pu2024smoothed}.  

\textbf{Methods for nonconvex constraints.} Suppose $g_i$'s are non-convex or $h_j$'s are nonlinear in \eqref{eq:gco_general_new}. In this case, even finding a feasible solution for \eqref{eq:gco_general_new} is generally intractable. When $r$ is the characteristic function of a box and $\bg\equiv\bzero$,~\citet{curtis2015adaptive,curtis2016adaptive} proposed an augmented Lagrangian trust-region method with the global asymptotic convergence guarantee. Their methods may produce an infeasible stationary point due to the non-convex constraints. In general, additional assumptions are needed to obtain a feasible stationary point. For example, assuming a (nearly) feasible solution can be easily obtained, the iPPP method~\cite{lin2022complexity} and a scaled dual descent alternating direction method of multipliers (SDD-ADMM) proposed by \citet{sun2024dual} can find an $\epsilon$-KKT point with an oracle complexity of $O(\epsilon^{-4})$. 
Another common assumption made in literature is the regularity condition, which assumes that there exists $\zeta>0$ such that
\begin{align}
\label{eq:PLconstraint}
    \zeta\sqrt{\|[\bg(\bx)]_+\|^2+\|\bh(\bx)\|^2}
    \leq \text{dist}\left(-\nabla\bg(\bx)[\bg(\bx)]_+-\nabla\bh(\bx)\bh(\bx),\partial r(\bx)\right)
\end{align}
for any $\bx\in\text{dom}(r)$. When \eqref{eq:PLconstraint} holds on the constraints and $g_i$'s are not present, \citet{sahin2019inexact} showed that iALM can find an $\epsilon$-KKT point with a complexity of $\tilde O(\epsilon^{-4})$.  This complexity is reduced to $\tilde O(\epsilon^{-3})$ by the iALM~\cite{li2021rate}\footnote{Although the problem considered in~\citet{li2021rate} involves only equality constraints, their results can be extended to the general problem \eqref{eq:gco_general_new}, as noted in~\citet[Remark 6]{li2021rate}.} and the iPPP method~\cite{lin2022complexity} with the inequality constraints, and to $O(\epsilon^{-3})$ by the SDD-ADMM method~\cite{sun2024dual} without the inequality constraints when there is only one block in $\bx$. Moreover, the SDD-ADMM method~\cite{sun2024dual} can achieve a complexity of $O(\epsilon^{-2})$ under addtional stronger assumptions. 

\section{Preliminaries}
\label{sec:preliminaries}
For vectors, $\|\cdot\|$ denotes the Euclidean ($\ell_2$) norm and $\|\cdot\|_1$ denotes the $\ell_1$ norm. For matrices, $\|\cdot\|$ denotes the spectral norm. 
For a differentiable function $h$ on $\mathbb{R}^n$, we say $h$ is \emph{$L$-smooth ($L\geq0$)} on $\X$ if 
$   \|\nabla h(\bx) - \nabla h(\bx')\| \leq L \|\bx-\bx'\|$
for any $(\bx,\bx')\in\X\times\X$, \emph{$\mu$-strongly convex ($\mu\geq0$)} on $\X$ if 
$	h(\bx)\geq h(\bx')+\langle\nabla h(\bx'),\bx-\bx'\rangle+\frac{\mu }{2}\|\bx-\bx'\|^2$ 
for any $(\bx,\bx')\in\X\times\X$, and \emph{$\rho$-weakly convex} ($\rho\geq0$) on $\X$ if
$	h(\bx)\geq h(\bx')+\langle\nabla h(\bx'),\bx-\bx'\rangle-\frac{\rho}{2}\|\bx-\bx'\|^2$ 
for any $(\bx,\bx')\in\X\times\X$. Note that the $L$-smoothness of $h$ implies 
that 
$h$ is $L$-weakly convex but $h$ can have a smaller weak convexity constant.

For any $\bx\in\mathbb{R}^n$ and $\delta>0$, we denote by $\mathbb{B}(\bx,\delta)$ the Euclidean ball centered at $\bx$ with radius $\delta$. For a closed convex set $\B$, we denote its normal cone at $\bx$ by $\mathcal{N}_\B(\bx)$. 
Let $\delta_{\B}(\bx)$ be the zero-infinity characteristic function of $\B$, i.e., $\delta_{\B}(\bx)=0$ if $\bx\in \B$ and $+\infty$ otherwise. Let $\text{proj}_{\B}(\cdot)$ be the projection mapping onto $\B$ and $\text{dist}(\bx,\A):=\min_{\by\in \A}\|\bx-\by\|$ for a nonempty closed set $\A$. Define $[\cdot]_+:=\max\{\cdot,0\}$. Denote $[m]=\{1,\dots,m\}$ and $[l]=\{1,\dots,l\}$.


The following assumption is made throughout the paper.
\begin{assumption}
\label{assume:general}
The following statements hold:
\begin{itemize}
    \item[\textup{A}.] $\X=\big\{\bx\in\mathbb{R}^n:\bA\bx\leq \bb\big\}$ with $\bA\in\mathbb{R}^{l\times n}$ and $\bb\in\mathbb{R}^{l}$ is compact with diameter $D_\X:=\max_{\bx,\bx'\in\X}\|\bx-\bx'\|<+\infty$.
    \item[\textup{B}.] $\min_{\bx\in\X}f(\bx)\geq \underline{f}$ for some $\underline{f}\in\mathbb{R}$.
    \item[\textup{C}.] \vskip 0.5pt $f$ and $g_i$, $i\in[m]$, are $L$-smooth, and each $g_i$ is convex.
    \item[\textup{D}.] $\max_{\bx\in\X}\{\|\nabla f(\bx)\|\}\leq B_f$ for some $B_f\geq0$. 
    \item[\textup{E}.] $\max_{\bx\in\X}\max\{\|\bg(\bx)\|,\|\nabla \bg(\bx)\|\}\leq B_g$ for some $B_g\geq0$. 
    \item[\textup{F}.] (Slater's condition) There exists $\bx_{\textup{feas}}\in\X$ such that $g_i(\bx_{\textup{feas}})<0$, $i\in[m]$. 
\end{itemize}
\end{assumption}


We denote the Lagrangian function of \eqref{eq:gco_general} by
\begin{align}
\label{eq:Lagrange}
    \mathcal{L}(\bx,\blambda): =f(\bx)+\textstyle\sum_{i=1}^m \lambda_i g_i(\bx),
\end{align}
where $\blambda=(\lambda_1,\dots,\lambda_m)^\top$ is a vector of Lagrangian multipliers. We consider a numerical method for finding an $\epsilon$-KKT point of \eqref{eq:gco_general} defined as follows.

\begin{definition}
\label{dfn:epsilon_KKT}
A point $\bx\in\X$ is an $\epsilon$-KKT point of problem \eqref{eq:gco_general} if there exists $\blambda\geq\bzero$ such that 
\begin{align}
\label{dfn:epsilon_stat}
    \textup{dist}\big( -\nabla f(\bx) -\textstyle\sum_{i=1}^m \lambda_i\cdot \nabla g_i(\bx), \N_\X(\bx)\big) \leq\epsilon,\\
    \label{dfn:epsilon_feas}
    \|[\bg(\bx)]_+\|\leq\epsilon,\\
    \label{dfn:epsilon_comple_slack}
\textstyle\sum_{i=1}^m |\lambda_ig_i(\bx)|\leq\epsilon.
\end{align}
\end{definition}

For any $I_g\subseteq[m]$ and $I_{A}\subseteq[l]$, we define an active feasible subset as
\begin{align}
\label{eq:act_feas_subset} 
\hspace{-0.05in}
    \mathcal{S}(I_g, I_A) 
     :=\left\{\bu\in\X:
    \begin{array}{l}
    ~g_i({\bu})=0,\;i\in I_g,
    ~g_i({\bu})\leq 0,\;i\in[m]\backslash I_g,\\[1ex]
    ~[\bA{\bu}-\bb]_j=0,\;j\in I_A,
    ~[\bA{\bu}-\bb]_j\leq0,\;j\in[l]\backslash I_A
\end{array}
\right\}.
\end{align}
When $\mathcal{S}(I_g, I_A)\neq\emptyset$, we call $\mathcal{S}(I_g, I_A)$ the active feasible subset induced by 
$I_g$ and $I_A$. For any $\bx\in\X$, let 
$$
    J_g(\bx):=
    \{i\in[m]:  g_i(\bx)=0\},\quad
    J_A(\bx):=
    \{j\in[l]: [\bA\bx-\bb]_j=0\}
$$
be the index sets of the constraints active at $\bx$. Let $\X^*$ denote the set of KKT points of \eqref{eq:gco_general}. We assume the constraints of  \eqref{eq:gco_general} satisfy a uniform local EBC near each KKT point. 
\begin{assumption}[Local EBC with exponent $d$]
\label{assume:local_EBC_exponent_d}
There exist $\delta>0$, $\gamma>0$, and $d\in[1,2]$ such that, for any $\bx^*\in\X^*$ and any $\bx\in\mathbb{B}(\bx^*,\delta)$, it holds for every index pair $I=(I_g,I_A)$ with $I_g\subseteq J_g(\bx^*)$ and $I_A\subseteq J_A(\bx^*)$ that
\begin{align}
\label{eq:local_EBC_exponent_d}
    \textup{dist} \big( \bx,\mathcal{S}(I_g, I_A) \big)^{2d}
    \leq
    \gamma\left(
    \begin{array}{l}
    \sum_{i\in I_g} |g_i(\bx)|^2
    +\sum_{i\in[m] \backslash I_g} ([g_i(\bx)]_+)^2\\[1ex]
    +\sum_{j\in I_A}\left|[\bA\bx-\bb]_j\right|^2\\
    +\sum_{j\in[l] \backslash I_A} \big(\left[[\bA\bx-\bb]_j\right]_+\big)^2
    \end{array}
    \right).
\end{align}
\end{assumption}
Inequality~\eqref{eq:local_EBC_exponent_d} states that each active feasible set $\mathcal{S}(I_g, I_A)$ defined in~\eqref{eq:act_feas_subset} satisfies a local EBC with exponent $d\in[1,2]$, where the constants $\delta$ and $\gamma$ are uniform over all $\bx$ sufficiently close to $\X^*$ and all admissible index pairs $(I_g,I_A)$. Similar EBCs have been studied in constrained optimization; see, e.g.,~\citet{yang2018rsg,lin2025adaptive}. We further show that Assumption~\ref{assume:local_EBC_exponent_d} holds whenever the constraints in~\eqref{eq:gco_general} satisfy the \textbf{local linear independence constraint qualification} (LICQ) assumed in \citet[Assumption 2.2]{pu2024smoothed}. Specifically, the local LICQ requires that there exists $\zeta>0$ such that, for any $\bx^*\in\X^*$,
\begin{align}
\label{eq:LICQ_sigmamin}
    \sigma_{\min}\big(\big[\nabla \bg_{J_g(\bx^*)}(\bx^*), \bA_{ J_A(\bx^*)}^\top\big]\big)\geq \zeta,
\end{align}
where $\sigma_{\min}(\cdot)$ denotes the smallest singular value. Under \eqref{eq:LICQ_sigmamin}, each set $\mathcal{S}(I_g,I_A)$ satisfies a local EBC with $d=1$, which in turn implies Assumption~\ref{assume:local_EBC_exponent_d} for any $d\in[1,2]$. On the contrary,  Assumption~\ref{assume:local_EBC_exponent_d} does not imply the local LICQ and thus is strictly weaker. See the following proposition whose formal proof is provided in Appendix~\ref{sec:LICQ_implies_local_EBC}.
\begin{proposition}
\label{thm:LICQ_implies_local_EBC}
If the constraints of \eqref{eq:gco_general} satisfy the local LICQ in \eqref{eq:LICQ_sigmamin}, then Assumption~\ref{assume:local_EBC_exponent_d} holds, but not vice versa.
\end{proposition}

Throughout the paper, we assume that $p$ is a constant satisfying $p>L$. Following~\citet{zhang2020proximal,zhang2022global,zhang2022iteration}, given any $\bz\in\X$, we define 
\begin{align}
    \label{eq:v(z)}    
    v(\bz):=&~ \min_{\bx\in\X, \,g_i(\bx)\leq0,\,i\in[m]} \left\{f(\bx)+\frac{p}{2}\|\bx-\bz\|^2\right\},\\
    \label{eq:x(z)}
    \bx(\bz):=&~ \argmin_{\bx\in\X, \,g_i(\bx)\leq0,\,i\in[m]}\left\{f(\bx)+\frac{p}{2}\|\bx-\bz\|^2\right\},
\end{align}
and $\blambda(\bz)$ is the optimal Lagrangian multipliers corresponding to $\bx(\bz)$.
Note that \eqref{eq:v(z)} is a convex optimization problem as $p>L$ and $v(\bz)$ is differentiable by Danskin's Theorem (see, e.g.~\citet[Proposition B.25]{bertsekas1999nonlinear}) with
\begin{align}
\label{eq:nablav}
    \nabla v(\bz)=p(\bz-\bx(\bz)).
\end{align}
Like \eqref{eq:Lagrange}, the Lagrangian function of \eqref{eq:v(z)} is
\begin{align}
\label{eq:imela}
\hspace{-0.05in}
    \mathcal{L}_p (\bx,\bz,\blambda) :=f(\bx)+{\textstyle \sum_{i=1}^m\lambda_i g_i(\bx)}+\frac{p}{2}\|\bx-\bz\|^2.
\end{align}
Similar to the notations in~\citet{zhang2020proximal,zhang2022global,zhang2022iteration}, we define
\begin{align}
    \label{eq:d(lambda,z)}    
    d(\blambda,\bz):=&~ \min_{\bx\in\X} \mathcal{L}_p (\bx,\bz,\blambda),\\
    \bx(\blambda,\bz):=&~ \argmin_{\bx\in\X} \mathcal{L}_p (\bx,\bz,\blambda).
\end{align}
By the strong duality (see, e.g.~\citet[Sec. 28--30]{rockafellar1970convex}), it holds that
\begin{align}
\label{eq:str-dual}
    v(\bz)= \max_{\blambda\geq\bzero} d(\blambda,\bz).
\end{align} 

\section{Inexact Moreau Envelope Lagrangian Method}

We present the inexact Moreau envelope Lagrangian (iMELa) method, shown in Algorithm~\ref{alg:imela}, for finding an $\epsilon$-KKT point of~\eqref{eq:gco_general}. At iteration $t$, the method first updates the dual variable $\blambda^{(t+1)}$ via a projected gradient ascent step 
on function $\mathcal{L}_p(\bx^{(t)},\bz^{(t)},\blambda)$ at $\blambda=\blambda^{(t)}$ with step-size $\tau_t$. Then, given the strongly convex subproblem $\min_{\bx\in\X}\mathcal{L}_p(\bx,\bz^{(t)},\blambda^{(t+1)})$, we denote its  optimal solution by 
$$
\tilde{\bx}^{(t+1)}:= \bx(\blambda^{(t+1)},\bz^{(t)})
$$
in the rest of the paper for simplicity of notation. We approximately solve this subproblem by computing an $\epsilon_t$-stationary point $\bx^{(t+1)}\in\X$ with $\epsilon_t\in(0,1)$, satisfying
\begin{small}
\begin{equation}
\label{eq:subprob_optim_imela}
    \text{dist}\big( -\nabla_\bx\mathcal{L}_p(\bx^{(t+1)},\bz^{(t)},\blambda^{(t+1)}),\N_\X(\bx^{(t+1)})\big)\leq\epsilon_t.
\end{equation}
\end{small}Such a point can be computed using the accelerated projected gradient (APG) method by~\citet[Eq.\,(2.2.63)]{nesterov2018lectures} with complexity provided shortly. Moreover, condition~\eqref{eq:subprob_optim_imela} can be numerically verified and thus used as a stopping criterion for the APG method, provided that $\N_\X(\bx)$ has a simple structure. Finally, the proximal center $\bz$ is updated via 
a convex combination of $\bz^{(t)}$ and $\bx^{(t+1)}$ with a weight $\theta_t\in[0,1]$, ensuring $\bz^{(t)}\in\X$ for all $t\geq0$.

\begin{algorithm}[t]
\caption{Inexact Moreau Envelope Lagrangian (iMELa) Method}
\label{alg:imela}
\begin{algorithmic}[1]
    \STATE {\bfseries Input:} total number of iterations $T$, proximal parameter $p>0$, step-sizes $\{\tau_t,\theta_t\}_{t\geq0}$.
    \STATE {\bfseries Initialization:} $\bx^{(0)}=\bz^{(0)}\in\X$, $\blambda^{(0)}=\bzero$.
    \FOR{iteration $t=0,1,\dots,T-1$}
    \STATE 
    $\lambda_i^{(t+1)} =[\lambda_i^{(t)}+\tau_t\cdot g_i(\bx^{(t)})]_+,\; \forall\,i\in[m]$.
    \STATE 
    $\bx^{(t+1)} \approx \tilde{\bx}^{(t+1)}:=\argmin_{\bx\in\X}\mathcal{L}_p(\bx,\bz^{(t)},\blambda^{(t+1)})$ that satisfies \eqref{eq:subprob_optim_imela}.
    \STATE \vskip 0.5pt 
    $\bz^{(t+1)}= \bz^{(t)}+\theta_t\cdot(\bx^{(t+1)}-\bz^{(t)})$.
    \ENDFOR
\end{algorithmic}
\end{algorithm}


At iteration $t$ of the iMELa method, the complexity for the APG method by~\citet[Eq.\,(2.2.63)]{nesterov2018lectures} to ensure \eqref{eq:subprob_optim_imela} increases with the smoothness parameter of $\mathcal{L}_p(\bx,\bz^{(t)},\blambda^{(t+1)})$, which is $L+L \|\blambda^{(t+1)}\|_1+p$. Fortunately, the lemma below shows that a uniform upper bound on $\{\|\blambda^{(t)}\|\}_{t\geq 0}$ is available, whose proof is given in
Section~\ref{sec:bound_lambda_lm}.
\begin{lemma}
\label{lmm:bound_lambda}
The sequence $\{\blambda^{(t)}\}_{t\geq0}$ generated by Algorithm~\ref{alg:imela} satisfies, for any $t\geq0$,
\begin{align}
\label{eq:bound_lambda}
    \|\blambda^{(t)}\|\leq M_{\blambda}:=\max
    \left\{ 2\overline{\tau}B_g,\frac{2\overline{\tau}(C_{\blambda}+\overline{\tau}B_g^2)}{2\underline{\tau}\min_{i\in[m]}[-g_i(\bx_{\textup{feas}})]}\right\},
\end{align}
where $\overline{\tau}=\sup_{t\geq0}\tau_t$, $\underline{\tau}=\inf_{t\geq0}\tau_t$, and $C_{\blambda}:=(B_f+p D_{\X}+1)\cdot D_{\X}$.
\end{lemma}

Thanks to Lemma~\ref{lmm:bound_lambda} and the fact that $\|\blambda^{(t+1)}\|_1\leq\sqrt{m}\|\blambda^{(t+1)}\|\leq\sqrt{m}M_{\blambda}$, function $\mathcal{L}_p(\bx,\bz^{(t)},\blambda^{(t+1)})$ is $(p-L)$-strongly convex and $K$-smooth in $\bx$ over $\X$ for any $t\geq0$, where
\begin{align}
\label{eq:al_Lips_cont_grad}
    K:=&~L+L \sqrt{m} M_{\blambda}+p\geq L+L \|\blambda^{(t+1)}\|_1+p.
\end{align}
Given this result, the complexity for the APG method by~\citet[Eq.\,(2.2.63)]{nesterov2018lectures} to ensure \eqref{eq:subprob_optim_imela} is $O(\ln(\epsilon_t^{-1}))$. This is stated in the following lemma whose proof is given in Section~\ref{sec:subprob_complexity}.

\begin{lemma}
\label{lmm:subprob_complexity}
Suppose the APG method by~\citet[Eq.\,(2.2.63)]{nesterov2018lectures} is applied to the subproblem $\min_{\bx\in\X}\mathcal{L}_p(\bx,\bz^{(t)},\blambda^{(t+1)})$ for $k_t$ iterations with any initial solution in $\X$ and returns a solution $\bu^{(k_t)}$. Let 
$$
\bx^{(t+1)}=\textup{proj}_{\X} \big(\bu^{(k_t)}-K^{-1}\cdot\nabla \mathcal{L}_p(\bu^{(k_t)},\bz^{(t)},\blambda^{(t+1)})\big).
$$
Then $\bx^{(t+1)}$ satisfies \eqref{eq:subprob_optim_imela} with 
$$
\textstyle
k_t=O\left(\sqrt{\frac{K}{p-L}}\ln(\epsilon_t^{-1})\right).
$$
In other words, \eqref{eq:subprob_optim_imela} is ensured with complexity $O(\ln(\epsilon_t^{-1}))$ for any $t\geq0$.
\end{lemma}
The procedure described in Lemma~\ref{lmm:subprob_complexity} is formally stated in Algorithm~\ref{alg:apg} in the appendix and is implemented in our numerical experiments.

\section{Convergence Analysis for the iMELa Method}
For simplicity of notation, let 
$$
\sigma=(p-L)/p.
$$
The following lemma from~\citet[Lemma 3.5]{zhang2020proximal} characterizes the Lipschitz continuity of $\bx(\bz)$ in~\eqref{eq:x(z)} and $\bx(\blambda,\bz)$ in~\eqref{eq:d(lambda,z)} under Assumption~\ref{assume:general}.
\begin{lemma}
\label{lmm:Lips_continuity}
For any $\bz,\bz'\in\X$ and $\blambda\geq\bzero$, we have
\begin{align}
\label{eq:Lips_continuity_x(z)}
    \|\bz-\bz'\|\geq&~ \sigma\|\bx(\bz)-\bx(\bz')\|,\\
\label{eq:Lips_continuity_x(lambda,z)}
    \|\bz-\bz'\|\geq&~ \sigma\|\bx(\blambda,\bz)-\bx(\blambda,\bz')\|.
\end{align}
\end{lemma}

Similar to the analysis in~\citet{zhang2020proximal,zhang2022global,zhang2022iteration}, consider the potential function
\begin{align}
\label{eq:phi^t}
\hspace{-0.1in}
    \phi^t:=\mathcal{L}_p (\bx^{(t)},\bz^{(t)},\blambda^{(t)})-2d(\blambda^{(t)},\bz^{(t)})+2v(\bz^{(t)}).
\end{align}
By Assumption~\ref{assume:general}B, \eqref{eq:v(z)},  
\eqref{eq:d(lambda,z)} and \eqref{eq:str-dual}, we have
\begin{align}
    \nonumber
    \phi^t=&~
    (\mathcal{L}_p (\bx^{(t)},\bz^{(t)},\blambda^{(t)})-d(\blambda^{(t)},\bz^{(t)}))\\\nonumber
    &+(v(\bz^{(t)})-d(\blambda^{(t)},\bz^{(t)}))+v(\bz^{(t)})\\
    \label{eq:phi^t_lowerbound}
    \geq
    &~
    v(\bz^{(t)})\geq \underline{f}.
\end{align}

Following the analysis in~\citet{zhang2020proximal,zhang2022global,zhang2022iteration}, we characterize the change of $\phi^t$ after each iteration as follows. The proof is given in Section~\ref{sec:descent}.
\begin{proposition}
\label{thm:difference_phi}
Suppose that in Algorithm~\ref{alg:imela}, 
\begin{align*}
    \tau_t=\tau=\frac{p-L}{4B_g^2}\text{ and }\theta_t=\theta \le\frac{p-L}{18p}, ~\forall\,t\geq0.
\end{align*}
The sequence $\{(\bx^{(t)},\blambda^{(t)},\bz^{(t)})\}_{t\geq0}$ generated by Algorithm~\ref{alg:imela} satisfies that
\begin{align}
\nonumber
    \phi^t-\phi^{t+1}
    \geq&~ 
    \frac{p-L}{4}\|\bx^{(t)}-\tilde{\bx}^{(t+1)}\|^2-6p\theta\|\tilde{\bx}^{(t+1)}-\bx(\bz^{(t)})\|^2\\
    \label{eq:difference_phi}
    &+\frac{p}{6\theta}\|\bz^{(t)}-\bz^{(t+1)}\|^2+\tau\, \textup{dist}\big( \bg(\tilde{\bx}^{(t+1)}), \mathcal{N}_{\mathbb{R}_+^m}(\blambda^{(t+1)})\big)^2
    -\frac{\epsilon_t^2}{2(p-L)}.
\end{align}
\end{proposition}

 
According to Proposition~\ref{thm:difference_phi}, to ensure a sufficient decrease of $\phi_t$, it is necessary to control the negative term $-6p\theta\|\tilde{\bx}^{(t+1)}-\bx(\bz^{(t)})\|^2$ appearing in~\eqref{eq:difference_phi}. To this end, we employ the local EBC with exponent $d$ in Assumption~\ref{assume:local_EBC_exponent_d}. Specifically, we show in the following proposition that, when $\tilde{\bx}^{(t+1)}$ is a near-KKT point, the local EBC allows us to upper bound $\|\tilde{\bx}^{(t+1)}-\bx(\bz^{(t)})\|^2$ using a multiple of $\text{dist}( \bg(\tilde{\bx}^{(t+1)}), \mathcal{N}_{\mathbb{R}_+^m}(\blambda^{(t+1)}))^{2/d}$. The proof is presented in Section~\ref{sec:regularity}.

\begin{proposition}
\label{thm:regularity_exponent_d} 
Suppose that Assumption~\ref{assume:local_EBC_exponent_d} holds. 
There exists a constant $R(\delta)$, depending only on the constant $\delta$ in Assumption~\ref{assume:local_EBC_exponent_d}, such that, 
under the same conditions as in Proposition~\ref{thm:difference_phi}, if
\begin{align}
\nonumber
     &\textup{dist}\big(-\nabla f(\tilde{\bx}^{(t+1)}) - \nabla\bg(\tilde{\bx}^{(t+1)})\blambda^{(t+1)}, \N_\X(\tilde{\bx}^{(t+1)})\big)^2 \\
     \label{eq:regular_cond_KKT_near_R_exponent_d}
    &\quad+~\textup{dist}\big( \bg(\tilde{\bx}^{(t+1)}), \mathcal{N}_{\mathbb{R}_+^m}(\blambda^{(t+1)})\big)^2  \leq R^2(\delta),
\end{align}
it holds that
\begin{align}
\label{eq:regularity_result_exponent_d} 
    \|\bx(\bz^{(t)})-\tilde{\bx}^{(t+1)}\|^{2}
    \leq&~ \frac{\gamma^{1/d} K}{p-L} \textup{dist}\big( \bg(\tilde{\bx}^{(t+1)}), \mathcal{N}_{\mathbb{R}_+^m}(\blambda^{(t+1)})\big)^{2/d}.
    %
\end{align}
\end{proposition}

Thanks to Proposition~\ref{thm:regularity_exponent_d}, when the local regularity condition \eqref{eq:regular_cond_KKT_near_R_exponent_d} holds, the negative term $-6p\theta\|\tilde{\bx}^{(t+1)}-\bx(\bz^{(t)})\|^2$ in~\eqref{eq:difference_phi} can be bounded above by a multiple of the quantity $\text{dist}( \bg(\tilde{\bx}^{(t+1)}), \mathcal{N}_{\mathbb{R}_+^m}(\blambda^{(t+1)}))^{2/d}$ as long as $\theta_t=\theta$ is chosen sufficiently small. This bound follows directly from the local EBC with exponent $d$ in Assumption~\ref{assume:local_EBC_exponent_d}.

However, when \eqref{eq:regular_cond_KKT_near_R_exponent_d} does not hold, the local EBC is unavailable. In this case, we resort to the following global result, which holds without any error bound assumption. See Section~\ref{sec:weak_dual_bound} for its proof.

\begin{lemma}
\label{lmm:weak_dual_bound}
With $M_{\blambda}$ defined in \eqref{eq:bound_lambda}, it holds that
\begin{align}
\label{eq:weak_dual_bound}
    \|\tilde{\bx}^{(t+1)} - \bx(\bz^{(t)})\|^2\leq \frac{2M_{\blambda}}{p-L}
    \textup{dist} \big( \bg(\tilde{\bx}^{(t+1)}), \mathcal{N}_{\mathbb{R}_+^m}(\blambda^{(t+1)})\big). 
\end{align}
\end{lemma}
The result of Lemma~\ref{lmm:weak_dual_bound} can be viewed as a global analogue of Assumption~\ref{assume:local_EBC_exponent_d} with exponent $d=2$, and serves as a fallback estimate when condition~\eqref{eq:regular_cond_KKT_near_R_exponent_d} fails.
Using Proposition~\ref{thm:regularity_exponent_d} and Lemma~\ref{lmm:weak_dual_bound}, we can ensure a sufficient decrease of $\phi_t$ by choosing $\theta_t=\theta$ and $\tau_t$ in Algorithm~\ref{alg:imela} carefully. This result is presented below and is the key to establish the total complexity. The proof is given in Section~\ref{sec:bounded_phi_proof}.

\begin{proposition}
\label{thm:bounded_phi_exponent_d}
Suppose that Assumptions~\ref{assume:general}
and~\ref{assume:local_EBC_exponent_d} hold and, in Algorithm~\ref{alg:imela}, we set 
\begin{align}
\label{eq:tilde_theta_d}
\begin{aligned}
    \tau_t=\tau=&~ \frac{p-L}{4B_g^2},\\
    \theta_t=\theta\leq &~\tilde{\theta}_d:= \min\left\{\frac{p-L}{18 p}, 6p\tau,  \frac{R^2(\delta) (p-L)^2\tau}{2\cdot 12^3 p^3 M_{\blambda}^2}, \hat{\theta}_{d}\right\},
\end{aligned}
\end{align}
for $\forall\, t\geq0$, where $C_d:=\frac{\gamma^{1/d}K}{p-L}$ and 
\begin{equation}
\label{eq:hat_theta_d}
    \hat{\theta}_{d}:=
    \frac{\tau}{12p C_d}
    \left(\frac{\tau}{12p}\right)^{1-1/d}.
\end{equation} 
Then the sequence $\{(\bx^{(t)},\blambda^{(t)},\bz^{(t)})\}_{t\geq0}$ generated by Algorithm~\ref{alg:imela} satisfies that
\begin{align}
\nonumber
    \phi^t-\phi^{t+1}
    \geq&~\frac{p-L}{4}\|\bx^{(t)}-\bx(\blambda^{(t+1)},\bz^{(t)})\|^2+\frac{p}{12\theta}\|\bz^{(t)}-\bz^{(t+1)}\|^2\\
    \label{eq:bounded_phi_exponent_d}
    &+\frac{\tau}{2} \textup{dist}\big( \bg(\bx(\blambda^{(t+1)},\bz^{(t)})), \mathcal{N}_{\mathbb{R}_+^m}(\blambda^{(t+1)})\big)^2
    - \Xi_{d}(\theta)- \frac{p\theta\epsilon_t^2}{12(p-L)^2},
\end{align}
where
\begin{equation}
\label{eq:Xi_d_def}
    \Xi_d(\theta)
    :=
    \begin{cases}
    \left(1-\dfrac{1}{d}\right) \left(\dfrac{2}{\tau d}\right)^{\frac{1}{d-1}} \big( 6p\theta\, C_d \big)^{\frac{d}{d-1}},
    & d\in(1,2],\\[1ex]
    0,
    & d=1.
    \end{cases}
\end{equation}
\end{proposition}

By Proposition~\ref{thm:bounded_phi_exponent_d}, we establish the complexity of Algorithm~\ref{alg:imela} in the following theorem whose proof is in Section~\ref{sec:main_result}.
\begin{theorem}
\label{thm:main_result_exponent_d}
Suppose that the assumptions in Proposition~\ref{thm:bounded_phi_exponent_d} hold with a fixed exponent $d\in[1,2]$.
For a given target accuracy $\epsilon>0$, choose $\theta_t=\theta=\theta(\epsilon)$, where
\begin{align*}
    \theta(\epsilon)
    :=
    \min\left\{
    \tilde{\theta}_d,\,
    \frac{6p\tau}{M_{\blambda}^2+1},\,
    \left(\frac{\epsilon^2}{2A_d}\right)^{d-1}
    \right\},~
    \text{ with }~
    A_d:=24p\left(1-\frac{1}{d}\right)
    \left(\frac{2}{\tau d}\right)^{\frac{1}{d-1}}
    (6pC_d)^{\frac{d}{d-1}},
\end{align*}
and $\tilde{\theta}_d$, $C_d$ (and $\tau$) defined in Proposition~\ref{thm:bounded_phi_exponent_d}. If  $\epsilon_t=\frac{c}{t+1}$ for any constant $c>0$ and
\begin{align}
\label{eq:C0C1_def_exponent_d}
\begin{aligned}
    T :=
    \left\lceil \frac{2\mathcal C_0(\theta(\epsilon))}{\epsilon^{2}} \right\rceil, \text{ where }~
    \mathcal C_0(\theta)
    :=&~
    \left(
    \begin{array}{l}
        \left( 1+\frac{(M_{\blambda}^2+1)B_g^2}{(p-L)^2} \right)\frac{c\pi^2}{3}\\
        +\;\mathcal C_1(\theta) 
        \left( \phi^0-\underline{f}+\frac{cp\theta\pi^2}{72(p-L)^2} \right)
    \end{array}
    \right),\\
    \mathcal C_1(\theta)
    :=&~ \max\left\{\frac{24p}{\theta},\frac{4(M_{\blambda}^2+1)}{\tau}\right\},
\end{aligned}
\end{align}
Algorithm~\ref{alg:imela} finds an $\epsilon$-KKT solution for \eqref{eq:gco_general} with oracle complexity $O(\epsilon^{-2d}\ln(\epsilon^{-1}))$.
\end{theorem}

\section{Numerical Experiment}
\label{sec:experiment}
We demonstrate the performance of the iMELa method on a fairness-aware classification problem, which is an instance of \eqref{eq:gco_general}. We compare it 
with the inexact proximal point penalty (iPPP) method~\cite{lin2022complexity}, the damped proximal augmented Lagrangian method (DPALM) \citep{dahal2023damped},
the smoothed proximal Lagrangian method (SP-LM)~\cite{pu2024smoothed}, and the switching subgradient (SSG) method~\cite{huang2023oracle}. 

First, we define the instance of~\eqref{eq:gco_general} solved in our experiments. Given a feature vector $\ba\in\mathbb{R}^d$ and a class label $b\in\{1,-1\}$, the goal of linear binary classification is to learn a model $\bx\in\mathbb{R}^d$ to predict $b$ based on the score $\bx^\top\ba$ with a larger score indicating a higher chance of a positive label. Suppose that a training set $\D=\{(\ba_i,b_i)\}_{i=1}^n$ is available. A model $\bx$ can be obtained by solving 
\begin{equation}
\label{eq:ermL_linear}
    \mathcal{L}^* := \min_{\bx\in\X}\Big\{\mathcal{L}(\bx):=\frac{1}{n}\sum_{i=1}^n \ell(b_i\cdot\bx^\top\ba_i)\Big\},
\end{equation}
where $\ell(z)=\log(1+\exp(-z))$ and $\X$ is a convex and compact set.

However, solving~\eqref{eq:ermL_linear} only optimizes the classification accuracy of the resulting model but does not guarantee its fairness. Suppose that a data point $\ba$ is classified as positive if $\bx^\top\ba\geq0$ and as negative otherwise. Also, suppose that model $\bx$ is applied to two groups of data, a protected group $\mathcal{D}_p=\{\ba_i^p\}_{i=1}^{n_p}$ and an unprotected group $\mathcal{D}_u=\{\ba_i^u\}_{i=1}^{n_u}$. A measure of the fairness of $\bx$ between $\mathcal{D}_p$ and $\mathcal{D}_u$ is defined as 
$$
    \Big|\frac{1}{n_p}\sum_{i=1}^{n_p} \mathbb{I}(\bx^\top\ba_i^p\geq0)-\frac{1}{n_u}\sum_{i=1}^{n_u} \mathbb{I}(\bx^\top\ba_i^u\geq0) \Big|,
$$
where $\mathbb{I}(\cdot)$ is the one-zero indicator function. This measure is known as the demographic parity ~\cite{feldman2015certifying}. When its value is small, the model $\bx$ produces similar predicted positive rates in both groups, indicating the fairness of the model. However, this measure is computationally challenging because of the discontinuity of $\mathbb{I}(\cdot)$. Hence, we approximate this measure by a continuous function of $\bx$ defined as
\begin{equation}
\label{eq:ROC_linear}
    \mathcal{R}(\bx):=
    \frac{1}{n_p}\sum_{i=1}^{n_p} \sigma(\bx^\top\ba_i^p)-\frac{1}{n_u}\sum_{i=1}^{n_u} \sigma(\bx^\top\ba_i^u),
\end{equation}
where $\sigma(z)=\exp(z)/(1+\exp(z))$. 

\begin{figure*}[!th]
     \begin{tabular}{@{}c|ccc@{}}
      & a9a & bank & COMPAS \\
		\hline \vspace*{-0.1in}\\
		\raisebox{10ex}{\small{\rotatebox[origin=c]{90}{Objective}}}
		& \hspace*{-0.06in}\includegraphics[width=0.30\textwidth]{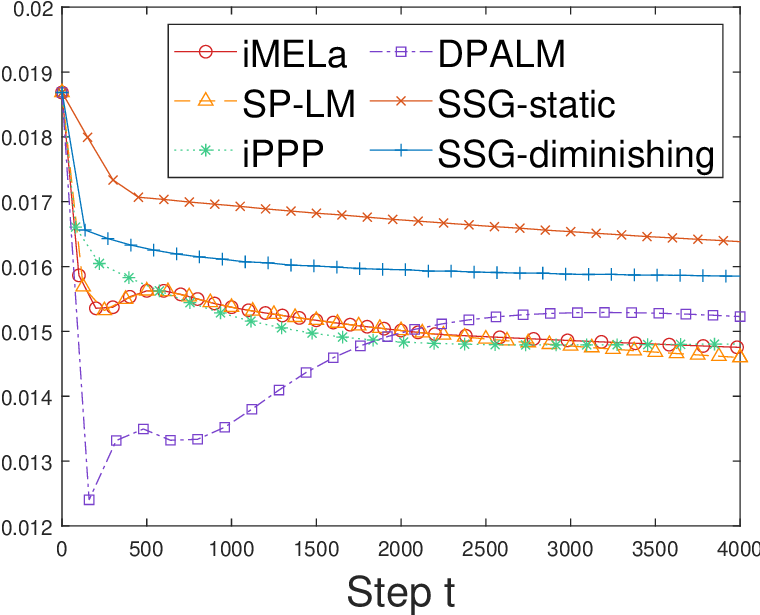}
		& \hspace*{-0.06in}\includegraphics[width=0.295\textwidth]{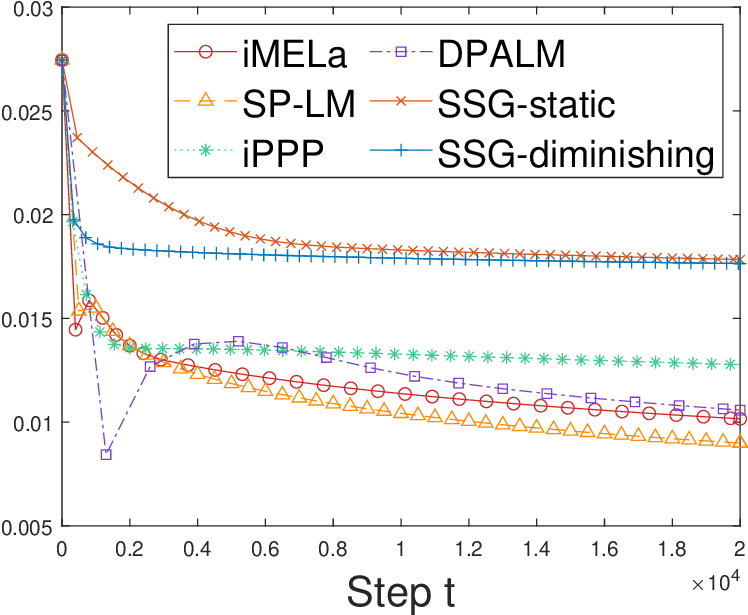}
		& \hspace*{-0.06in}\includegraphics[width=0.29\textwidth]{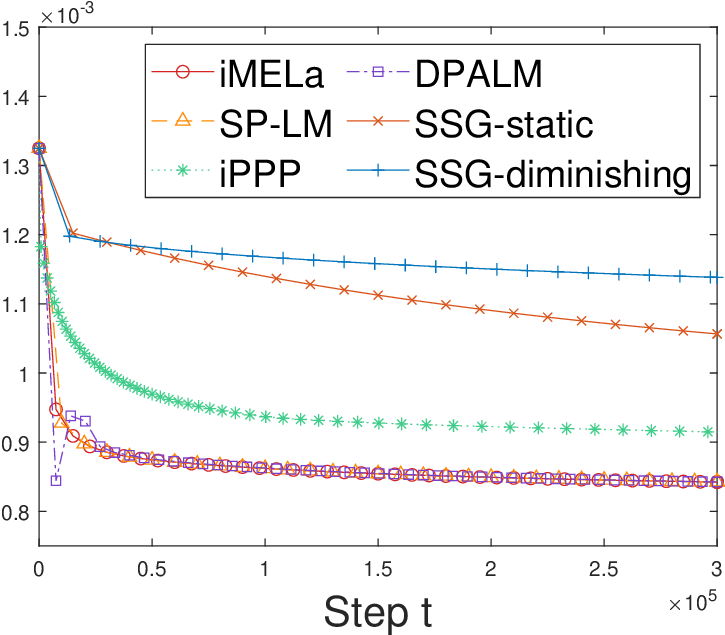}
          \\
		\raisebox{10ex}{\small{\rotatebox[origin=c]{90}{Infeasibility}}}
		& \hspace*{-0.06in}\includegraphics[width=0.30\textwidth]{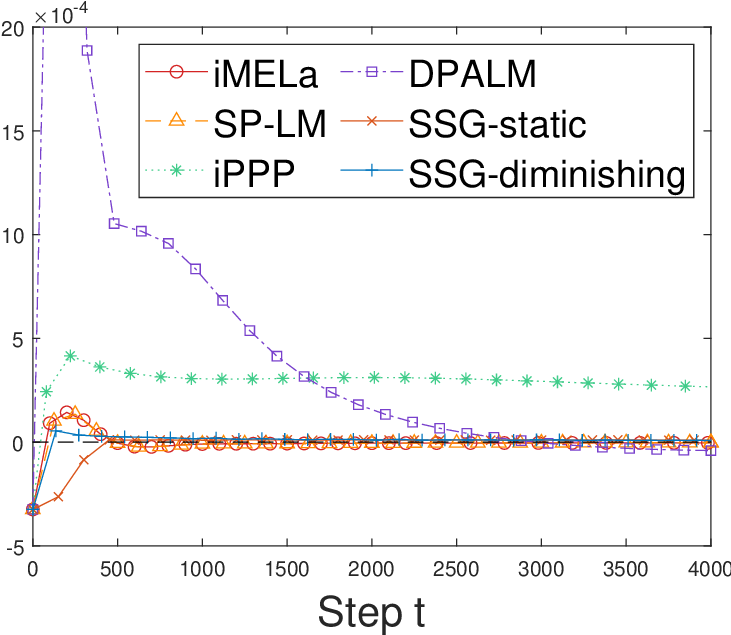}
		& \hspace*{-0.06in}\includegraphics[width=0.295\textwidth]{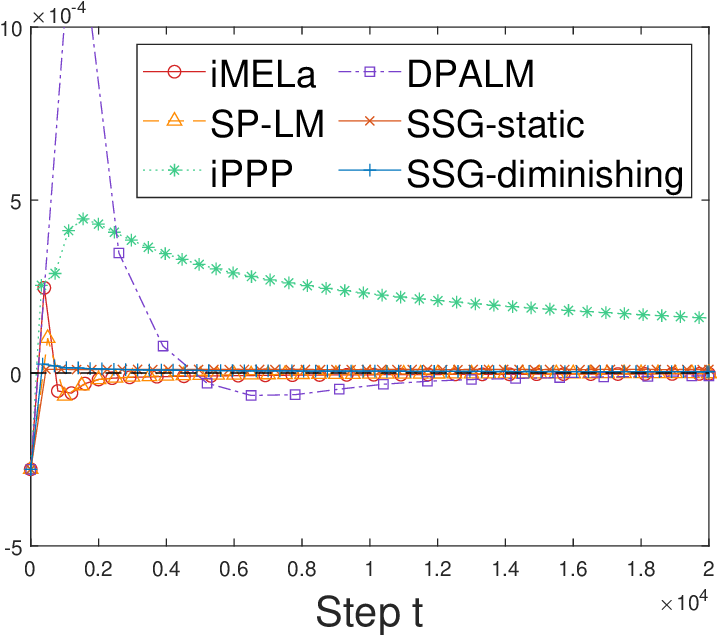}
		& \hspace*{-0.06in}\includegraphics[width=0.295\textwidth]{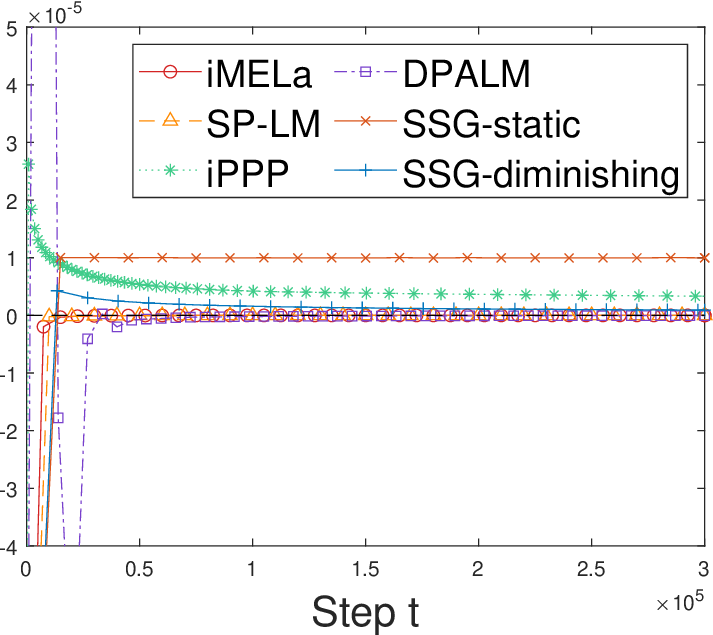}
           \\
        \raisebox{10ex}{\small{\rotatebox[origin=c]{90}{Stationarity}}}
		& \hspace*{-0.06in}\includegraphics[width=0.305\textwidth]{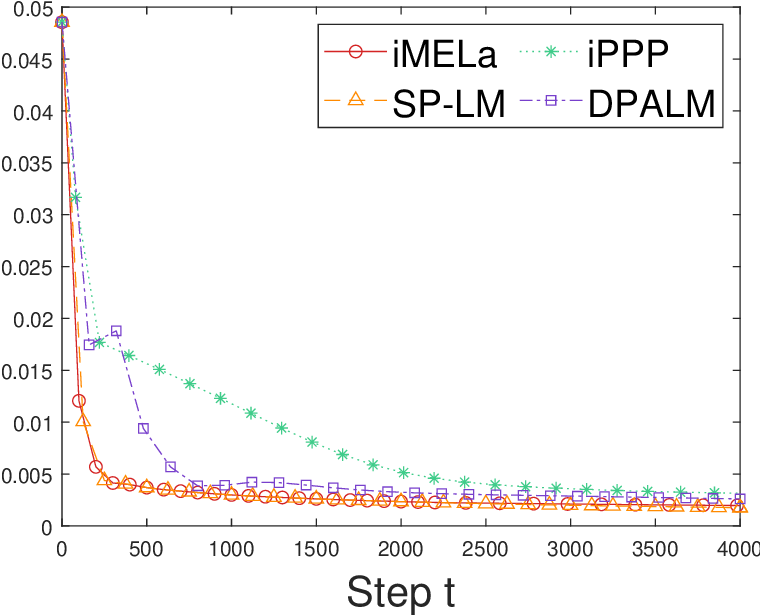}
		& \hspace*{-0.06in}\includegraphics[width=0.295\textwidth]{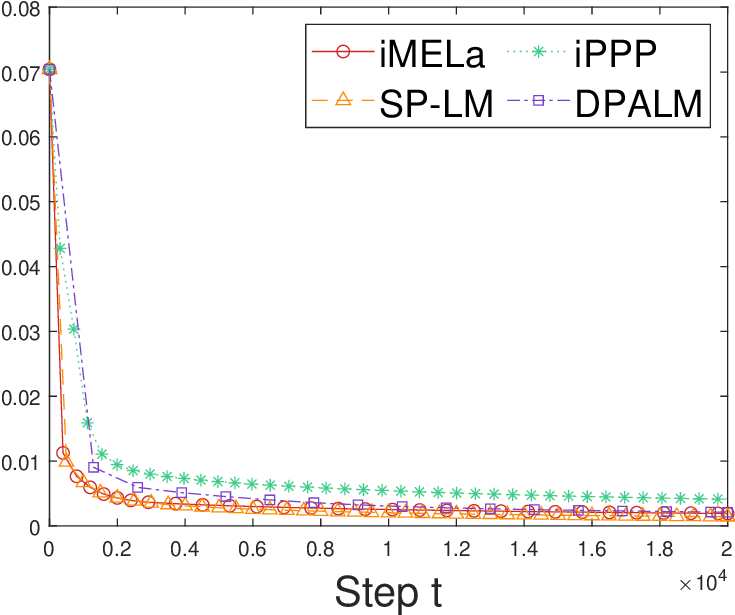}
		& \hspace*{-0.06in}\includegraphics[width=0.295\textwidth]{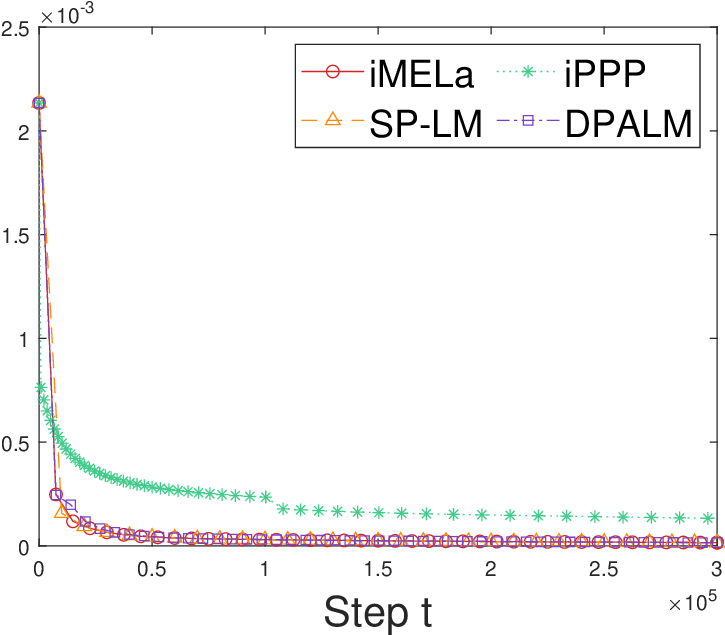}
          \\
        \raisebox{10ex}{\small{\rotatebox[origin=c]{90}{Complementary Slackness}}}
		& \hspace*{-0.06in}\includegraphics[width=0.30\textwidth]{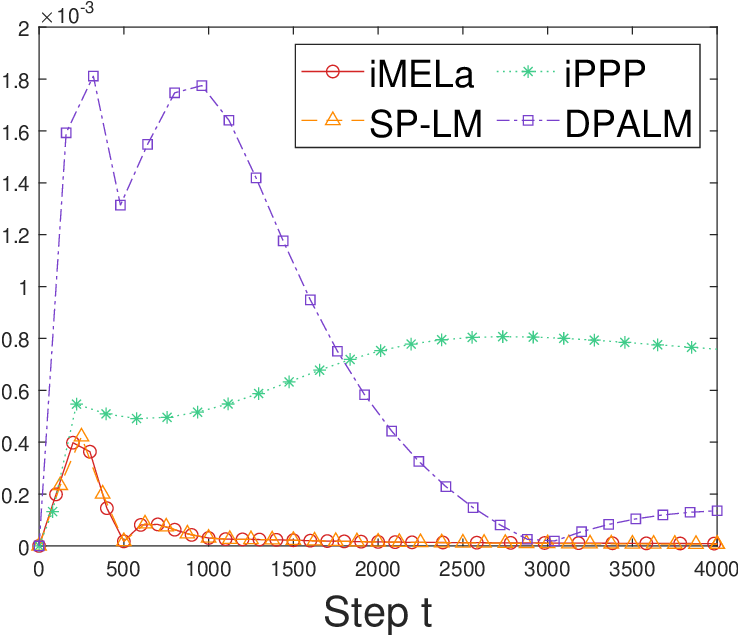}
		& \hspace*{-0.06in}\includegraphics[width=0.295\textwidth]{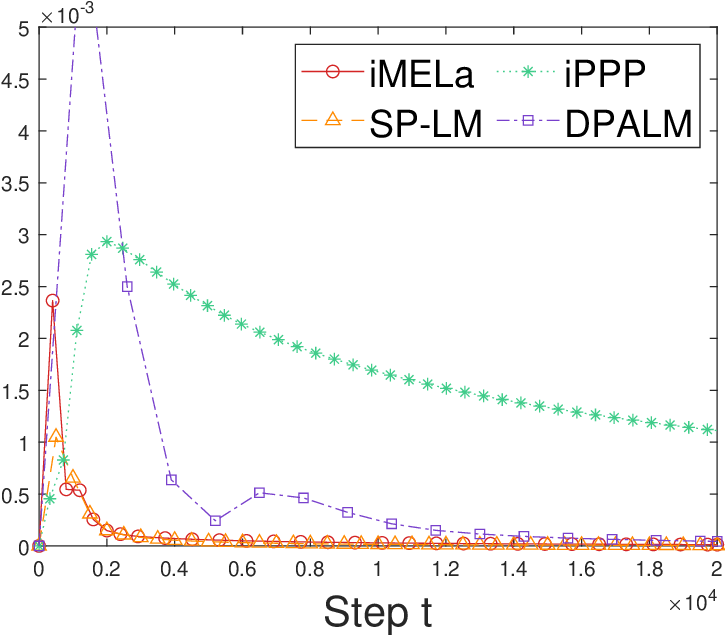}
		& \hspace*{-0.06in}\includegraphics[width=0.295\textwidth]{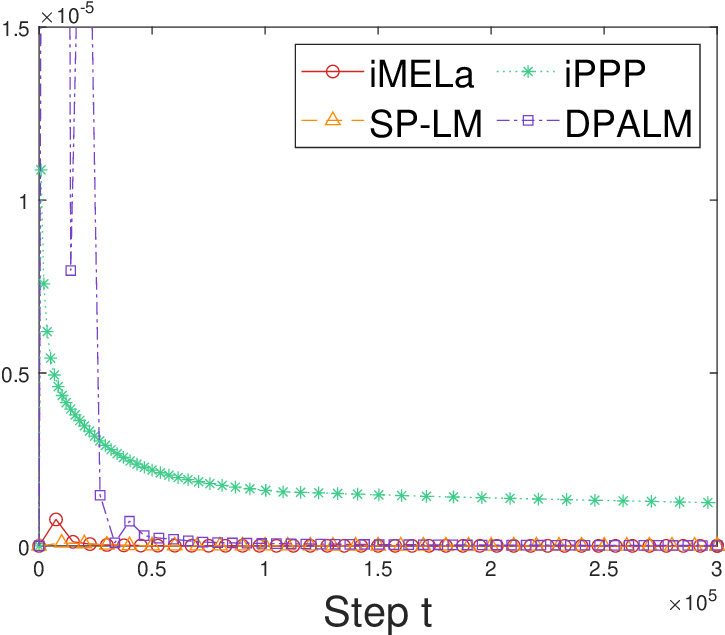}
        \end{tabular}
	\caption{Performances of $\text{iMELa}$, $\text{SP}$-$\text{LM}$, $\text{iPPP}$, $\text{DPALM}$, and $\text{SSG}$ with static and diminishing step-sizes vs. number of steps on classification problems with demographic parity fairness.} 
    \label{fig:figure_convex_experiment_step}
	\vspace{-0.1in}
\end{figure*}

To obtain a fair $\bx$, we balance the loss $\mathcal{L}(\bx)$ in \eqref{eq:ermL_linear} and the fairness measure $|\mathcal{R}(\bx)|$ in \eqref{eq:ROC_linear} by solving 
\begin{equation}
\label{eq:DPfairnessclassification_linear}
    \min_{\bx\in\X} 
    \frac{1}{2}\big(\mathcal{R}(\bx)\big)^2\text{ s.t. } \mathcal{L}(\bx)\leq \mathcal{L}^*+\kappa,
\end{equation}
where $\kappa$ is the slackness parameter indicating how much we are willing to increase the loss in order to reduce $|\mathcal{R}(\bx)|$ to obtain a more fair model, and $\X$ is the same 
as that in \eqref{eq:ermL_linear}. In the experiments, we set $\X=\{\bx\in\mathbb{R}^d:\|\bx\|_1\leq r\}$ with $r=6\|\ba_i\|_1$ to satisfy Assumption~\ref{assume:general}A.

We now briefly justify that problem \eqref{eq:DPfairnessclassification_linear} satisfies Assumption~\ref{assume:local_EBC_exponent_d} with exponent $d=1$. Specifically, we show that the constraint system satisfies the local LICQ condition at all KKT points, which implies the desired local EBC with exponent $d=1$ via Proposition~\ref{thm:LICQ_implies_local_EBC}. The formal statement and proof are deferred to Section~\ref{sec:DP_LICQ}.

We solve problem \eqref{eq:DPfairnessclassification_linear} on three datasets: \textit{a9a}~\cite{kohavi1996scaling}, \textit{bank}~\cite{moro2014data} and \textit{COMPAS}~\cite{angwin2016compas}.  Details about these datasets and how they are used to formulate \eqref{eq:DPfairnessclassification_linear} are provided in Section~\ref{sec:dataset}. The implementation details of each method, including the choices of tuning parameters, are provided in Section~\ref{sec:alg_details}.

We compare the methods based on the objective value, infeasibility, stationarity, and complementary slackness (see \eqref{dfn:epsilon_stat}, \eqref{dfn:epsilon_feas} and \eqref{dfn:epsilon_comple_slack}) achieved by $\bx^{(t)}$ and the corresponding Lagrangian multipliers.\footnote{Since the SSG method does not generate a Lagrangian multiplier, we do not report its measures of stationarity and complementary slackness.} 
While the calculations for most of these metrics are straightforward, additional details are provided for computing the stationarity metric. Observe that $\X=\{\bx\in\mathbb{R}^d:\bB\bx\leq r\cdot\mathbf{1}\}$, where $\bB\in\mathbb{R}^{2^d\times d}$ is a matrix whose rows are all vectors in $\{1,-1\}^{d}$, and $\mathbf{1}$ denotes the all-one vector in $\mathbb{R}^{2^d}$. If $\|\bx\|_1<r$, then $\N_{\X}(\bx)=\{\bzero\}$ so that the stationarity measure in \eqref{dfn:epsilon_stat} reduces to $\|\mathcal{R}(\bx)\cdot\nabla \mathcal{R}(\bx)+\lambda\cdot\nabla \mathcal{L}(\bx)\|$, where $\lambda\geq0$ is the Lagrange multiplier produced by the algorithm to pair with $\bx$. If $\|\bx\|_1=r$, then by~\citet[Example A.5.2.6 (b)]{hiriart2001fundamentals}, there exists an index set of the active constraints at $\bx\in\X$ defined by
\begin{align}
\label{eq:activeset}
    J(\bx)=\{j=1,\dots,2^d:[\bB]_j\bx=r\}
\end{align}
where $[\bB]_j$ denotes the $j$-th row of $\bB$. Then the normal cone to $\X$ at $\bx$ in this case is
\begin{align}
\nonumber
    \N_{\X}(\bx)=&~
    \text{cone}\{[\bB]_j^\top:j\in J(\bx)\} \\
    \label{eq:normalcone}
    =&~ 
    \textstyle\{\sum_{j\in J(\bx)}a_j[\bB_j]^\top:a_j\geq0,\,j\in J(\bx)\}.
\end{align}
Therefore, the stationarity measure in \eqref{dfn:epsilon_stat} can be calculated by solving 
\begin{align*}
    \min_{a_j\geq0,\,j\in J(\bx)}\left\|\mathcal{R}(\bx)\cdot\nabla \mathcal{R}(\bx)+\lambda\cdot\nabla \mathcal{L}(\bx)+\textstyle \sum_{j\in J(\bx)}a_j[\bB_j]^\top\right\|.
\end{align*}
We use \texttt{mpcActiveSetSolver} in MATLAB to solve this quadratic program.

To compare the methods in their oracle complexity, in Figure~\ref{fig:figure_convex_experiment_step}, we plot the four aforementioned performance metrics against the accumulated number of gradient steps performed by the APG method called within the iMELa method, the iPPP method, and DPALM, and the number of gradient steps performed by SP-LM and the SSG method. According to Figure~\ref{fig:figure_convex_experiment_step}, on all the instances studied, the iMELa method performs similarly to SP-LM and outperforms the iPPP method in all performance metrics. Moreover, DPALM exhibits performance comparable to iMELa and SP-LM on \textit{COMPAS}, while on \textit{a9a} and \textit{bank} it converges more slowly and attains worse infeasibility metrics. These observations are consistent with our theoretical finding that the iMELa method has an oracle complexity similar to SP-LM and better than those of the iPPP method and DPALM. In addition, the iMELa method and SP-LM reduce the infeasibility quickly and maintain a (nearly) feasible solution throughout the remaining iterations. Finally, we find that for all instances, the number of gradient steps performed by the APG method for solving \eqref{eq:imela_subprob} within the iMELa method and for solving \eqref{eq:dpalm_subprob} within DPALM is generally much smaller than that for solving \eqref{eq:ippp_subprob}, which is the main reason for the better empirical performance of the iMELa method and DPALM than the iPPP method.

\section{Conclusion}
We analyze the oracle complexity of the iMELa method for finding an $\epsilon$-KKT point in smooth non-convex functional constrained optimization. Each iteration consists of approximately solving a strongly convex subproblem, followed by updates to the Lagrange multipliers and the proximal center. Our results show that the oracle complexity of the iMELa method matches the best known result in the literature, up to a logarithmic factor, under a local error bound condition with exponent $d=1$, which is a strictly weaker assumption than the local LICQ condition required to achieve the best complexity guarantee.

\bibliography{references}

\begin{thebibliography}{55}
\providecommand{\natexlab}[1]{#1}
\providecommand{\url}[1]{\texttt{#1}}
\expandafter\ifx\csname urlstyle\endcsname\relax
  \providecommand{\doi}[1]{doi: #1}\else
  \providecommand{\doi}{doi: \begingroup \urlstyle{rm}\Url}\fi

\bibitem[Angwin et~al.(2016)Angwin, Larson, Mattu, and Kirchner]{angwin2016compas}
Julia Angwin, Jeff Larson, Surya Mattu, and Lauren Kirchner.
\newblock Machine bias.
\newblock \emph{ProPublica, May 23}, 2016.

\bibitem[Bertsekas(1999)]{bertsekas1999nonlinear}
Dimitri~P. Bertsekas.
\newblock \emph{Nonlinear Programming}.
\newblock Athena Scientific, 2nd edition, 1999.

\bibitem[Boob et~al.(2023)Boob, Deng, and Lan]{boob2023stochastic}
Digvijay Boob, Qi~Deng, and Guanghui Lan.
\newblock Stochastic first-order methods for convex and nonconvex functional constrained optimization.
\newblock \emph{Mathematical Programming}, 197\penalty0 (1):\penalty0 215--279, 2023.

\bibitem[Carmon et~al.(2020)Carmon, Duchi, Hinder, and Sidford]{carmon2020lower}
Yair Carmon, John~C. Duchi, Oliver Hinder, and Aaron Sidford.
\newblock Lower bounds for finding stationary points {I}.
\newblock \emph{Mathematical Programming}, 184\penalty0 (1--2):\penalty0 71--120, 2020.

\bibitem[Carmon et~al.(2021)Carmon, Duchi, Hinder, and Sidford]{carmon2021lower}
Yair Carmon, John~C. Duchi, Oliver Hinder, and Aaron Sidford.
\newblock Lower bounds for finding stationary points {II}: First-order methods.
\newblock \emph{Mathematical Programming}, 185\penalty0 (1--2):\penalty0 315--355, 2021.

\bibitem[Curtis et~al.(2015)Curtis, Jiang, and Robinson]{curtis2015adaptive}
Frank~E. Curtis, Hao Jiang, and Daniel~P. Robinson.
\newblock An adaptive augmented {L}agrangian method for large-scale constrained optimization.
\newblock \emph{Mathematical Programming}, 152\penalty0 (1--2):\penalty0 201--245, 2015.

\bibitem[Curtis et~al.(2016)Curtis, Gould, Jiang, and Robinson]{curtis2016adaptive}
Frank~E. Curtis, Nicholas~I.M. Gould, Hao Jiang, and Daniel~P. Robinson.
\newblock Adaptive augmented {L}agrangian methods: Algorithms and practical numerical experience.
\newblock \emph{Optimization Methods and Software}, 31\penalty0 (1):\penalty0 157--186, 2016.

\bibitem[Dahal et~al.(2023)Dahal, Liu, and Xu]{dahal2023damped}
Hari Dahal, Wei Liu, and Yangyang Xu.
\newblock Damped proximal augmented {L}agrangian method for weakly-convex problems with convex constraints.
\newblock \emph{arXiv preprint arXiv:2311.09065}, 2023.

\bibitem[Feldman et~al.(2015)Feldman, Friedler, Moeller, Scheidegger, and Venkatasubramanian]{feldman2015certifying}
Michael Feldman, Sorelle~A. Friedler, John Moeller, Carlos Scheidegger, and Suresh Venkatasubramanian.
\newblock Certifying and removing disparate impact.
\newblock In \emph{KDD'15: Proceedings of the 21th ACM SIGKDD International Conference on Knowledge Discovery and Data Mining}, pages 259--268, 2015.

\bibitem[Hajinezhad and Hong(2019)]{hajinezhad2019perturbed}
Davood Hajinezhad and Mingyi Hong.
\newblock Perturbed proximal primal-dual algorithm for nonconvex nonsmooth optimization.
\newblock \emph{Mathematical Programming}, 176\penalty0 (1--2):\penalty0 207--245, 2019.

\bibitem[He and Yuan(2010)]{he2010acceleration}
Bingsheng He and Xiaoming Yuan.
\newblock On the acceleration of augmented {L}agrangian method for linearly constrained optimization.
\newblock \emph{Optimization Online}, 2010.

\bibitem[Hestenes(1969)]{hestenes1969multiplier}
Magnus~R. Hestenes.
\newblock Multiplier and gradient methods.
\newblock \emph{Journal of Optimization Theory and Applications}, 4\penalty0 (5):\penalty0 303--320, 1969.

\bibitem[Hiriart-Urruty and Lemar{\'e}chal(2001)]{hiriart2001fundamentals}
Jean-Baptiste Hiriart-Urruty and Claude Lemar{\'e}chal.
\newblock \emph{Fundamentals of Convex Analysis}.
\newblock Grundlehren Text Editions. Springer, 2001.

\bibitem[Hong(2016)]{hong2016decomposing}
Mingyi Hong.
\newblock Decomposing linearly constrained nonconvex problems by a proximal primal dual approach: Algorithms, convergence, and applications.
\newblock \emph{arXiv preprint arXiv:1604.00543}, 2016.

\bibitem[Horn and Johnson(2012)]{horn2012matrix}
Roger~A. Horn and Charles~R. Johnson.
\newblock \emph{Matrix Analysis}.
\newblock Cambridge University Press, 2nd edition, 2012.

\bibitem[Huang and Lin(2023)]{huang2023oracle}
Yankun Huang and Qihang Lin.
\newblock Oracle complexity of single-loop switching subgradient methods for non-smooth weakly convex functional constrained optimization.
\newblock \emph{Advances in Neural Information Processing Systems}, 36:\penalty0 61327--61340, 2023.

\bibitem[Jia and Grimmer(2025)]{jia2025first}
Zhichao Jia and Benjamin Grimmer.
\newblock First-order methods for nonsmooth nonconvex functional constrained optimization with or without {S}later points.
\newblock \emph{SIAM Journal on Optimization}, 35\penalty0 (2):\penalty0 1300--1329, 2025.

\bibitem[Kohavi(1996)]{kohavi1996scaling}
Ron Kohavi.
\newblock Scaling up the accuracy of naive-{B}ayes classifiers: A decision-tree hybrid.
\newblock In \emph{KDD'96: Proceedings of the 2nd International Conference on Knowledge Discovery and Data Mining}, pages 202--207, 1996.

\bibitem[Kong et~al.(2019)Kong, Melo, and Monteiro]{kong2019complexity}
Weiwei Kong, Jefferson~G. Melo, and Renato~D.C. Monteiro.
\newblock Complexity of a quadratic penalty accelerated inexact proximal point method for solving linearly constrained nonconvex composite programs.
\newblock \emph{SIAM Journal on Optimization}, 29\penalty0 (4):\penalty0 2566--2593, 2019.

\bibitem[Kong et~al.(2020)Kong, Melo, and Monteiro]{kong2020efficient}
Weiwei Kong, Jefferson~G. Melo, and Renato~D.C. Monteiro.
\newblock An efficient adaptive accelerated inexact proximal point method for solving linearly constrained nonconvex composite problems.
\newblock \emph{Computational Optimization and Applications}, 76\penalty0 (2):\penalty0 305--346, 2020.

\bibitem[Kong et~al.(2023{\natexlab{a}})Kong, Melo, and Monteiro]{kong2023iteration-MathOR}
Weiwei Kong, Jefferson~G. Melo, and Renato~D.C. Monteiro.
\newblock Iteration complexity of a proximal augmented {L}agrangian method for solving nonconvex composite optimization problems with nonlinear convex constraints.
\newblock \emph{Mathematics of Operations Research}, 48\penalty0 (2):\penalty0 1066--1094, 2023{\natexlab{a}}.

\bibitem[Kong et~al.(2023{\natexlab{b}})Kong, Melo, and Monteiro]{kong2023iteration-SIOPT}
Weiwei Kong, Jefferson~G. Melo, and Renato~D.C. Monteiro.
\newblock Iteration complexity of an inner accelerated inexact proximal augmented {L}agrangian method based on the classical {L}agrangian function.
\newblock \emph{SIAM Journal on Optimization}, 33\penalty0 (1):\penalty0 181--210, 2023{\natexlab{b}}.

\bibitem[Lan and Monteiro(2016)]{lan2016iteration}
Guanghui Lan and Renato~D.C. Monteiro.
\newblock Iteration-complexity of first-order augmented {L}agrangian methods for convex programming.
\newblock \emph{Mathematical Programming}, 155\penalty0 (1):\penalty0 511--547, 2016.

\bibitem[Lewis and Pang(1998)]{lewis1998error}
Adrian~S. Lewis and Jong-Shi Pang.
\newblock Error bounds for convex inequality systems.
\newblock In Jean-Pierre Crouzeix, Juan-Enrique Martinez-Legaz, and Michel Volle, editors, \emph{Generalized Convexity, Generalized Monotonicity: Recent Results}, volume~27 of \emph{Nonconvex Optimization and Its Applications}, pages 75--110. Kluwer Academic Publishers, 1998.

\bibitem[Li(1997)]{li1997abadie}
Wu~Li.
\newblock Abadie's constraint qualification, metric regularity, and error bounds for differentiable convex inequalities.
\newblock \emph{SIAM Journal on Optimization}, 7\penalty0 (4):\penalty0 966--978, 1997.

\bibitem[Li and Xu(2021)]{li2021augmented}
Zichong Li and Yangyang Xu.
\newblock Augmented {L}agrangian-based first-order methods for convex-constrained programs with weakly convex objective.
\newblock \emph{INFORMS Journal on Optimization}, 3\penalty0 (4):\penalty0 373--397, 2021.

\bibitem[Li et~al.(2021)Li, Chen, Liu, Lu, and Xu]{li2021rate}
Zichong Li, Pin-Yu Chen, Sijia Liu, Songtao Lu, and Yangyang Xu.
\newblock Rate-improved inexact augmented {L}agrangian method for constrained nonconvex optimization.
\newblock In \emph{Proceedings of the 24th International Conference on Artificial Intelligence and Statistics}, volume 130, pages 2170--2178. PMLR, 2021.

\bibitem[Lin et~al.(2022)Lin, Ma, and Xu]{lin2022complexity}
Qihang Lin, Runchao Ma, and Yangyang Xu.
\newblock Complexity of an inexact proximal-point penalty method for constrained smooth non-convex optimization.
\newblock \emph{Computational Optimization and Applications}, 82\penalty0 (1):\penalty0 175--224, 2022.

\bibitem[Lin et~al.(2025)Lin, Soheili, Ma, and Nadarajah]{lin2025adaptive}
Qihang Lin, Negar Soheili, Runchao Ma, and Selvaprabu Nadarajah.
\newblock An adaptive parameter-free and projection-free restarting level set method for constrained convex optimization under the error bound condition.
\newblock \emph{Journal of Machine Learning Research}, 26\penalty0 (200):\penalty0 1--45, 2025.

\bibitem[Liu and Xu(2025)]{liu2025single}
Wei Liu and Yangyang Xu.
\newblock A single-loop {SPIDER}-type stochastic subgradient method for expectation-constrained nonconvex nonsmooth optimization.
\newblock \emph{arXiv preprint arXiv:2501.19214}, 2025.

\bibitem[Liu et~al.(2019)Liu, Liu, and Ma]{liu2019nonergodic}
Ya-Feng Liu, Xin Liu, and Shiqian Ma.
\newblock On the nonergodic convergence rate of an inexact augmented {L}agrangian framework for composite convex programming.
\newblock \emph{Mathematics of Operations Research}, 44\penalty0 (2):\penalty0 632--650, 2019.

\bibitem[Luo and Pang(1994)]{luo1994error}
Zhi-Quan Luo and Jong-Shi Pang.
\newblock Error bounds for analytic systems and their applications.
\newblock \emph{Mathematical Programming}, 67\penalty0 (1--3):\penalty0 1--28, 1994.

\bibitem[Ma et~al.(2020)Ma, Lin, and Yang]{ma2020quadratically}
Runchao Ma, Qihang Lin, and Tianbao Yang.
\newblock Quadratically regularized subgradient methods for weakly convex optimization with weakly convex constraints.
\newblock In \emph{Proceedings of the 37th International Conference on Machine Learning}, volume 119, pages 6554--6564. PMLR, 2020.

\bibitem[Melo et~al.(2020)Melo, Monteiro, and Wang]{melo2020iteration}
Jefferson~G. Melo, Renato~D.C. Monteiro, and Hairong Wang.
\newblock Iteration-complexity of an inexact proximal accelerated augmented {L}agrangian method for solving linearly constrained smooth nonconvex composite optimization problems.
\newblock \emph{arXiv preprint arXiv:2006.08048}, 2020.

\bibitem[Moro et~al.(2014)Moro, Cortez, and Rita]{moro2014data}
S{\'e}rgio Moro, Paulo Cortez, and Paulo Rita.
\newblock A data-driven approach to predict the success of bank telemarketing.
\newblock \emph{Decision Support Systems}, 62:\penalty0 22--31, 2014.

\bibitem[Nesterov(2013)]{nesterov2013gradient}
Yurii Nesterov.
\newblock Gradient methods for minimizing composite functions.
\newblock \emph{Mathematical Programming}, 140\penalty0 (1):\penalty0 125--161, 2013.

\bibitem[Nesterov(2018)]{nesterov2018lectures}
Yurii Nesterov.
\newblock \emph{Lectures on Convex Optimization}, volume 137 of \emph{Springer Optimization and Its Applications}.
\newblock Springer, 2nd edition, 2018.

\bibitem[Pang(1997)]{pang1997error}
Jong-Shi Pang.
\newblock Error bounds in mathematical programming.
\newblock \emph{Mathematical Programming}, 79\penalty0 (1--3):\penalty0 299--332, 1997.

\bibitem[Powell(1969)]{powell1969method}
Michael~J.D. Powell.
\newblock A method for nonlinear constraints in minimization problems.
\newblock In Roger Fletcher, editor, \emph{Optimization}, pages 283--298. Academic Press, 1969.

\bibitem[Pu et~al.(2024)Pu, Sun, and Zhang]{pu2024smoothed}
Wenqiang Pu, Kaizhao Sun, and Jiawei Zhang.
\newblock Smoothed proximal {L}agrangian method for nonlinear constrained programs.
\newblock \emph{arXiv preprint arXiv:2408.15047}, 2024.

\bibitem[Rockafellar(1970)]{rockafellar1970convex}
R.~Tyrrell Rockafellar.
\newblock \emph{Convex Analysis}, volume~28 of \emph{Princeton Mathematical Series}.
\newblock Princeton University Press, 1970.

\bibitem[Rockafellar(1973{\natexlab{a}})]{rockafellar1973dual}
R.~Tyrrell Rockafellar.
\newblock A dual approach to solving nonlinear programming problems by unconstrained optimization.
\newblock \emph{Mathematical Programming}, 5\penalty0 (1):\penalty0 354--373, 1973{\natexlab{a}}.

\bibitem[Rockafellar(1973{\natexlab{b}})]{rockafellar1973multiplier}
R.~Tyrrell Rockafellar.
\newblock The multiplier method of {H}estenes and {P}owell applied to convex programming.
\newblock \emph{Journal of Optimization Theory and Applications}, 12\penalty0 (6):\penalty0 555--562, 1973{\natexlab{b}}.

\bibitem[Rockafellar(1974)]{rockafellar1974augmented}
R.~Tyrrell Rockafellar.
\newblock Augmented {L}agrange multiplier functions and duality in nonconvex programming.
\newblock \emph{SIAM Journal on Control}, 12\penalty0 (2):\penalty0 268--285, 1974.

\bibitem[Rockafellar(1976)]{rockafellar1976augmented}
R.~Tyrrell Rockafellar.
\newblock Augmented {L}agrangians and applications of the proximal point algorithm in convex programming.
\newblock \emph{Mathematics of Operations Research}, 1\penalty0 (2):\penalty0 97--116, 1976.

\bibitem[Rockafellar and Wets(1998)]{rockafellar1998variational}
R.~Tyrrell Rockafellar and Roger J-B Wets.
\newblock \emph{Variational Analysis}, volume 317 of \emph{Grundlehren der mathematischen Wissenschaften}.
\newblock Springer, 1998.

\bibitem[Sahin et~al.(2019)Sahin, Eftekhari, Alacaoglu, Latorre, and Cevher]{sahin2019inexact}
Mehmet~Fatih Sahin, Armin Eftekhari, Ahmet Alacaoglu, Fabian Latorre, and Volkan Cevher.
\newblock An inexact augmented {L}agrangian framework for nonconvex optimization with nonlinear constraints.
\newblock \emph{Advances in Neural Information Processing Systems}, 32:\penalty0 13966--13978, 2019.

\bibitem[Sun and Sun(2024)]{sun2024dual}
Kaizhao Sun and Xu~Andy Sun.
\newblock Dual descent augmented {L}agrangian method and alternating direction method of multipliers.
\newblock \emph{SIAM Journal on Optimization}, 34\penalty0 (2):\penalty0 1679--1707, 2024.

\bibitem[Xu(2021{\natexlab{a}})]{xu2021first}
Yangyang Xu.
\newblock First-order methods for constrained convex programming based on linearized augmented {L}agrangian function.
\newblock \emph{INFORMS Journal on Optimization}, 3\penalty0 (1):\penalty0 89--117, 2021{\natexlab{a}}.

\bibitem[Xu(2021{\natexlab{b}})]{xu2021iteration}
Yangyang Xu.
\newblock Iteration complexity of inexact augmented {L}agrangian methods for constrained convex programming.
\newblock \emph{Mathematical Programming}, 185\penalty0 (1--2):\penalty0 199--244, 2021{\natexlab{b}}.

\bibitem[Yang and Lin(2018)]{yang2018rsg}
Tianbao Yang and Qihang Lin.
\newblock {RSG}: Beating subgradient method without smoothness and strong convexity.
\newblock \emph{Journal of Machine Learning Research}, 19\penalty0 (6):\penalty0 1--33, 2018.

\bibitem[Zeng et~al.(2022)Zeng, Yin, and Zhou]{zeng2022moreau}
Jinshan Zeng, Wotao Yin, and Ding-Xuan Zhou.
\newblock Moreau envelope augmented {L}agrangian method for nonconvex optimization with linear constraints.
\newblock \emph{Journal of Scientific Computing}, 91\penalty0 (2):\penalty0 1--36,{ no. }61, 2022.

\bibitem[Zhang and Luo(2020)]{zhang2020proximal}
Jiawei Zhang and Zhi-Quan Luo.
\newblock A proximal alternating direction method of multiplier for linearly constrained nonconvex minimization.
\newblock \emph{SIAM Journal on Optimization}, 30\penalty0 (3):\penalty0 2272--2302, 2020.

\bibitem[Zhang and Luo(2022)]{zhang2022global}
Jiawei Zhang and Zhi-Quan Luo.
\newblock A global dual error bound and its application to the analysis of linearly constrained nonconvex optimization.
\newblock \emph{SIAM Journal on Optimization}, 32\penalty0 (3):\penalty0 2319--2346, 2022.

\bibitem[Zhang et~al.(2022)Zhang, Pu, and Luo]{zhang2022iteration}
Jiawei Zhang, Wenqiang Pu, and Zhi-Quan Luo.
\newblock On the iteration complexity of smoothed proximal {ALM} for nonconvex optimization problem with convex constraints.
\newblock \emph{arXiv preprint arXiv:2207.06304}, 2022.

\end{thebibliography}

\appendix
\section{Proof of Proposition~\ref{thm:LICQ_implies_local_EBC}}
\label{sec:LICQ_implies_local_EBC}
In this section, we provide the detailed proof for Proposition~\ref{thm:LICQ_implies_local_EBC}, starting with a few technical lemmas.
\begin{lemma}
\label{lmm:bound_KKT_set}
The set $\X^*$ is compact.
\end{lemma}
\begin{proof}
By the structure of $\X$ and the Slater's condition in Assumption~\ref{assume:general}A and F, $\X^*$ is non-empty (see, e.g.~\citet[Sec. 3.3]{bertsekas1999nonlinear}). We then prove the lemma by two steps.

First, we show that, for any KKT point $\bx^*\in\X^*$, its any associated vector of Lagrangian multipliers $\blambda^*$ satisfies that $\|\blambda^*\|\leq B_{\blambda}:=\frac{B_f D_{\X}}{\min_{i\in[m]}[-g_i(\bx_{\text{feas}})]}$. By the KKT conditions at $(\bx^*,\blambda^*)$ in \eqref{eq:gco_general}, we have
\begin{align}
\label{eq:x*_KKT_stat_comple_slack}
    \textstyle  \bxi^*:= -\nabla f(\bx^*) - \sum_{i=1}^m \lambda_i^* \nabla g_i(\bx^*) \in \N_\X(\bx^*)\text{ and }
    \lambda_i^* g_i(\bx^*) = 0,\;i\in [m].
\end{align} 
Hence, by the convexity of $g_i$'s and $\blambda^*\geq\bzero$, it follows that
\begin{align*}
    \|\blambda^*\|_1\cdot \min_{i\in[m]}[-g_i(\bx_{\text{feas}})]
    \leq&
    -\textstyle\sum_{i=1}^m \lambda_i^* g_i(\bx_{\text{feas}})
    \leq-\textstyle\sum_{i=1}^m \lambda_i^* \left(g_i(\bx^*) + \langle \nabla g_i(\bx^*), \bx_{\text{feas}}-\bx^* \rangle\right)\\
    =&~
    \big\langle \bxi^* + \nabla f(\bx^*), 
    \bx_{\text{feas}}-\bx^*\big\rangle
    \leq
    \big\langle\nabla f(\bx^*), \bx_{\text{feas}}-\bx^*\big\rangle
    \leq B_f D_{\X},
\end{align*}
where $\bx_{\text{feas}}$ is given in Assumption~\ref{assume:general}F, the equality holds by~\eqref{eq:x*_KKT_stat_comple_slack}, the third inequality holds by $\bxi^*\in \N_\X(\bx^*)$, and the last one is from Assumption~\ref{assume:general}A, D, and F. Thus, we obtain $\|\blambda^*\|\leq \|\blambda^*\|_1\leq B_{\blambda}$,
which proves the claimed result.

Second, by Assumption~\ref{assume:general}A, $\X^*$ is bounded. Hence, it suffices to show that $\X^*$ is closed. Suppose there is a sequence of KKT points and their associated vector of Lagrangian multipliers, denote by $\{(\bar{\bx}^{(k)},\bar{\blambda}^{(k)})\}_{k\geq1}\subset\X^*\times\mathbb{R}_+^m$, such that $\bar{\bx}^{(k)}\rightarrow\bar\bx^*$ as $k\rightarrow \infty$ for some $\bar\bx^*\in\X$. By definition, there exists $\bar{\bxi}^{(k)}\in \N_\X(\bar{\bx}^{(k)})$ such that  
\begin{align*}
    \nabla f(\bar{\bx}^{(k)})+\textstyle\sum_{i=1}^m \bar{\lambda}_i^{(k)}\nabla g_i(\bar{\bx}^{(k)})+\bar{\bxi}^{(k)}=\bzero,
    ~g_i(\bar{\bx}^{(k)})\leq0, ~\bar{\lambda}_i^{(k)} g_i(\bar{\bx}^{(k)})=0, ~\forall\, i\in[m].
\end{align*}
Since $\|\bar{\blambda}^{(k)}\|\leq B_{\blambda}$, passing to a subsequence if necessary, it holds that $\lim_{k\rightarrow \infty}\bar{\blambda}^{(k)}=\bar{\blambda}^*$ for some $\bar{\blambda}^*\geq\bzero$. As $k\rightarrow\infty$, we conclude from the continuity of $\nabla f(\cdot)$ and $\nabla\bg(\cdot)$ and the outer semi-continuity of $\N_\X(\cdot)$~\citep[Proposition 6.6]{rockafellar1998variational} that $\bar{\bxi}^{(k)}\rightarrow\bar{\bxi}^*$ for some $\bar{\bxi}^*\in \N_\X(\bar{\bx}^*)$ and
\begin{align*}
    \nabla f(\bar{\bx}^*)+\textstyle\sum_{i=1}^m \bar{\lambda}_i^*\nabla g_i(\bar{\bx}^*)+\bar{\bxi}^*=\bzero,
    ~g_i(\bar{\bx}^*)\leq0, ~\bar{\lambda}_i^* g_i(\bar{\bx}^*)=0, ~\forall\,i\in[m],
\end{align*}
meaning that $\bar{\bx}^*\in\X^*$, which completes the proof. 
\end{proof}

\begin{lemma}
\label{lmm:nearKKT_activeset}
For any $\vartheta>0$, there exists $r(\vartheta)>0$ such that, if $\bx$ is feasible to \eqref{eq:gco_general} and $\textup{dist}(\bx,\X^*)\leq r(\vartheta)$, then there must exist $\bx^*\in\X^*$ such that $\|\bx-\bx^*\|\leq\vartheta$ and $J_g(\bx)\subset J_g(\bx^*)\text{ and }J_A(\bx)\subset J_A(\bx^*)$.
\end{lemma}
\begin{proof}
We prove this by contradiction. Suppose that the claim does not hold. Then there exists $\vartheta>0$ and a sequence $\{\bx^{(k)}\}_{k\geq1}$ such that $\bx^{(k)}$ is feasible to \eqref{eq:gco_general} and $\text{dist}(\bx^{(k)},\X^*)\leq r_k:=\frac{\vartheta}{k}\leq\vartheta$ but $J_g(\bx^{(k)})\not\subset J_g(\bx^*)\text{ or }J_A(\bx^{(k)})\not\subset J_A(\bx^*)$ for any $k$ and any $\bx^*\in\X^*$
satisfying $\|\bx^{(k)}-\bx^*\|\leq\vartheta$. 

Due to the compactness of $\X$ and $\X^*$ (see Lemma~\ref{lmm:bound_KKT_set}),  by passing to a subsequence if necessary, there exists some $\bar{\bx}^*\in\X^*$ such that  
$\bx^{(k)}\rightarrow\bar{\bx}^*$ as $k\rightarrow\infty$. By the continuity of $\bg(\cdot)$ and $\bA(\cdot)-\bb$, it holds that for any $i\in[m]$, if $g_i(\bar{\bx}^*)<0$, then $g_i(\bx^{(k)})<0$ as well when $k$ is large enough, and for any $j\in[l]$, if $[\bA\bar{\bx}^*-\bb]_j<0$, then $[\bA\bx^{(k)}-\bb]_j<0$ as well when $k$ is large enough. That means $J_g(\bx^{(k)})\subset J_g(\bar\bx^*)$ and $ J_A(\bx^{(k)})\subset J_A(\bar\bx^*)$ for a large enough $k$. This contradicts to the assumption.  
Therefore, the conclusion holds.
\end{proof}
\begin{lemma}
\label{lmm:nearKKT_LICQ}
Suppose the constraints of \eqref{eq:gco_general} satisfy the local LICQ in \eqref{eq:LICQ_sigmamin}.
Then there exists $\delta'>0$ such that, if $\bx$ is feasible to \eqref{eq:gco_general} and $\|\bx-\bx^*\|\leq\delta'$ for some $\bx^*\in\X^*$, then $\sigma_{\min}\big(\big[\nabla \bg_{I_g}(\bx), \bA_{I_A}^\top\big]\big)\geq\frac{\zeta}{2}$ for any $I_g\subset J_g(\bx^*)$ and $
I_A\subset J_A(\bx^*)$.
\end{lemma}
\begin{proof}
Fix $\bx^*\in\X^*$ and any $I_g\subset J_g(\bx^*)$, $I_A\subset J_A(\bx^*)$. By the continuity of $\sigma_{\min}(\cdot)$ (see e.g.~\citet[Corollary 7.3.5 (a)]{horn2012matrix}), for any $\bx$ we have
\begin{align*}
    \sigma_{\min}\big(\big[ \nabla \bg_{I_g}(\bx), \bA_{I_A}^\top\big]\big)
    \geq&~ \sigma_{\min}\big(\big[ \nabla \bg_{I_g}(\bx^*), \bA_{I_A}^\top\big]\big)
    -\big\|\big[ \nabla \bg_{I_g}(\bx),\bA_{I_A}^\top\big]- \big[\nabla \bg_{I_g}(\bx^*),\bA_{I_A}^\top \big]\big\|\\
    =&~ \sigma_{\min}\big(\big[ \nabla \bg_{I_g}(\bx^*), \bA_{I_A}^\top\big]\big)
    -\big\|\big[ \nabla \bg_{I_g}(\bx)- \nabla \bg_{I_g}(\bx^*),\; \bO\big]\big\|\\
    \geq&~ \sigma_{\min}\big(\big[ \nabla \bg_{I_g}(\bx^*), \bA_{I_A}^\top\big]\big)
    -\big\|\nabla \bg_{I_g}(\bx)- \nabla \bg_{I_g}(\bx^*)\big\|.
\end{align*}
Since the functions $g_i$, $i\in[m]$, are $L$-smooth by Assumption~\ref{assume:general}C, we have
$$
\big\|\nabla \bg_{I_g}(\bx)- \nabla \bg_{I_g}(\bx^*)\big\|
\leq L\sqrt{|I_g|}\,\|\bx-\bx^*\|
\leq L\sqrt{m}\,\|\bx-\bx^*\|.
$$
Define $L_H:=L\sqrt{m}$ and set
$$
\delta':=\frac{\zeta}{2L_H}=\frac{\zeta}{2L\sqrt{m}}.
$$
Then for all feasible $\bx$ satisfying $\|\bx-\bx^*\|\leq\delta'$, we obtain
\begin{align*}
    \sigma_{\min}\big(\big[\nabla \bg_{I_g}(\bx), \bA_{I_A}^\top\big]\big)
    \geq&~
    \sigma_{\min}\big(\big[ \nabla \bg_{I_g}(\bx^*), \bA_{I_A}^\top\big]\big)
    -L_H\|\bx-\bx^*\| \\
    \geq&~
    \zeta-L_H\delta'
    =\zeta-\tfrac{\zeta}{2}
    =\tfrac{\zeta}{2}.
\end{align*}
\end{proof}
\begin{proof}[Proof of Proposition~\ref{thm:LICQ_implies_local_EBC}]
Let $\delta=\min\left\{\frac{r(\delta')}{2},\frac{\zeta}{2L\sqrt{m}}\right\}$, where $r(\cdot)$ is from Lemma \ref{lmm:nearKKT_activeset} and $\delta'$ is from Lemma~\ref{lmm:nearKKT_LICQ}. Suppose that $\bx\in\X$ and $\|\bx-\bx^*\|\leq \delta$ for some $\bx^*\in\X^*$, we then need to show that there exists $\gamma>0$ such that \eqref{eq:local_EBC_exponent_d} holds with $d=1$ for any $I_g\subset J_g(\bx^*)$ and $I_A\subset J_A(\bx^*)$.
For such $I_g$ and $I_A$,
consider the minimization problem 
\begin{align}
\label{eq:projbarx}
    \bar{\bx}=\argmin_{\bz} \frac{1}{2}\|\bz-\bx\|^2
    \text{ s.t.}&\begin{array}{l}
    g_i(\bz)=0,\,i\in I_g,\,g_i(\bz)\leq 0,\,i\in[m]\backslash I_g,\\[0.1ex]
    [\bA\bz-\bb]_j=0,\,j\in I_A,\,[\bA\bz- \bb]_j\leq0,\,j\in[l]\backslash I_A.
    \end{array}
\end{align}
Note that $\bar{\bx}$ is well-defined because the feasible set of \eqref{eq:projbarx} is compact and non-empty as it at least contains $\bx^*$. As a result, $\text{dist} (\bx,\mathcal{S}(I_g, I_A))=\|\bx-\bar{\bx}\|\leq \|\bx-\bx^*\|\leq \delta$ and thus $\text{dist}(\bar{\bx},\X^*)\leq\|\bar{\bx}-\bx^*\|\leq 2\delta \leq r(\delta')$. Since $\bar\bx$ is also feasible to \eqref{eq:gco_general}, by invoking Lemma~\ref{lmm:nearKKT_activeset}, there exists $\bar\bx^*\in\X^*$ such that $\|\bar\bx-\bar\bx^*\|\leq\delta'$ and $J_g(\bar\bx)\subset J_g(\bar\bx^*)$ and $J_A(\bar\bx)\subset J_A(\bar\bx^*)$. By Lemma~\ref{lmm:nearKKT_LICQ} at $\bar\bx$ and $\delta'$, we have $\sigma_{\min}\big(\big[\nabla \bg_{J_g(\bar\bx)}(\bar\bx), \bA_{J_A(\bar\bx)}^\top\big]\big)\geq\frac{\zeta}{2}$. This means \eqref{eq:projbarx} satisfies the LICQ at $\bar\bx$. Hence, there exist a vector of Lagrangian multipliers 
$\bar\blambda\in\mathbb{R}^{m}$ and $\bar\bw\in\mathbb{R}^{l}$ such that 
\begin{align}
\label{eq:KKTprojection_stat}
    &\bar{\bx}-\bx+\textstyle\sum_{i=1}^m \bar\lambda_i \cdot\nabla g_i(\bar{\bx})
    +\textstyle\sum_{j=1}^l  \bar{w}_j\cdot[\bA]_j^\top=\bzero,\\
    \label{eq:KKTprojection_feas}
    &\bar\lambda_i \geq 0,\; i\in[m]\backslash I_g,\quad\bar{w}_j \geq 0,\;j\in[l]\backslash I_{A},\\
    \label{eq:KKTprojection_comple_slack}
    &\bar\lambda_i\cdot g_i(\bar{\bx}) = 0,\; i\in[m]\backslash I_g,\quad\bar{w}_j\cdot [\bA\bar{\bx}-\bb]_j = 0,\;j\in[l]\backslash I_{A}.
\end{align}
Notice $I_g\subset J_g(\bar{\bx})$ and $I_A\subset J_A(\bar{\bx})$.
By \eqref{eq:KKTprojection_comple_slack}, we must have $\bar\lambda_i=0$ for $i\in[m]\backslash J_g(\bar{\bx})$ and $\bar{w}_j=0$ for $j\in[l]\backslash J_A(\bar{\bx})$. Thus, \eqref{eq:KKTprojection_stat} and $\sigma_{\min}([\nabla \bg_{J_g(\bar{\bx})}(\bar{\bx}), \bA_{ J_A(\bar{\bx})}^\top])\geq \textstyle \frac{\zeta}{2}$ imply 
\begin{align}
\label{eq:KKTprojection_zeta}
    \|\bar{\bx}-\bx\|=\left\|\textstyle\sum_{i\in J_g(\bar{\bx}) } \bar\lambda_i \cdot \nabla g_i(\bar{\bx})
    +\textstyle\sum_{j\in J_A(\bar{\bx})} \bar{w}_j\cdot [\bA]_j^\top\right\|\geq \textstyle \frac{\zeta}{2}\sqrt{\rule{0pt}{2ex}\|\bar\blambda\|^2+\|\bar{\bw}\|^2}.
\end{align}
Moreover, for $i\in I_g$, by the $L$-smoothness of $g_i$ and the fact that $g_i(\bar{\bx})=0$, we have
\begin{align}
    \nonumber
    \bar\lambda_i\cdot \left\langle\nabla g_i(\bar{\bx}),\bx-\bar{\bx}\right\rangle
    \leq&~|\bar\lambda_i|\cdot|g_i(\bx)-g_i(\bar{\bx})|+\frac{|\bar\lambda_i|L}{2}\|\bx-\bar{\bx}\|^2\\
    \label{eq:lambda_bar_in_I_g}
    =&~ |\bar\lambda_i|\cdot|g_i(\bx)|+\frac{|\bar\lambda_i|L}{2}\|\bx-\bar{\bx}\|^2.
\end{align}
For $i\in[m]\backslash I_g$, using \eqref{eq:KKTprojection_feas}, \eqref{eq:KKTprojection_comple_slack} and the convexity of $g_i$, we have 
\begin{align}
\label{eq:lambda_bar_notin_I_g}
    \bar\lambda_i\cdot\left\langle\nabla g_i(\bar{\bx}),\bx-\bar{\bx}\right\rangle\leq \bar\lambda_i\cdot\left(g_i(\bx)-g_i(\bar{\bx})\right)
    \leq |\bar\lambda_i|\cdot[g_i(\bx)]_+.
\end{align}
Following similar arguments, we have
\begin{align}
\label{eq:w_bar_in_[l]}
    \bar{w}_j\cdot [\bA(\bx-\bar{\bx}) ]_j
    \leq\begin{cases}
        |\bar{w}_j|\cdot|[\bA\bx-\bb ]_j| & \text{ for }j\in I_A\\
        |\bar{w}_j|\cdot [[\bA\bx-\bb ]_j]_+ & \text{ for }j\in[l]\backslash I_{A}.
    \end{cases}
\end{align}

Taking the inner product of the left-hand side of \eqref{eq:KKTprojection_stat} and $\bar{\bx}-\bx$ and applying \eqref{eq:lambda_bar_in_I_g}, \eqref{eq:lambda_bar_notin_I_g}, and \eqref{eq:w_bar_in_[l]} lead to
\begin{small}
\begin{align*}
    \|\bar{\bx}-\bx\|^2
    \leq&\left(\begin{array}{l}
    \sum_{i\in I_g} |\bar\lambda_i|\cdot|g_i(\bx)|
    +\sum_{i\in[m]\backslash I_g} |\bar\lambda_i|\cdot[g_i(\bx)]_+\\[1ex]
    +\sum_{j\in I_{A}} |\bar{w}_j|\cdot\left|[\bA\bx-\bb]_j\right|
    +\sum_{j\in[l]\backslash I_{A}} |\bar{w}_j|\cdot\left[[\bA\bx-\bb]_j\right]_+
    \end{array}
    \right)+\sum_{i\in I_g}\frac{|\bar\lambda_i|L}{2}\|\bx-\bar{\bx}\|^2\\
    \leq&\left(\begin{array}{l}
    \sum_{i\in I_g} |g_i(\bx)|^2
    +\sum_{i\in[m]\backslash I_g} \left([g_i(\bx)]_+\right)^2\\[0.5ex]
    +\sum_{j\in I_{A}}\left|[\bA\bx-\bb]_j\right|^2
    +\sum_{j\in[l]\backslash I_{A}} \left(\left[[\bA\bx-\bb]_j\right]_+\right)^2
    \end{array}
    \right)^{1/2}\cdot \sqrt{\|\bar\blambda\|^2+\|\bar{\bw}\|^2}\\
    &+\sqrt{m} \sqrt{\|\bar\blambda\|^2+\|\bar{\bw}\|^2}\cdot\frac{L}{2}\|\bx-\bar{\bx}\|^2.
\end{align*}
\end{small}Recall \eqref{eq:KKTprojection_zeta} and the fact that $\|\bx-\bar{\bx}\|\leq \delta\leq \frac{\zeta}{2L\sqrt{m}}$. The inequality above implies
\begin{small}
\begin{align*}
    \|\bar{\bx}-\bx\|^2
    \leq&\left(\begin{array}{l}
    \sum_{i\in I_g} |g_i(\bx)|^2
    +\sum_{i\in[m]\backslash I_g} \left([g_i(\bx)]_+\right)^2\\[0.5ex]
    +\sum_{j\in I_{A}}\left|[\bA\bx-\bb]_j\right|^2
    +\sum_{j\in[l]\backslash I_{A}} \left(\left[[\bA\bx-\bb]_j\right]_+\right)^2
    \end{array}
    \right)^{1/2}\cdot \frac{2\|\bar{\bx}-\bx\|}{\zeta}\\
    &+\sqrt{m}\|\bar{\bx}-\bx\|\cdot\frac{L}{\zeta}\|\bx-\bar{\bx}\|^2\\
    \leq&\left(\begin{array}{l}
    \sum_{i\in I_g} |g_i(\bx)|^2
    +\sum_{i\in[m]\backslash I_g} \left([g_i(\bx)]_+\right)^2\\[0.5ex]
    +\sum_{j\in I_{A}}\left|[\bA\bx-\bb]_j\right|^2
    +\sum_{j\in[l]\backslash I_{A}} \left(\left[[\bA\bx-\bb]_j\right]_+\right)^2
    \end{array}
    \right)^{1/2}\cdot \frac{2\|\bar{\bx}-\bx\|}{\zeta}+\frac{1}{2}\|\bx-\bar{\bx}\|^2,
\end{align*}
\end{small}which implies \eqref{eq:local_EBC_exponent_d} with $d=1$ and $\gamma=16/\zeta^2$.

On the contrary, we can next show that Assumption~\ref{assume:local_EBC_exponent_d} does not imply the local LICQ in \eqref{eq:LICQ_sigmamin}, indicating that Assumption~\ref{assume:local_EBC_exponent_d} is weaker than the local LICQ. 

Consider an instance with an objective $f(\bx)=x_1^2+x_2^2$ and $\X:=\{\bx\in\mathbb{R}^2:\ -2\leq x_1\leq 2,\ -2\leq x_2\leq 2\}$, and introduce two nonlinear convex inequality constraints $g_1(\bx):=x_1$, $g_2(\bx):=x_1+x_1^2$. Clearly, both $g_1$ and $g_2$ are smooth and convex on $\mathbb{R}^2$, and the system admits a Slater point, e.g., $\bx_{\text{feas}}:=(-1/2,0)^\top$ satisfying $g_1(\bx_{\text{feas}})<0$ and $g_2(\bx_{\text{feas}})<0$. It is immediate that $\bx^*=\bzero$ is the unique KKT point of this instance. Moreover, we have $J_g(\bx^*)=\{1,2\}$ and $J_A(\bx^*)=\emptyset$. However, the local LICQ fails at $\bx^*$. Indeed, both constraints are active at $\bx^*$, while their gradients are linearly dependent:
$$
\nabla g_1(\bx^*)=(1,0)^\top,\quad 
\nabla g_2(\bx^*)=(1+2x_1,0)^\top\big|_{\bx=\bx^*}=(1,0)^\top,
$$
so the active gradients cannot be linearly independent. 

On the other hand, \eqref{eq:local_EBC_exponent_d} holds around $\bx^*$ with exponent $d=1$ and an explicit constant $\gamma=4$. Fix $\delta:=1/2$. For any $\bx\in\X\cap\mathbb{B}(\bx^*,\delta)$, we have $x_1\in[-1/2,1/2]$. Now take any $I_g\subseteq J_g(\bx^*)=\{1,2\}$ and $I_A\subseteq J_A(\bx^*)=\emptyset$.

Case 1: $I_g=\emptyset$.
Then, by definition, $\mathcal{S}(\emptyset,\emptyset)=\{\bx\in\X:\ g_1(\bx)\leq 0,\ g_2(\bx)\leq 0\}$. For any $\bx\in\X\cap\mathbb{B}(\bx^*,\delta)$, if $x_1\leq 0$ then $x_1\in[-1/2,0]$ and hence $1+x_1\geq 1/2$, which yields $g_2(\bx)=x_1+x_1^2=x_1(1+x_1)\leq 0$. Hence, for any $\bx\in\X\cap\mathbb{B}(\bx^*,\delta)$, we have $g_1(\bx)\leq 0 \Rightarrow g_2(\bx)\leq 0$, and thus $\{\bx\in\X:\ x_1\leq 0\}\cap\mathbb{B}(\bx^*,\delta)\subseteq \mathcal{S}(\emptyset,\emptyset)$. Conversely, $\bx\in \mathcal{S}(\emptyset,\emptyset)$ implies $g_1(\bx)\leq 0$, i.e., $x_1\leq 0$. Therefore, for all $\bx\in\X\cap\mathbb{B}(\bx^*,\delta)$,
$$
\text{dist}(\bx,\mathcal{S}(\emptyset,\emptyset))^2
=([g_1(\bx)]_+)^2
\leq 4\big(([g_1(\bx)]_+)^2+([g_2(\bx)]_+)^2\big).
$$
Case 2: $I_g\neq\emptyset$. We show that for all $\bx\in\X\cap\mathbb{B}(\bx^*,\delta)$, $\text{dist}(\bx,\mathcal{S}(I_g,\emptyset))=|x_1|$.

If $1\in I_g$, then every $\bz\in\mathcal{S}(I_g,\emptyset)$ satisfies $g_1(\bz)=0$, and hence $z_1=0$.
Therefore,
$$
\text{dist}(\bx,\mathcal{S}(I_g,\emptyset))
\geq
\text{dist}(\bx,\{\bz\in\X:\ z_1=0\})
=
|x_1|.
$$
On the other hand, the point $\hat{\bx}:=(0,x_2)^\top$ belongs to $\mathcal{S}(I_g,\emptyset)$, which yields
$$
\text{dist}(\bx,\mathcal{S}(I_g,\emptyset))
\leq
\|\bx-\hat{\bx}\|
=
|x_1|.
$$
Combining the two inequalities gives
$\text{dist}(\bx,\mathcal{S}(I_g,\emptyset))=|x_1|$.

If $1\notin I_g$, then $2\in I_g$ and any $\bz\in\mathcal{S}(I_g,\emptyset)$ must satisfy
$g_2(\bz)=0$, i.e., $z_1\in\{0,-1\}$.
Since $x_1\in[-1/2,1/2]$ on $\mathbb{B}(\bx^*,\delta)$, we have
$|x_1-z_1|\geq |x_1|$ for all $z_1\in\{0,-1\}$, and hence
$$
\text{dist}(\bx,\mathcal{S}(I_g,\emptyset))
\geq
|x_1|.
$$
Moreover, the point $\hat{\bx}:=(0,x_2)^\top$ satisfies $g_2(\hat{\bx})=0$ and $g_1(\hat{\bx})=0$ and thus
$\hat{\bx}\in\mathcal{S}(I_g,\emptyset)$, which implies
$$
\text{dist}(\bx,\mathcal{S}(I_g,\emptyset))
\leq
\|\bx-\hat{\bx}\|
=
|x_1|.
$$
Combining the upper and lower bounds yields
$\text{dist}(\bx,\mathcal{S}(I_g,\emptyset))=|x_1|$.

With the distance identity established, we verify \eqref{eq:local_EBC_exponent_d} with $\gamma=4$.
If $x_1\geq 0$, then $[g_1(\bx)]_+=x_1$ and thus
$$
\text{dist}(\bx,\mathcal{S}(I_g,\emptyset))^2=|x_1|^2=([g_1(\bx)]_+)^2
\leq 4\big(|g_2(\bx)|^2+([g_1(\bx)]_+)^2\big).
$$
If $x_1\leq 0$, then $[g_1(\bx)]_+=0$ and $x_1\in[-1/2,0]$, so $1+x_1\geq 1/2$ and hence
$$
|g_2(\bx)|=|x_1+x_1^2|=|x_1|\,(1+x_1)\geq \frac12|x_1|.
$$
Thus,
$$
\text{dist}(\bx,\mathcal{S}(I_g,\emptyset))^2=|x_1|^2
\leq 4|g_2(\bx)|^2
\leq 4\big(|g_2(\bx)|^2+([g_1(\bx)]_+)^2\big).
$$

Combining Case 1 and Case 2, we conclude that for all $\bx\in\X\cap\mathbb{B}(\bx^*,\delta)$, all $I_g\subseteq J_g(\bx^*)$, and all $I_A\subseteq J_A(\bx^*)$ (with $J_A(\bx^*)=\emptyset$), \eqref{eq:local_EBC_exponent_d} holds with $d=1$ and $\gamma=4$. This example shows that \eqref{eq:local_EBC_exponent_d} with $d=1$ may hold even when the local LICQ fails at $\bx^*$. This completes the proof.

%
\end{proof}

\section{Proof of Technical Results}
In this section, we present the proofs of the technical results for the complexity analysis. 

\subsection{Technical Results for Subproblem Accuracy}
\label{sec:subprob_optim}
The following lemma shows that, when the inexact solution $\bx^{(t+1)}$ of the 
subproblem $\min_{\bx\in\X}\mathcal{L}_p(\bx,\bz^{(t)},\blambda^{(t+1)})$ satisfies~\eqref{eq:subprob_optim_imela}, it also produces a small objective gap for this subproblem.
\begin{lemma}
\label{lmm:subprob_descent}
Suppose $\bx^{(t+1)}\in\X$ satisfies~\eqref{eq:subprob_optim_imela}.
Then it holds that
\begin{equation}
\label{eq:subprob_optimy_imela_corollary}
    \mathcal{L}_p (\bx^{(t+1)},\bz^{(t)},\blambda^{(t+1)})
    - \mathcal{L}_p (\tilde{\bx}^{(t+1)},\bz^{(t)},\blambda^{(t+1)})
    \leq \frac{\epsilon_t^2}{2(p-L)}.
\end{equation}
\end{lemma}
\begin{proof}
Condition~\eqref{eq:subprob_optim_imela} is equivalent to the existence of $\bv^{(t+1)}\in\mathbb{R}^n$ and $\bxi^{(t+1)}\in\N_\X(\bx^{(t+1)})$ such that $\|\bv^{(t+1)}\|\leq\epsilon_t$ and
\begin{equation}
\label{eq:subprob_optim_imela_equiv} 
   \nabla_\bx\mathcal{L}_p (\bx^{(t+1)},\bz^{(t)},\blambda^{(t+1)})
   =-\,\bxi^{(t+1)}-\bv^{(t+1)}.
\end{equation}
By the definition of the normal cone, this implies that for any $\bx\in\X$,
\begin{equation}
\label{eq:subprob_optim_imela_xi}
    \langle \nabla_\bx\mathcal{L}_p (\bx^{(t+1)},\bz^{(t)}, \blambda^{(t+1)})+\bv^{(t+1)},\bx-\bx^{(t+1)} \rangle
    = \langle \bxi^{(t+1)},\bx^{(t+1)}-\bx \rangle\geq 0.
\end{equation}
In addition, 
since $\mathcal{L}_p(\cdot,\bz^{(t)},\blambda^{(t+1)})$ is $(p-L)$-strongly convex in $\bx$, we have
\begin{align*}
    &~\mathcal{L}_p (\bx^{(t+1)},\bz^{(t)},\blambda^{(t+1)}) - \mathcal{L}_p (\tilde{\bx}^{(t+1)},\bz^{(t)},\blambda^{(t+1)})\\
    \leq&-\langle\nabla_\bx\mathcal{L}_p (\bx^{(t+1)},\bz^{(t)},\blambda^{(t+1)}), \tilde{\bx}^{(t+1)}-\bx^{(t+1)}\rangle \textstyle-\frac{p-L}{2}\|\tilde{\bx}^{(t+1)}-\bx^{(t+1)}\|^2\\
    \leq&
    ~\langle\bv^{(t+1)}, \tilde{\bx}^{(t+1)}-\bx^{(t+1)}\rangle 
    -\frac{p-L}{2}\|\tilde{\bx}^{(t+1)}-\bx^{(t+1)}\|^2\\ 
    \leq&  
    ~\frac{1}{2(p-L)}\|\bv^{(t+1)}\|^2
    \leq 
    \frac{\epsilon_t^2}{2(p-L)},
\end{align*}
where 
the second inequality follows from~\eqref{eq:subprob_optim_imela_xi}
with $\bx=\tilde{\bx}^{(t+1)}$, and the last inequality follows from Young’s inequality. The final result yields \eqref{eq:subprob_optimy_imela_corollary}.
\end{proof}

\subsection{Proof of Lemma~\ref{lmm:bound_lambda}}
\label{sec:bound_lambda_lm}
\begin{proof}[Proof of Lemma~\ref{lmm:bound_lambda}]
Let $\bv^{(t+1)}$ be the vector from \eqref{eq:subprob_optim_imela_equiv}. Recall that $\bx_{\text{feas}}\in\X$ is from Assumption~\ref{assume:general}F. By \eqref{eq:subprob_optim_imela_xi} at $\bx=\bx_{\text{feas}}$, we have
\begin{align}
    \label{eq:subroutine_epsilon_lm}
    \langle\nabla_\bx \mathcal{L}_p (\bx^{(t+1)},\bz^{(t)},\blambda^{(t+1)})+\bv^{(t+1)},\bx^{(t+1)}-\bx_{\text{feas}}\rangle\leq 0.
\end{align}  
Using the convexity of $g_i$'s and $\blambda^{(t+1)}\geq\bzero$, we have
\begin{align}
    \nonumber
    \langle\blambda^{(t+1)}, \bg(\bx^{(t+1)})- \bg(\bx_{\text{feas}})\rangle
    \leq&~\textstyle\sum_{i=1}^m \langle\lambda_i^{(t+1)}\nabla g_i(\bx^{(t+1)}),\bx^{(t+1)}-\bx_{\text{feas}}\rangle\\   
    \nonumber
    =&~\langle \nabla_\bx\mathcal{L}_p(\bx^{(t+1)},\bz^{(t)},\blambda^{(t+1)})+\bv^{(t+1)},\bx^{(t+1)}-\bx_{\text{feas}}\rangle\\   
    \nonumber
    &+\langle \nabla f(\bx^{(t+1)})+p(\bx^{(t+1)}-\bz^{(t)})+\bv^{(t+1)},\bx_{\text{feas}}-\bx^{(t+1)}\rangle\\   
    \nonumber
    \leq&~\langle \nabla f(\bx^{(t+1)})+p(\bx^{(t+1)}-\bz^{(t)})+\bv^{(t+1)},\bx_{\text{feas}}-\bx^{(t+1)}\rangle\\
    \label{eq:lambda_result_C_lambda}
    \leq&~(B_f+p D_{\X}+1)\cdot D_{\X}=C_{\blambda},
\end{align}
where the second inequality is by \eqref{eq:subroutine_epsilon_lm} and the last one is from Assumption~\ref{assume:general}A and D, and the fact that $\|\bv^{(t+1)}\|\leq\epsilon_t<1$. Since $\blambda^{(0)}=\bzero$ in Algorithm~\ref{alg:imela}, \eqref{eq:lambda_result_C_lambda} also holds for $t=-1$. So, $\langle\blambda^{(t)}, \bg(\bx^{(t)})- \bg(\bx_{\text{feas}})\rangle\leq C_{\blambda},~\forall\,t\geq0$, or equivalently, 
\begin{align}
\label{eq:lambda_result_C_lambda_equiv}
   \langle\blambda^{(t)}, \bg(\bx^{(t)})\rangle\leq C_{\blambda}-\langle\blambda^{(t)}, - \bg(\bx_{\text{feas}})\rangle,~\forall\,t\geq0.
\end{align}

We then prove the result in~\eqref{eq:bound_lambda} by induction. It clearly holds for $t=0$ since $\blambda^{(0)}=\bzero$. 
Suppose that the inequality in \eqref{eq:bound_lambda} holds 
for some $t\geq0$. By 
$\|\bg(\bx_{\text{feas}})\|\leq B_g$ (see Assumption~\ref{assume:general}E) and $\overline{\tau}=\sup_{t\geq0}\tau_t$, we have $\|\tau_t \bg(\bx_{\text{feas}})\|\leq\overline{\tau} B_g\leq M_{\blambda}/2$. Since $\|\blambda^{(t)}\|\leq M_{\blambda}$, we have
\begin{align}
    \nonumber
    \|\blambda^{(t)}\|^2 -2\tau_t\langle\blambda^{(t)}, - \bg(\bx_{\text{feas}})\rangle
    \leq&~\max_{\blambda\geq\bzero ,\,\|\blambda\|\leq M_{\blambda}} \left\{\|\blambda\|^2-2\tau_t\langle\blambda,-\bg(\bx_{\text{feas}})\rangle\right\}\\
    \label{eq:bound_lambda_t_lm}
    =&~\max_{\theta\in[0,1],\, \blambda\geq\bzero,\,\|\blambda\|= M_{\blambda}} \left\{\|\theta\blambda\|^2-2\theta\tau_t\langle\blambda,-\bg(\bx_{\text{feas}})\rangle\right\}.
\end{align} 
Fixing any $\blambda$ satisfying $\blambda\geq\bzero $ and $\|\blambda\|= M_{\blambda}$, consider the maximization in \eqref{eq:bound_lambda_t_lm} only over $\theta\in[0,1]$. Since $\tau_t\langle\blambda,-\bg(\bx_{\text{feas}})\rangle/\|\blambda\|^2\leq \|\tau_t \bg(\bx_{\text{feas}})\|/M_{\blambda}\leq1/2$,  $\theta=1$ is an optimal solution for \eqref{eq:bound_lambda_t_lm}. As a consequence,
\begin{align}
\nonumber
    \|\blambda^{(t)}\|^2 -2\tau_t\langle\blambda^{(t)}, - \bg(\bx_{\text{feas}})\rangle
    \leq&~\max_{\blambda\geq\bzero ,\,\|\blambda\|= M_{\blambda}} \left\{\|\blambda\|^2-2\tau_t\langle\blambda,-\bg(\bx_{\text{feas}})\rangle\right\}\\
    \label{eq:bound_lambda_M_lm}
    =&~\textstyle M_{\blambda}^2 -2\tau_t M_{\blambda} \cdot\min_{i\in[m]}[-g_i(\bx_{\text{feas}})].
\end{align} 
According to the updating equation of $\blambda^{(t+1)}$ in Step 4 of Algorithm~\ref{alg:imela}, we have 
\begin{align}
\nonumber
    \|\blambda^{(t+1)}\|^2 
    =&~\textstyle\sum_{i=1}^m ([\lambda_i^{(t)}+\tau_t\cdot g_i(\bx^{(t)})]_+)^2\\
    \nonumber
    \leq&~\textstyle\sum_{i=1}^m (\lambda_i^{(t)}+\tau_t\cdot g_i(\bx^{(t)}))^2
    =\|\blambda^{(t)}\|^2 +2\tau_t\langle\blambda^{(t)}, \bg(\bx^{(t)})\rangle +\tau_t^2\|\bg(\bx^{(t)})\|^2\\
    \nonumber
    \leq&~\|\blambda^{(t)}\|^2 +2\tau_t\big(C_{\blambda}-\langle\blambda^{(t)}, - \bg(\bx_{\text{feas}})\rangle\big) +\tau_t^2\|\bg(\bx^{(t)})\|^2\\
    \nonumber
    \leq&~\textstyle M_{\blambda}^2+2\overline{\tau} C_{\blambda}-2\underline{\tau} M_{\blambda}\cdot\min_{i\in[m]}[-g_i(\bx_{\text{feas}})] +\overline{\tau}^2B_g^2
    \leq M_{\blambda}^2,
\end{align} 
where the second inequality is from \eqref{eq:lambda_result_C_lambda_equiv}, and the last one is from \eqref{eq:bound_lambda_M_lm}, the definitions of $\overline{\tau}$, $\underline{\tau}$ and $M_{\blambda}$ in the settings. 
\eqref{eq:bound_lambda} then follows by the inductive result.
\end{proof}

\subsection{Proof of Lemma~\ref{lmm:subprob_complexity}}
\label{sec:subprob_complexity}
\begin{proof}[Proof of Lemma~\ref{lmm:subprob_complexity}]
Recall that function $\mathcal{L}_p(\bx,\bz^{(t)},\blambda^{(t+1)})$ is $(p-L)$-strongly convex and $K$-smooth in $\bx$ over $\X$ for any $t\geq0$. Suppose the APG method by~\citet[Eq.\,(2.2.63)]{nesterov2018lectures} is applied to $\min_{\bx\in\X}\mathcal{L}_p(\bx,\bz^{(t)},\blambda^{(t+1)})$ for $k_t$ iterations with any initial solution in $\bu^{(0)}\in\X$ and returns a solution $\bu^{(k_t)}$.
According to \citet[Theorem 2.2.1--2.2.3 \& Lemma 2.2.4]{nesterov2018lectures}, Lemma~\ref{lmm:bound_lambda}, and Assumption~\ref{assume:general}A, D, and E, we have 
\begin{align*}
    &~\mathcal{L}_p (\bu^{(k_t)},\bz^{(t)},\blambda^{(t+1)})
    -\mathcal{L}_p (\tilde{\bx}^{(t+1)},\bz^{(t)},\blambda^{(t+1)})\\
    \leq & 
    \left(1-\sqrt{\frac{p-L}{K}}\right)^{k_t}
    \left(\mathcal{L}_p (\bu^{(0)},\bz^{(t)},\blambda^{(t+1)})
    -\mathcal{L}_p (\tilde{\bx}^{(t+1)},\bz^{(t)},\blambda^{(t+1)})+\frac{p-L}{2}\|\bu^{(0)}-\tilde{\bx}^{(t+1)}\|^2\right)\\
    \leq & 
    \left(1-\sqrt{\frac{p-L}{K}}\right)^{k_t}\left((B_f+M_{\blambda} B_g)D_\X+\frac{2p-L}{2}D_\X^2\right).
\end{align*}
Hence, it must hold that 
\begin{equation}
\label{eq:x_primal_descent}
\begin{aligned}
    \mathcal{L}_p (\bu^{(k_t)},\bz^{(t)},\blambda^{(t+1)})-
    \mathcal{L}_p (\tilde{\bx}^{(t+1)},\bz^{(t)},\blambda^{(t+1)})
    \leq \frac{\epsilon_t^2}{8K}
\end{aligned}    
\end{equation}
for some $k_t$ satisfying
\begin{align}
\nonumber
    k_t=&~ \sqrt{\frac{K}{p-L}}\ln\left(\frac{(B_f+M_{\blambda} B_g)D_\X+\frac{2p-L}{2}D_\X^2}{\epsilon_t^2/(8K)}\right)\\
    \label{eq:subprob_complexity}
    =&~ O\left(\sqrt{\frac{K}{p-L}}\ln(\epsilon_t^{-1})\right).
\end{align}
By the definition of $\bx^{(t+1)}$, we have
\begin{equation}
\label{eq:xt_opt_condition}
\nabla_\bx\mathcal{L}_p(\bu^{(k_t)},\bz^{(t)},\blambda^{(t+1)})+K(\bx^{(t+1)}-\bu^{(k_t)})+\bv^{(t+1)}=0
\end{equation}
for some $\bv^{(t+1)}\in\N_\X(\bx^{(t+1)})$.
Finally, according to~\citet[Theorem 1]{nesterov2013gradient}, 
\begin{align*}
    &\text{dist}\big( -\nabla_\bx\mathcal{L}_p(\bx^{(t+1)},\bz^{(t)},\blambda^{(t+1)}),\N_\X(\bx^{(t+1)})\big)\\
    \leq&~ 2K\|\bu^{(k_t)}-\bx^{(t+1)}\|\\
    \leq&~ 2\sqrt{2K \left(\mathcal{L}_p (\bu^{(k_t)},\bz^{(t)},\blambda^{(t+1)})-
    \mathcal{L}_p (\tilde{\bx}^{(t+1)},\bz^{(t)},\blambda^{(t+1)}) \right)}\leq \epsilon_t,
\end{align*}
where the first inequality is because of \eqref{eq:xt_opt_condition} and the $K$-smoothness of $\mathcal{L}_p(\bx,\bz^{(t)},\blambda^{(t+1)})$, the second inequality is by \citet[Theorem 1]{nesterov2013gradient}, and
the last one is because of \eqref{eq:x_primal_descent}. This means \eqref{eq:subprob_optim_imela} is ensured with complexity $O(\ln(\epsilon_t^{-1}))$ for any $t\geq0$.
\end{proof}

\subsection{Proof of Proposition~\ref{thm:difference_phi}}
\label{sec:descent}
To prove Proposition~\ref{thm:difference_phi}, we need to first present three descent lemmas similar to the ones in~\citet{zhang2020proximal,zhang2022global,zhang2022iteration}. 
\begin{lemma}
\label{lmm:primal_descent}
The sequence $\{(\bx^{(t)},\blambda^{(t)},\bz^{(t)})\}_{t\geq0}$ generated by Algorithm~\ref{alg:imela} satisfies that
\begin{align}
\nonumber
    &~\mathcal{L}_p (\bx^{(t)},\bz^{(t)},\blambda^{(t)})-\mathcal{L}_p (\bx^{(t+1)},\bz^{(t+1)},\blambda^{(t+1)})\\
    \nonumber
    \geq&~ \frac{p-L}{2}\|\bx^{(t)}-\tilde{\bx}^{(t+1)}\|^2-\langle\blambda^{(t+1)}-\blambda^{(t)},\bg(\bx^{(t)})\rangle\\
    \label{eq:primal_descent}
    &+\frac{p}{2}(2/\theta_t-1)\|\bz^{(t)}-\bz^{(t+1)}\|^2-\frac{\epsilon_t^2}{2(p-L)}.
\end{align}
\end{lemma}

\begin{proof}
According to the definition of $\mathcal{L}_p$ in \eqref{eq:imela}, 
we have
\begin{align}
    \label{eq:primal_descent_lambda}
    \mathcal{L}_p (\bx^{(t)},\bz^{(t)},\blambda^{(t)})-\mathcal{L}_p (\bx^{(t)},\bz^{(t)},\blambda^{(t+1)})
    =-\langle\blambda^{(t+1)}-\blambda^{(t)},\bg(\bx^{(t)})\rangle.
\end{align}
By the $(p-L)$-strong convexity of $\mathcal{L}_p$ in $\bx$, we have
\begin{align}
\label{eq:primal_descent_x}
    \mathcal{L}_p (\bx^{(t)},\bz^{(t)},\blambda^{(t+1)}) -\mathcal{L}_p (\tilde{\bx}^{(t+1)},\bz^{(t)},\blambda^{(t+1)})
    \geq\frac{p-L}{2}\|\bx^{(t)}-\tilde{\bx}^{(t+1)}\|^2.
\end{align}
Moreover, we have from Step 6 of Algorithm~\ref{alg:imela} that
\begin{align}
    \nonumber
    &~\mathcal{L}_p (\bx^{(t+1)},\bz^{(t)},\blambda^{(t+1)})-\mathcal{L}_p (\bx^{(t+1)},\bz^{(t+1)},\blambda^{(t+1)})\\
    \nonumber
    =&~(p/2)\cdot(\|\bx^{(t+1)}-\bz^{(t)}\|^2-\|\bx^{(t+1)}-\bz^{(t+1)}\|^2)\\
    \nonumber
    =&~(p/2)\cdot\langle\bz^{(t+1)}-\bz^{(t)},(\bx^{(t+1)}-\bz^{(t)})+(\bx^{(t+1)}-\bz^{(t+1)})\rangle\\
    \label{eq:primal_descent_z}
    =&~(p/2)\cdot(2/\theta_t-1)\|\bz^{(t)}-\bz^{(t+1)}\|^2 
\end{align}
for $\theta_t\leq1$. Combining \eqref{eq:subprob_optimy_imela_corollary}, \eqref{eq:primal_descent_lambda}, \eqref{eq:primal_descent_x} and \eqref{eq:primal_descent_z}  yields \eqref{eq:primal_descent}.
\end{proof}

\begin{lemma}
\label{lmm:dual_descent}
The sequence $\{(\blambda^{(t)},\bz^{(t)})\}_{t\geq0}$ generated by Algorithm~\ref{alg:imela} satisfies that
\begin{align}
    \nonumber
    d(\blambda^{(t+1)},\bz^{(t+1)})-d(\blambda^{(t)},\bz^{(t)})
    \geq&~\langle\blambda^{(t+1)}-\blambda^{(t)}, \bg(\tilde{\bx}^{(t+1)})\rangle\\
    \label{eq:dual_descent}
    &+\frac{p}{2}\langle\bz^{(t+1)}-\bz^{(t)},\bz^{(t+1)}+\bz^{(t)}-2\bx(\blambda^{(t+1)},\bz^{(t+1)})\rangle.
\end{align}
\end{lemma}
\begin{proof}
By the definition of $\mathcal{L}_p$ in \eqref{eq:imela} and the definition of $d$ in \eqref{eq:d(lambda,z)}, we have
\begin{align*}
    d(\blambda^{(t+1)},\bz^{(t)})-d(\blambda^{(t)},\bz^{(t)})
    =&~\mathcal{L}_p (\tilde{\bx}^{(t+1)},\bz^{(t)},\blambda^{(t+1)})-\mathcal{L}_p (\bx(\blambda^{(t)},\bz^{(t)}),\bz^{(t)},\blambda^{(t)})\\
    \geq&~\mathcal{L}_p (\tilde{\bx}^{(t+1)},\bz^{(t)},\blambda^{(t+1)})-\mathcal{L}_p(\tilde{\bx}^{(t+1)},\bz^{(t)},\blambda^{(t)})\\
    =&~\langle\blambda^{(t+1)}-\blambda^{(t)},\bg(\tilde{\bx}^{(t+1)})\rangle.
\end{align*}

Next, using the same technique, we have
\begin{align*}
    &~ d(\blambda^{(t+1)},\bz^{(t+1)})-d(\blambda^{(t+1)},\bz^{(t)})\\
    =&~\mathcal{L}_p (\bx(\blambda^{(t+1)},\bz^{(t+1)}),\bz^{(t+1)},\blambda^{(t+1)})-\mathcal{L}_p(\tilde{\bx}^{(t+1)},\bz^{(t)},\blambda^{(t+1)})\\
    \geq&~\mathcal{L}_p (\bx(\blambda^{(t+1)},\bz^{(t+1)}),\bz^{(t+1)},\blambda^{(t+1)})-\mathcal{L}_p(\bx(\blambda^{(t+1)},\bz^{(t+1)}),\bz^{(t)},\blambda^{(t+1)})\\
    =&~
    \frac{p}{2}  (\|\bx(\blambda^{(t+1)},\bz^{(t+1)})-\bz^{(t+1)}\|^2-\|\bx(\blambda^{(t+1)},\bz^{(t+1)})-\bz^{(t)}\|^2)\\
    =&~
    \frac{p}{2} \langle\bz^{(t+1)}-\bz^{(t)},\bz^{(t+1)}+\bz^{(t)}-2\bx(\blambda^{(t+1)},\bz^{(t+1)})\rangle.
\end{align*}
Combining the above inequalities yields~\eqref{eq:dual_descent}.
\end{proof}
\begin{lemma}
\label{lmm:proximal_descent}
The sequence $\{\bz^{(t)}\}_{t\geq0}$ generated by Algorithm~\ref{alg:imela} satisfies that
\begin{align}
\label{eq:proximal_descent}
    v(\bz^{(t+1)})-v(\bz^{(t)})
    \leq p\langle\bz^{(t+1)}-\bz^{(t)},\bz^{(t)}-\bx(\bz^{(t)})\rangle
    +\frac{p(1/\sigma+1)}{2}\|\bz^{(t)}-\bz^{(t+1)}\|^2.
\end{align}
\end{lemma}
\begin{proof}
According to \eqref{eq:nablav} and Lemma~\ref{lmm:Lips_continuity}, for any $\bz$ and $\bz'$ in $\X$, it holds that
\begin{align*}
    \|\nabla v(\bz)-\nabla v(\bz')\|
    =&~\|p(\bz-\bx(\bz))-p(\bz'-\bx(\bz'))\|\\
    \leq&~p(\|\bz-\bz'\|+\|\bx(\bz)-\bx(\bz')\|)\\
    \leq&~ p(1/\sigma+1)\|\bz-\bz'\|,
\end{align*}
meaning that $v(\bz)$ is $p(1/\sigma+1)$-smooth, which implies \eqref{eq:proximal_descent}.
\end{proof}

Using the three descent lemmas above,  Proposition~\ref{thm:difference_phi} can be proved as follows. 
\begin{proof}[Proof Proposition~\ref{thm:difference_phi}]
Applying \eqref{eq:primal_descent}, \eqref{eq:dual_descent} and \eqref{eq:proximal_descent} to \eqref{eq:phi^t} and rearranging terms, we have
\begin{align}
    \nonumber
    \phi^t-\phi^{t+1}=&~\mathcal{L}_p (\bx^{(t)},\bz^{(t)},\blambda^{(t)})-\mathcal{L}_p (\bx^{(t+1)},\bz^{(t+1)},\blambda^{(t+1)})\\
    \nonumber
    &-2(d(\blambda^{(t)},\bz^{(t)})-d(\blambda^{(t+1)},\bz^{(t+1)})) +2(v(\bz^{(t)})-v(\bz^{(t+1)}))\\
    \nonumber
    \geq&~\frac{p-L}{2}\|\bx^{(t)}-\tilde{\bx}^{(t+1)}\|^2-\frac{\epsilon_t^2}{2(p-L)}+\frac{p}{2}(2/\theta_t-1)\|\bz^{(t)}-\bz^{(t+1)}\|^2\\
    \nonumber
    &-\langle\blambda^{(t+1)}-\blambda^{(t)},\bg(\bx^{(t)})\rangle +2\langle\blambda^{(t+1)}-\blambda^{(t)}, \bg(\tilde{\bx}^{(t+1)})\rangle\\
    \nonumber
    &+p\langle\bz^{(t+1)}-\bz^{(t)},\bz^{(t+1)}-\bz^{(t)}-2(\bx(\blambda^{(t+1)},\bz^{(t+1)})-\bx(\bz^{(t)}))\rangle\\
    \label{eq:difference_phi_descent_result}
    &-p(1/\sigma+1)\|\bz^{(t)}-\bz^{(t+1)}\|^2.
\end{align}
For any $\alpha>0$, by Young's inequality, it holds that
\begin{align}
    \label{eq:difference_phi_alpha}
    2\langle\bz^{(t+1)}-\bz^{(t)},\tilde{\bx}^{(t+1)}-\bx(\bz^{(t)})\rangle
    \leq\|\bz^{(t+1)}-\bz^{(t)}\|^2/\alpha+\alpha\|\tilde{\bx}^{(t+1)}-\bx(\bz^{(t)})\|^2.
\end{align}
By Cauchy-Schwartz inequality and \eqref{eq:Lips_continuity_x(z)}, we have 
\begin{align}
\nonumber
    &~2\langle\bz^{(t+1)}-\bz^{(t)},\bx(\blambda^{(t+1)},\bz^{(t+1)})-\tilde{\bx}^{(t+1)}\rangle\\
    \label{eq:difference_phi_sigma}
    \leq&~2\|\bz^{(t+1)}-\bz^{(t)}\|\cdot\|\bx(\blambda^{(t+1)},\bz^{(t+1)})-\tilde{\bx}^{(t+1)}\|
    \leq(2/\sigma)\cdot \|\bz^{(t+1)}-\bz^{(t)}\|^2.
\end{align}
Applying \eqref{eq:difference_phi_alpha} and \eqref{eq:difference_phi_sigma} to \eqref{eq:difference_phi_descent_result}, we have
\begin{align}
    \nonumber
    \phi^t-\phi^{t+1}
    \geq&~\frac{p-L}{2}\|\bx^{(t)}-\tilde{\bx}^{(t+1)}\|^2-\frac{\epsilon_t^2}{2(p-L)}-p\alpha\|\tilde{\bx}^{(t+1)}-\bx(\bz^{(t)})\|^2\\
    \nonumber
    &+p(1/\theta_t-{3/2}-1/\alpha-3/\sigma)\cdot\|\bz^{(t)}-\bz^{(t+1)}\|^2\\
    \label{eq:difference_phi_descent_terms}
    &+2\langle\blambda^{(t+1)}-\blambda^{(t)}, \bg(\tilde{\bx}^{(t+1)})-\bg(\bx^{(t)})\rangle+\langle\blambda^{(t+1)}-\blambda^{(t)},\bg(\bx^{(t)})\rangle.
\end{align}
According to the updating equation of $\blambda^{(t+1)}$ in Step 4 of Algorithm~\ref{alg:imela}, we have 
\begin{align}
\label{eq:imela_g_lambda_norm}
    -\tau_t \langle\blambda^{(t)},\bg(\bx^{(t)})\rangle \geq -\tau_t\langle\blambda^{(t+1)},\bg(\bx^{(t)})\rangle +\|\blambda^{(t+1)}-\blambda^{(t)}\|^2,&\\
    \label{eq:imela_g_lambda_cone}
    -\blambda^{(t+1)} +\blambda^{(t)}+\tau_t\bg(\bx^{(t)})\in \mathcal{N}_{\mathbb{R}_+^m}(\blambda^{(t+1)}).&
\end{align}
Therefore, we have
\begin{align}
    \nonumber
    &~2\langle\blambda^{(t+1)}-\blambda^{(t)}, \bg(\tilde{\bx}^{(t+1)})-\bg(\bx^{(t)})\rangle+\langle\blambda^{(t+1)}-\blambda^{(t)},\bg(\bx^{(t)})\rangle\\  \nonumber
    \geq&~2\langle\blambda^{(t+1)}-\blambda^{(t)}, \bg(\tilde{\bx}^{(t+1)})-\bg(\bx^{(t)})\rangle+\tau_t^{-1}\cdot\|\blambda^{(t+1)}-\blambda^{(t)}\|^2\\\nonumber
    =&-\tau_t\|\bg(\tilde{\bx}^{(t+1)})-\bg(\bx^{(t)})\|^2+ \frac{1}{\tau_t}\|\blambda^{(t+1)}-\blambda^{(t)}+\tau_t(\bg(\tilde{\bx}^{(t+1)})-\bg(\bx^{(t)}))\|^2\\
    \nonumber
    \geq&-\tau_t\|\bg(\tilde{\bx}^{(t+1)})-\bg(\bx^{(t)})\|^2+\frac{1}{\tau_t}\cdot \text{dist}\big(\tau_t\bg(\tilde{\bx}^{(t+1)}), \mathcal{N}_{\mathbb{R}_+^m}(\blambda^{(t+1)})\big)^2\\
    \label{eq:difference_phi_result_tau_lm}
    =&-\tau_t\|\bg(\tilde{\bx}^{(t+1)})-\bg(\bx^{(t)})\|^2+\tau_t\cdot \text{dist}\big(\bg(\tilde{\bx}^{(t+1)}), \mathcal{N}_{\mathbb{R}_+^m}(\blambda^{(t+1)})\big)^2,
\end{align}
where the two inequalities are by \eqref{eq:imela_g_lambda_norm} and \eqref{eq:imela_g_lambda_cone}, respectively, and the last equality is from the structure of $\mathcal{N}_{\mathbb{R}_+^m}(\blambda^{(t+1)})$.


Applying \eqref{eq:difference_phi_result_tau_lm} and Assumption~\ref{assume:general}E to \eqref{eq:difference_phi_descent_terms}, we have 
\begin{align}
\nonumber
    \phi^t-\phi^{t+1}
    \nonumber
    \geq&\left(\frac{p-L}{2}-\tau_tB_g^2\right)\|\bx^{(t)}-\tilde{\bx}^{(t+1)}\|^2-\frac{\epsilon_t^2}{2(p-L)}-p\alpha\|\tilde{\bx}^{(t+1)}-\bx(\bz^{(t)})\|^2\\
    \nonumber
    &+p(1/\theta_t-{3/2}-1/\alpha-3/\sigma)\cdot\|\bz^{(t)}-\bz^{(t+1)}\|^2\\
    \label{eq:difference_phi_descent_coeffs}
    &+\tau_t\cdot \text{dist}\big(\bg(\tilde{\bx}^{(t+1)}), \mathcal{N}_{\mathbb{R}_+^m}(\blambda^{(t+1)})\big)^2.
\end{align}
Let $\alpha=6\theta$. Recall that $\theta_t=\theta \leq \frac{p-L}{18p}<\frac{1}{18}$ and $\sigma=\frac{p-L}{p}$. We have
\begin{align}
\label{eq:difference_phi_theta_coefficient}
   \frac{1}{\theta_t}-{3/2}-1/\alpha-3/\sigma\geq\frac{1}{2\theta_t}-1/\alpha-3/\sigma=\frac{1}{3\theta}-\frac{3p}{p-L}\ge\frac{1}{6\theta}.
\end{align}
Recall that $\tau_t=\tau=\frac{p-L}{4B_g^2}$. We have $\frac{p-L}{2}-\tau_tB_g^2=\frac{p-L}{4}$. 
Applying this equation, 
\eqref{eq:difference_phi_theta_coefficient}, 
and $\alpha=6\theta$ to \eqref{eq:difference_phi_descent_coeffs} leads to \eqref{eq:difference_phi}.
\end{proof}

\section{Sufficient Decrease of Potential Function $\phi_t$}
\label{sec:bounded_phi}
In this section, we present a series of technical results that establish the sufficient decrease of the potential function $\phi_t$ after each iteration of Algorithm~\ref{alg:imela}. Such sufficient decrease is the key to establish the total complexity. 


\subsection{Proof of Proposition~\ref{thm:regularity_exponent_d}
}
\label{sec:regularity}
To prove Proposition~\ref{thm:regularity_exponent_d}, 
we first need the following lemma.
\begin{lemma}
\label{thm:nearKKT_impliesnearx} 
For any $\delta' >0$, there must exist $R(\delta') > 0$ such that, if 
\begin{equation}
\label{eq:regular_cond_delta_near}
\begin{aligned}  
     \textup{dist}\big(-\nabla f(\tilde{\bx}) - \nabla\bg(\tilde{\bx})\blambda, \N_\X(\tilde{\bx})\big)^2
    +\textup{dist}\big( \bg(\tilde{\bx}), \mathcal{N}_{\mathbb{R}_+^m}(\blambda)\big)^2 \leq R^2(\delta')
\end{aligned} 
\end{equation}
for some $\tilde\bx\in\X$ and some $\blambda\geq\bzero$ satisfying $\|\blambda\|\leq M_{\blambda}$, there must exist $\bx^*\in\X^*$ such that $\|\tilde\bx-\bx^*\|\leq\delta'$ and $\textup{supp}(\blambda)\subset J_g(\bx^*)\text{ and }J_A(\tilde\bx)\subset J_A(\bx^*)$. 
\end{lemma}
\begin{proof}
We prove this by contradiction. Suppose that the claim does not hold. Then there exists a sequence of $\{(\tilde{\bx}^{(k)},\tilde{\blambda}^{(k)})\}_{k\geq1}\subset\X\times\mathbb{R}_+^m$, $\|\tilde{\blambda}^{(k)}\|\leq M_{\blambda}$, such that
\begin{align}
\label{eq:regular_cond_delta_near_claim}
    \text{dist}\big(-\nabla f(\tilde{\bx}^{(k)}) - \nabla\bg(\tilde{\bx}^{(k)})\tilde{\blambda}^{(k)}, \N_\X(\tilde{\bx}^{(k)})\big)^2
    +\text{dist}\big( \bg(\tilde{\bx}^{(k)}), \mathcal{N}_{\mathbb{R}_+^m}(\tilde{\blambda}^{(k)})\big)^2
    \leq 1/k,
\end{align}
but for any $k$, there is no $\bx^*\in\X^*$ such that $\|\tilde\bx^{(k)}-\bx^*\|\leq\delta'$ and $\textup{supp}(\tilde\blambda^{(k)})\subset J_g(\bx^*)$ and $J_A(\tilde\bx^{(k)})\subset J_A(\bx^*)$.

Due to the compactness of $\X$ (see Assumption~\ref{assume:general}A) and the set $\{\blambda\in\mathbb{R}_+^m: \|\blambda\|\leq M_{\blambda}\}$, 
by passing to a subsequence if necessary, there exists some pair $(\tilde{\bx}^*, \tilde{\blambda}^*)\in\X\times\mathbb{R}_+^m$, $\|\tilde{\blambda}^*\|\leq M_{\blambda}$, such that 
$(\tilde{\bx}^{(k)},\tilde{\blambda}^{(k)})\rightarrow(\tilde{\bx}^*,\tilde{\blambda}^*)$ as $k\rightarrow\infty$. 
In addition, by the continuity of $\nabla f(\cdot)$ and $\nabla\bg(\cdot)$ and the outer semi-continuity of $\N_\X(\cdot)$~\citep[Proposition 6.6]{rockafellar1998variational}, it follows from \eqref{eq:regular_cond_delta_near_claim} that 
\begin{align*}
    \text{dist}\big(-\nabla f(\tilde{\bx}^*) - \nabla\bg(\tilde{\bx}^*)\tilde{\blambda}^*, \N_\X(\tilde{\bx}^*)\big)^2 
    +\text{dist}\big( \bg(\tilde{\bx}^*), \mathcal{N}_{\mathbb{R}_+^m}(\tilde{\blambda}^*)\big)^2 =0,
\end{align*}
which easily implies that $\tilde{\bx}^*\in\X^*$. Also, by the continuity of $\bA(\cdot)-\bb$, if $[\bA\tilde{\bx}^*-\bb]_j < 0$, then $[\bA\tilde{\bx}^{(k)}-\bb]_j < 0$ as well when $k$ is large enough, and thus $J_A(\tilde\bx^{(k)})\subset J_A(\tilde\bx^*)$. Moreover, let $\tilde{I}_g^{(k)}=\text{supp}(\tilde{\blambda}^{(k)})$. Then 
\begin{align*}
    \textstyle
    \sum_{i\in \tilde{I}_g^{(k)}} |g_i(\tilde{\bx}^{(k)})|^2 + \sum_{i\in [m]\backslash\tilde{I}_g^{(k)}} \big([g_i(\tilde{\bx}^{(k)})]_+\big)^2 
    =\text{dist}\big( \bg(\tilde{\bx}^{(k)}), \mathcal{N}_{\mathbb{R}_+^m}(\tilde{\blambda}^{(k)})\big)^2
    \leq 1/k.
\end{align*}
Now let $k$ be large enough such that for any $i\not\in J_g(\tilde{\bx}^*)$, i.e., $g_i(\tilde{\bx}^*) < 0$, it holds that $g_i(\tilde{\bx}^{(k)}) \leq \frac{1}{2}g_i(\tilde{\bx}^*)$ and $|g_i(\tilde{\bx}^*)|^2 > \frac{4}{k}$. This observation together with the above inequality implies $\text{supp}(\tilde{\blambda}^{(k)})=\tilde{I}_g^{(k)}\subset J_g(\tilde{\bx}^*)$. Hence, we find $\tilde{\bx}^*\in\X^*$ such that $\|\tilde\bx^{(k)}-\tilde\bx^*\|\leq\delta'$ and $\textup{supp}(\tilde\blambda^{(k)})\subset J_g(\tilde\bx^*)\text{ and }J_A(\tilde\bx^{(k)})\subset J_A(\tilde\bx^*)$ when $k$ is large enough. This contradicts to the assumption.
Therefore, the conclusion holds.
\end{proof}
With Lemma~\ref{thm:nearKKT_impliesnearx}, we are able to prove Proposition~\ref{thm:regularity_exponent_d}.
\begin{proof}[Proof of Proposition~\ref{thm:regularity_exponent_d}]
Let $\delta$ be given by Assumption~\ref{assume:local_EBC_exponent_d} and let $R(\cdot)$ be as in Lemma~\ref{thm:nearKKT_impliesnearx}. 
By \eqref{eq:d(lambda,z)}, we have
$\tilde{\bxi}^{(t+1)}:=-\nabla_\bx\mathcal{L}_p (\tilde{\bx}^{(t+1)},\bz^{(t)},\blambda^{(t+1)})\in \N_\X(\tilde{\bx}^{(t+1)}).$
Since $\mathcal{L}_p$ is $K$-smooth in $\bx$, it holds for any $\bx\in\X$ that
\begin{small}
\begin{align}
\label{eq:regularity_result_ksi_exponent_d}
    \mathcal{L}_p (\bx,\bz^{(t)},\blambda^{(t+1)}) - \mathcal{L}_p (\tilde{\bx}^{(t+1)},\bz^{(t)},\blambda^{(t+1)})
    \leq -\langle\tilde{\bxi}^{(t+1)},\bx-\tilde{\bx}^{(t+1)}\rangle
    +(K/2)\cdot \|\bx-\tilde{\bx}^{(t+1)}\|^2.
\end{align}
\end{small}By the $(p-L)$-strongly convexity of $\mathcal{L}_p$ in $\bx$, it holds that
\begin{align}
\label{eq:regularity_result_sc_exponent_d}
    \frac{p-L}{2}\|\bx-\tilde{\bx}^{(t+1)}\|^2
    \leq\mathcal{L}_p (\bx,\bz^{(t)},\blambda^{(t+1)})-\mathcal{L}_p (\tilde{\bx}^{(t+1)},\bz^{(t)},\blambda^{(t+1)}).
\end{align}
Adding \eqref{eq:regularity_result_ksi_exponent_d} and \eqref{eq:regularity_result_sc_exponent_d} with $\bx=\bx(\bz^{(t)})$ yields
\begin{align}
\nonumber
    \frac{p-L}{2}\|\bx(\bz^{(t)})-\tilde{\bx}^{(t+1)}\|^2
    \leq&~\mathcal{L}_p (\bx(\bz^{(t)}),\bz^{(t)},\blambda^{(t+1)})-\mathcal{L}_p  (\bx,\bz^{(t)},\blambda^{(t+1)})\\
    \label{eq:regularity_result_ksi_tilde_exponent_d}
    &-\langle\tilde{\bxi}^{(t+1)},\bx-\tilde{\bx}^{(t+1)}\rangle+(K/2)\cdot \|\bx-\tilde{\bx}^{(t+1)}\|^2.
\end{align}

Recall that $\|\blambda^{(t+1)}\|\leq M_{\blambda}$. According to Lemma~\ref{thm:nearKKT_impliesnearx} and \eqref{eq:regular_cond_KKT_near_R_exponent_d}, there must exist $\bx^*\in\X^*$ such that $\|\tilde\bx^{(t+1)}-\bx^*\|\leq\delta$, $I_g^{(t+1)}=\text{supp}(\blambda^{(t+1)})\subset J_g(\bx^*)$ and $I_A^{(t+1)}=J_A(\tilde\bx^{(t+1)})\subset J_A(\bx^*)$.
Let $\bar{\bx}^{(t+1)}$ be the projection of $\tilde\bx^{(t+1)}$ onto $\mathcal{S}(I_g^{(t+1)}, I_A^{(t+1)})$, so that  $\text{dist} (\tilde\bx^{(t+1)},\mathcal{S}(I_g^{(t+1)}, I_{A}^{(t+1)}))=\|\tilde\bx^{(t+1)}-\bar{\bx}^{(t+1)}\|$. As a result, it holds that
\begin{align}
\label{eq:regularity_xbar_g_exponent_d}
    &~g_i(\bar{\bx}^{(t+1)})=0,\;i\in I_g^{(t+1)},\;g_i(\bar{\bx}^{(t+1)})\leq 0,\;i\in[m]\backslash I_g^{(t+1)},\\
\label{eq:regularity_xbar_polyhedral_exponent_d}
    &~[\bA\bar{\bx}^{(t+1)}-\bb]_j=0,\;j\in I_{A}^{(t+1)},\;[\bA\bar{\bx}^{(t+1)}-\bb]_j\leq0,\;j\in[l]\backslash I_{A}^{(t+1)},
\end{align}
Applying Assumption~\ref{assume:local_EBC_exponent_d} and using $I_A^{(t+1)}=J_A(\tilde\bx^{(t+1)})$, we obtain
\begin{align}
\nonumber
    \|\bar{\bx}^{(t+1)}-\tilde{\bx}^{(t+1)}\|^{2d} \leq&~ 
    %
    \gamma^2 \cdot\text{dist} \big(\bg(\tilde{\bx}^{(t+1)}), \mathcal{N}_{\mathbb{R}_+^m}(\blambda^{(t+1)})\big)^2 \\
    \label{eq:regularity_result_norm_dist_exponent_d}
    \iff \|\bar{\bx}^{(t+1)}-\tilde{\bx}^{(t+1)}\|^{2} \leq&~\gamma^{2/d} \cdot\text{dist} \big(\bg(\tilde{\bx}^{(t+1)}), \mathcal{N}_{\mathbb{R}_+^m}(\blambda^{(t+1)})\big)^{2/d}.
\end{align}

Let $\bA(I_{A}^{(t+1)},:)$ be the matrix formed by the rows of $\bA$ with indices in $I_{A}^{(t+1)}$. 
Since $\tilde{\bxi}^{(t+1)}\in\N_\X(\tilde{\bx}^{(t+1)})$, there exists some $\bw^{(t+1)}\in\mathbb{R}_+^{|I_{A}^{(t+1)}|}$ such that 
    $\tilde{\bxi}^{(t+1)} 
    =[\bA(I_{A}^{(t+1)},:)]^\top \bw^{(t+1)},$
which, according to \eqref{eq:regularity_xbar_polyhedral_exponent_d}, implies
\begin{align}
\label{eq:regularity_result_ksi_equality_exponent_d}
    \langle\tilde{\bxi}^{(t+1)}, \bar{\bx}^{(t+1)}-\tilde{\bx}^{(t+1)}\rangle=0.
\end{align}
Setting $\bx=\bar{\bx}^{(t+1)}$ in \eqref{eq:regularity_result_ksi_tilde_exponent_d} and applying \eqref{eq:regularity_result_norm_dist_exponent_d} and \eqref{eq:regularity_result_ksi_equality_exponent_d}, 
we have
\begin{align}
    \nonumber
    \frac{p-L}{2}\|\bx(\bz^{(t)})-\tilde{\bx}^{(t+1)}\|^2 
    \leq&~\mathcal{L}_p (\bx(\bz^{(t)}),\bz^{(t)},\blambda^{(t+1)})-\mathcal{L}_p (\bar{\bx}^{(t+1)},\bz^{(t)},\blambda^{(t+1)})\\
    \label{eq:regularity_result_xbar_dist_exponent_d}
    &+(\gamma^{1/d} K/2)\cdot\text{dist}\big(\bg(\tilde{\bx}^{(t+1)}), \mathcal{N}_{\mathbb{R}_+^m}(\blambda^{(t+1)})\big)^{2/d}.
\end{align}
By \eqref{eq:x(z)} and \eqref{eq:regularity_xbar_g_exponent_d}, we must have $\lambda_i^{(t+1)}g_i(\bx(\bz^{(t)}))\leq0$ and $\lambda_i^{(t+1)}g_i(\bar{\bx}^{(t+1)})=0$ for $i\in[m]$. This implies
\begin{align}
\nonumber
    &~\mathcal{L}_p (\bx(\bz^{(t)}),\bz^{(t)},\blambda^{(t+1)})-\mathcal{L}_p (\bar{\bx}^{(t+1)},\bz^{(t)},\blambda^{(t+1)})\\
    \label{eq:regularity_result_imela_diff_exponent_d}
    \leq& \left(f(\bx(\bz^{(t)}))+\frac{p}{2}\|\bx(\bz^{(t)})-\bz^{(t)}\|^2\right)-\left(f(\bar{\bx}^{(t+1)})+\frac{p}{2}\|\bar{\bx}^{(t+1)}-\bz^{(t)}\|^2\right)
    \leq0,
\end{align}
where the last inequality is because $\bar{\bx}^{(t+1)}$ is a feasible solution to the minimization problem in \eqref{eq:x(z)} while
$\bx(\bz^{(t)})$ is its optimal solution. Combining \eqref{eq:regularity_result_imela_diff_exponent_d} with \eqref{eq:regularity_result_xbar_dist_exponent_d} yields the claim.
\end{proof}

\subsection{Proof of Lemma~\ref{lmm:weak_dual_bound}}
\label{sec:weak_dual_bound}
\begin{proof}[Proof of Lemma~\ref{lmm:weak_dual_bound}]
By the $(p-L)$-strong convexity of $\mathcal{L}_p$ in $\bx$, we have 
\begin{align*}
    &
    ~\frac{p-L}{2}\|\bx(\bz^{(t)})-\tilde{\bx}^{(t+1)}\|^2 \leq\mathcal{L}_p (\bx(\bz^{(t)}),\bz^{(t)}, \blambda^{(t+1)}) - \mathcal{L}_p (\tilde{\bx}^{(t+1)},\bz^{(t)},\blambda^{(t+1)}),\\
    &
    ~\frac{p-L}{2}\|\tilde{\bx}^{(t+1)}-\bx(\bz^{(t)})\|^2\leq\mathcal{L}_p
    (\tilde{\bx}^{(t+1)}, \bz^{(t)}, \blambda(\bz^{(t)})) - \mathcal{L}_p(\bx(\bz^{(t)}),\bz^{(t)}, \blambda(\bz^{(t)})),
\end{align*}
where $\blambda(\bz^{(t)})$ is a vector of Lagrangian multipliers corresponding to $\bx(\bz^{(t)})$ defined in \eqref{eq:x(z)}. Notice that by the strong duality (see, e.g.~\citet[Sec. 28--30]{rockafellar1970convex}), we have $\blambda(\bz^{(t)})\in\argmax_{\blambda\geq\bzero}\mathcal{L}_p(\bx(\bz^{(t)}),\bz^{(t)},\blambda)$. Hence, adding the above two inequalities gives
\begin{align}
\nonumber
    (p-L)\cdot \|\tilde{\bx}^{(t+1)}- \bx(\bz^{(t)})\|^2
    \nonumber
    \leq&~\mathcal{L}_p (\tilde{\bx}^{(t+1)},\bz^{(t)},\blambda(\bz^{(t)})) - \mathcal{L}_p (\tilde{\bx}^{(t+1)}, \bz^{(t)},\blambda^{(t+1)})\\
    \nonumber
    =&~\langle\blambda(\bz^{(t)})-\blambda^{(t+1)}, \bg(\tilde{\bx}^{(t+1)})\rangle\\
    \nonumber
    \leq&~\textstyle \sum_{i\in\mathrm{supp}(\blambda^{(t+1)})} |\lambda_i(\bz^{(t)}) - \lambda_i^{(t+1)}|\cdot |g_i(\tilde{\bx}^{(t+1)})|\\
    \nonumber
    &+\textstyle \sum_{i\in[m]\backslash \mathrm{supp}(\blambda^{(t+1)})} \lambda_i(\bz^{(t)}) \cdot [g_i(\tilde{\bx}^{(t+1)})]_+\\
    \label{eq:weak_dual_bound_lambda}
    \leq&~\|\blambda(\bz^{(t)}) - \blambda^{(t+1)}\| \cdot \text{dist}\big( \bg(\tilde{\bx}^{(t+1)}), \mathcal{N}_{\mathbb{R}_+^m}(\blambda^{(t+1)})\big),
\end{align}
where the last inequality follows from the Cauchy--Schwarz inequality.

In addition, by the KKT conditions at $(\bx(\bz^{(t)}),\blambda(\bz^{(t)}))$ in \eqref{eq:x(z)}, we have
\begin{align}
\label{eq:x(z)_KKT_stat}
    \textstyle  \bxi:= -\nabla f(\bx(\bz^{(t)})) - \sum_{i=1}^m \lambda_i(\bz^{(t)}) \nabla g_i(\bx(\bz^{(t)})) - p(\bx(\bz^{(t)})-\bz^{(t)}) \in \N_\X(\bx(\bz^{(t)})),&\\
    \label{eq:x(z)_KKT_comple_slack}
\lambda_i(\bz^{(t)}) g_i(\bx(\bz^{(t)})) = 0,\;i\in [m].&
\end{align} 
Hence, by the convexity of $g_i$'s and $\blambda(\bz^{(t)})\geq\bzero$, it follows that
\begin{align*}
    \textstyle\sum_{i=1}^m \lambda_i(\bz^{(t)}) \cdot g_i(\bx_{\text{feas}})
    \quad\geq\;\;\;\; &~\textstyle\sum_{i=1}^m \lambda_i(\bz^{(t)}) \left(g_i(\bx(\bz^{(t)})) + \big\langle \nabla g_i(\bx(\bz^{(t)})), \bx_{\text{feas}}-\bx(\bz^{(t)}) \big\rangle\right)\\
    \overset{ \eqref{eq:x(z)_KKT_stat}, \eqref{eq:x(z)_KKT_comple_slack}}
    {=}&~\big\langle \bxi + \nabla f(\bx(\bz^{(t)})) + p(\bx(\bz^{(t)})-\bz^{(t)}), 
    \bx(\bz^{(t)})-\bx_{\text{feas}}  \big\rangle\\
    \qquad\geq\;\;\;\; &~\big\langle\nabla f(\bx(\bz^{(t)})) + p(\bx(\bz^{(t)})-\bz^{(t)}), \bx(\bz^{(t)})-\bx_{\text{feas}}\big\rangle,
\end{align*}
where $\bx_{\text{feas}}$ is given in Assumption~\ref{assume:general}F, and the last inequality holds because $\bxi\in \N_\X(\bx(\bz^{(t)}))$. Now applying Assumption~\ref{assume:general}A, D, and F, we obtain
\begin{align*}
    &~\textstyle \|\blambda(\bz^{(t)})\|_1\cdot \min_{i\in[m]}[-g_i(\bx_{\text{feas}})]\\
    \leq&~\langle\nabla f(\bx(\bz^{(t)})) + p(\bx(\bz^{(t)})-\bz^{(t)}), \bx_{\text{feas}}-\bx(\bz^{(t)})\rangle\\
    \leq&~\big(\|\nabla f(\bx(\bz^{(t)}))\|+p\cdot\|\bx(\bz^{(t)})-\bz^{(t)}\|\big) \cdot\|\bx_{\text{feas}}-\bx(\bz^{(t)})\|\\
    \leq&~(B_f+p D_{\X})\cdot D_{\X},
\end{align*} 
and thus $\|\blambda(\bz^{(t)})\|\leq \|\blambda(\bz^{(t)})\|_1\leq \frac{(B_f+p D_{\X})\cdot D_{\X}}{\min_{i\in[m]}[-g_i(\bx_{\text{feas}})]} \leq M_{\blambda},$
which together with \eqref{eq:weak_dual_bound_lambda} and \eqref{eq:bound_lambda} gives the desired result. 
\end{proof}

\subsection{Proof of Proposition~\ref{thm:bounded_phi_exponent_d}}
\label{sec:bounded_phi_proof}
\begin{proof}[Proof of Proposition~\ref{thm:bounded_phi_exponent_d}]
By the choices of $\tau_t$ and $\theta_t$, inequality \eqref{eq:difference_phi} holds by 
Proposition~\ref{thm:difference_phi}. We then prove \eqref{eq:bounded_phi_exponent_d} by considering two cases depending on whether  \eqref{eq:regular_cond_KKT_near_R_exponent_d} holds. 

Case I: \eqref{eq:regular_cond_KKT_near_R_exponent_d} is satisfied.
In this case, inequality \eqref{eq:regularity_result_exponent_d} holds by
Proposition~\ref{thm:regularity_exponent_d}. Let
$$
s:=\text{dist}\big( \bg(\tilde{\bx}^{(t+1)}), \mathcal{N}_{\mathbb{R}_+^m}(\blambda^{(t+1)}) \big), \quad C_d:=\frac{\gamma^{1/d} K}{p-L}.
$$
Then \eqref{eq:regularity_result_exponent_d} implies
$$
\|\tilde{\bx}^{(t+1)}-\bx(\bz^{(t)})\|^2 \leq C_d s^{2/d},
$$
and hence
\begin{align*}
    &~ \tau\cdot\text{dist}\big( \bg(\tilde{\bx}^{(t+1)}), \mathcal{N}_{\mathbb{R}_+^m}(\blambda^{(t+1)})\big)^2-6p\theta\|\tilde{\bx}^{(t+1)}-\bx(\bz^{(t)})\|^2\\
    \geq&~
    \tau s^2-6p\theta C_d s^{2/d}.
\end{align*}
When $d\in(1,2]$, let
$$
r:=\frac{d}{d-1},
\quad
r':=d,
\quad
\frac{1}{r}+\frac{1}{r'}=1.
$$
We apply Young's inequality in the form
$$
AB \leq A^{r}/r+\frac{B^{r'}}{r'},
\quad A,B\geq0.
$$
Choose
$$
A:=\left(\frac{2}{\tau d}\right)^{1/d} 6p\theta C_d,
\quad
B:=\left(\frac{\tau d}{2}\right)^{1/d}s^{2/d}.
$$
Then 
\begin{align*}
    6p\theta C_d s^{2/d} = AB \leq&~ \frac{1}{r}
    \left(\left(\frac{2}{\tau d}\right)^{1/d} 6p\theta C_d\right)^{r}
    + \frac{1}{r'}
    \left(\left(\frac{\tau d}{2}\right)^{1/d} s^{2/d}\right)^{r'}\\
    =&~ \left(1-\frac{1}{d}\right)
    \left(\frac{2}{\tau d}\right)^{\frac{1}{d-1}}
    \big(6p\theta C_d\big)^{\frac{d}{d-1}}
    + \frac{\tau}{2}s^2,
\end{align*}
where the equality uses $\frac{1}{r}=1-\frac{1}{d}$ and $\frac{1}{r'}=\frac{1}{d}$. Rearranging yields
\begin{align*}
    \tau s^2-6p\theta C_d s^{2/d}
    \geq&~
    \frac{\tau}{2}s^{2}
    -
    \left(1-\frac{1}{d}\right)
    \left(\frac{2}{\tau d}\right)^{\frac{1}{d-1}}
    \big(6p\theta C_d\big)^{\frac{d}{d-1}}\\
    =&~\frac{\tau}{2}s^2-\Xi_d(\theta),
\end{align*}
where $\Xi_d(\theta)$ is defined in \eqref{eq:Xi_d_def}.

When $d=1$, inequality \eqref{eq:regularity_result_exponent_d} reduces to the linear bound
$$
\|\tilde{\bx}^{(t+1)}-\bx(\bz^{(t)})\|^2 \leq C_1 s^2, ~\text{ with }~ C_1:=\frac{\gamma K}{p-L}.
$$
Hence
$$
\tau s^2-6p\theta\|\tilde{\bx}^{(t+1)}-\bx(\bz^{(t)})\|^2
\geq
(\tau-6p\theta C_1)s^2.
$$
By the definition of $\hat{\theta}_d$ with $d=1$, we have $\theta\leq\hat{\theta}_d=\frac{\tau}{12pC_1}$, which yields
$$
6p\theta C_1 \leq \frac{\tau}{2}, ~\text{ and }~ \tau-6p\theta C_1 \geq \frac{\tau}{2}.
$$
Hence,
\begin{align*}
    &~\tau\cdot \text{dist}\big( \bg(\tilde{\bx}^{(t+1)}), \mathcal{N}_{\mathbb{R}_+^m}(\blambda^{(t+1)})\big)^2-6p\theta\|\tilde{\bx}^{(t+1)}-\bx(\bz^{(t)})\|^2\\
    \geq&~  \frac{\tau}{2} \text{dist}\big( \bg(\tilde{\bx}^{(t+1)}), \mathcal{N}_{\mathbb{R}_+^m}(\blambda^{(t+1)})\big)^2,
\end{align*}
which corresponds to $\Xi_1(\theta)=0$ given in \eqref{eq:Xi_d_def}.

In both subcases,
\begin{align}
\nonumber
    &~\tau\cdot\text{dist}\big( \bg(\tilde{\bx}^{(t+1)}), \mathcal{N}_{\mathbb{R}_+^m}(\blambda^{(t+1)})\big)^2-6p\theta\|\tilde{\bx}^{(t+1)}-\bx(\bz^{(t)})\|^2\\
    \label{eq:bounded_phi_regularity_result_1_exponent_d}
    \geq&~\frac{\tau}{2}\text{dist}\big( \bg(\tilde{\bx}^{(t+1)}), \mathcal{N}_{\mathbb{R}_+^m}(\blambda^{(t+1)})\big)^2 -\Xi_d(\theta).
\end{align}
Then applying \eqref{eq:bounded_phi_regularity_result_1_exponent_d} to \eqref{eq:difference_phi} implies \eqref{eq:bounded_phi_exponent_d}.

Case II: Condition \eqref{eq:regular_cond_KKT_near_R_exponent_d} does not hold. In this case, it holds that
\begin{equation}
\label{eq:regular_cond_delta_near_violated_exponent_d}
\begin{aligned}  
    \text{dist}\big(-\nabla f(\tilde{\bx}^{(t+1)}) - \nabla\bg(\tilde{\bx}^{(t+1)})\blambda^{(t+1)}, \N_\X(\tilde{\bx}^{(t+1)})\big)^2 & \\
    +\;\text{dist}\big( \bg(\tilde{\bx}^{(t+1)}), \mathcal{N}_{\mathbb{R}_+^m}(\blambda^{(t+1)})\big)^2
    & > R^2(\delta). 
\end{aligned} 
\end{equation}
By the triangle inequality and $-\nabla_\bx \mathcal{L}_p (\tilde{\bx}^{(t+1)},\bz^{(t)},\blambda^{(t+1)})\in\N_\X(\tilde{\bx}^{(t+1)})$, it holds that
\begin{align}
\nonumber
    &~\text{dist}\big(-\nabla f(\tilde{\bx}^{(t+1)}) - \nabla\bg(\tilde{\bx}^{(t+1)})\blambda^{(t+1)}, \N_\X(\tilde{\bx}^{(t+1)})\big)^2\\
    \nonumber
    \leq&~ p^2 \|\tilde{\bx}^{(t+1)} - \bz^{(t)}\|^2
    \leq 2p^2 \big( \|\tilde{\bx}^{(t+1)} - \bx^{(t+1)}\|^2 + \|\bx^{(t+1)} - \bz^{(t)}\|^2 \big)\\
    \nonumber
    =&~ 2p^2\big( \|\tilde{\bx}^{(t+1)} - \bx^{(t+1)}\|^2 + (1/\theta^2)\cdot\|\bz^{(t+1)} - \bz^{(t)}\|^2 \big)\\
    \label{eq:upper-bd-df_exponent_d}
    \leq&~ 2p^2 \left( \frac{\epsilon_t^2}{(p-L)^2} + \frac{1}{\theta^2}\|\bz^{(t+1)} - \bz^{(t)}\|^2 \right),
\end{align}
where the last inequality follows
from \eqref{eq:subprob_optimy_imela_corollary} and \eqref{eq:regularity_result_sc_exponent_d} with $\bx = \bx^{(t+1)}$. Plugging
\eqref{eq:upper-bd-df_exponent_d} into \eqref{eq:regular_cond_delta_near_violated_exponent_d}, we obtain
\begin{align}
\label{eq:regular_cond_delta_near_violated-2_exponent_d}
    \text{dist}\big( \bg(\tilde{\bx}^{(t+1)}), \mathcal{N}_{\mathbb{R}_+^m}(\blambda^{(t+1)})\big)^2
    +2p^2\left(\frac{\epsilon_t^2}{(p-L)^2} + \frac{1}{\theta^2}\|\bz^{(t+1)} - \bz^{(t)}\|^2\right)  > R^2(\delta). 
\end{align}
We then have
\begin{small}
\begin{align}
\nonumber
    \;&~
    \frac{p}{6\theta}\|\bz^{(t)}-\bz^{(t+1)}\|^2 + \tau\cdot\text{dist}\big( \bg(\tilde{\bx}^{(t+1)}), \mathcal{N}_{\mathbb{R}_+^m}(\blambda^{(t+1)})\big)^2 -6p\theta\|\tilde{\bx}^{(t+1)}-\bx(\bz^{(t)})\|^2\\
    \nonumber
    \overset{\eqref{eq:weak_dual_bound}}{\geq}&~
    \frac{p}{6\theta}\|\bz^{(t)}-\bz^{(t+1)}\|^2 + \tau\cdot\text{dist}\big( \bg(\tilde{\bx}^{(t+1)}), \mathcal{N}_{\mathbb{R}_+^m}(\blambda^{(t+1)})\big)^2 - \frac{12 p\theta M_{\blambda}}{p-L} \text{dist}\big( \bg(\tilde{\bx}^{(t+1)}), \mathcal{N}_{\mathbb{R}_+^m}(\blambda^{(t+1)})\big)\\ 
    \nonumber
    \geq\,&~
    \frac{p}{6\theta}\|\bz^{(t)}-\bz^{(t+1)}\|^2 + \frac{3\tau}{4}\text{dist}\big( \bg(\tilde{\bx}^{(t+1)}), \mathcal{N}_{\mathbb{R}_+^m}(\blambda^{(t+1)})\big)^2 - \frac{1}{\tau} \left(\frac{12 p\theta M_{\blambda}}{p-L}\right)^2\\
    \nonumber
    \geq\,&~
    \frac{p}{12\theta}\|\bz^{(t)}-\bz^{(t+1)}\|^2 + \frac{\tau}{2}\text{dist}\big( \bg(\tilde{\bx}^{(t+1)}), \mathcal{N}_{\mathbb{R}_+^m}(\blambda^{(t+1)})\big)^2 - \frac{1}{\tau} \left(\frac{12 p\theta M_{\blambda}}{p-L}\right)^2 + \frac{R^2(\delta) \theta}{24p} - \frac{p\theta\epsilon_t^2}{12(p-L)^2}  \\ 
    \label{eq:bounded_phi_regularity_result_2_exponent_d}
    \geq\,&~
    \frac{p}{12\theta}\|\bz^{(t)}-\bz^{(t+1)}\|^2 + \frac{\tau}{2}\text{dist}\big( \bg(\tilde{\bx}^{(t+1)}), \mathcal{N}_{\mathbb{R}_+^m}(\blambda^{(t+1)})\big)^2 - \frac{p\theta\epsilon_t^2}{12(p-L)^2},
\end{align}
\end{small}where the second inequality follows from the Young's inequality, the third one holds by \eqref{eq:regular_cond_delta_near_violated-2_exponent_d} and $\theta \leq6p\tau$, and the last one results from $\theta\leq\frac{R^2(\delta) (p-L)^2\tau}{2\cdot 12^3 p^3 M_{\blambda}^2}$.  
Applying \eqref{eq:bounded_phi_regularity_result_2_exponent_d} to \eqref{eq:difference_phi} gives \eqref{eq:bounded_phi_exponent_d}.
\end{proof}

\section{Proof of Theorem~\ref{thm:main_result_exponent_d}}
\label{sec:main_result}
\begin{proof}[Proof of Theorem~\ref{thm:main_result_exponent_d}]
Consider any $t\in\{0,1,\dots,T-1\}$. Firstly, by \eqref{eq:subprob_optim_imela} and the triangle inequality, we have
\begin{align}
\nonumber
    &~\text{dist}\big(-\nabla f(\bx^{(t+1)})-\textstyle\sum_{i=1}^m \lambda_i^{(t+1)}\cdot\nabla g_i(\bx^{(t+1)}), \N_\X(\bx^{(t+1)})\big)\\
    \nonumber
    \leq&~\text{dist}\big(-\nabla_\bx\mathcal{L}_p (\bx^{(t+1)},\bz^{(t)},\blambda^{(t+1)}), \N_\X(\bx^{(t+1)})\big) +p\|\bx^{(t+1)}-\bz^{(t)}\|\\
    \label{eq:epsilon_stat_result_imela_exponent_d}
    \leq&~\epsilon_t+p\|\bx^{(t+1)}-\bz^{(t)}\|
    = \epsilon_t+
    (p/\theta)\cdot\|\bz^{(t+1)}-\bz^{(t)}\|.
\end{align}
Secondly, recalling that $\tilde{\bx}^{(t+1)} = \bx(\blambda^{(t+1)},\bz^{(t)})$, we have
\begin{align}
\nonumber
    \|[\bg(\bx^{(t+1)})]_+\|
    \leq&~\|[\bg(\tilde{\bx}^{(t+1)})]_+\| +\|[\bg(\bx^{(t+1)})]_+-[\bg(\tilde{\bx}^{(t+1)})]_+\|\\
    \nonumber
    \leq&~\|[\bg(\tilde{\bx}^{(t+1)})]_+\| +B_g\|\bx^{(t+1)}-\tilde{\bx}^{(t+1)}\|\\
    \label{eq:epsilon_feas_xtilde_imela_exponent_d}
    \leq&~\|[\bg(\tilde{\bx}^{(t+1)})]_+\| +\frac{B_g\epsilon_t}{p-L}
    \leq\text{dist} \big(\bg(\tilde{\bx}^{(t+1)}), \mathcal{N}_{\mathbb{R}_+^m}(\blambda^{(t+1)})\big) +\frac{B_g\epsilon_t}{p-L},
\end{align}
where the second inequality is by Assumption~\ref{assume:general}E, the third one is from~\eqref{eq:subprob_optimy_imela_corollary} and \eqref{eq:regularity_result_sc_exponent_d} with $\bx = \bx^{(t+1)}$, and the last one holds by $\|[\bg(\tilde{\bx}^{(t+1)})]_+\|= \text{dist}(\bg(\tilde{\bx}^{(t+1)}), \mathbb{R}_-^m)$ and $\mathcal{N}_{\mathbb{R}_+^m}(\blambda^{(t+1)})\subset \mathbb{R}_-^m$.

Lastly, recalling that $I_{g}^{(t+1)}=\text{supp}(\blambda^{(t+1)})$, we have
\begin{align}
\nonumber
    \textstyle\sum_{i=1}^m |\lambda_i^{(t+1)}g_i(\bx^{(t+1)})|
    =&~
    \textstyle\sum_{i\in I_{g}^{(t+1)}} |\lambda_i^{(t+1)}g_i(\bx^{(t+1)})|
    \\
    \nonumber
    \leq&~\textstyle\sqrt{\sum_{i\in I_{g}^{(t+1)}}(\lambda_i^{(t+1)})^2}\sqrt{\sum_{i\in I_{g}^{(t+1)}}g^2_i(\bx^{(t+1)})} \\
    \nonumber
    \leq&~M_{\blambda}\cdot\text{dist} \big(\bg(\bx^{(t+1)}), \mathcal{N}_{\mathbb{R}_+^m}(\blambda^{(t+1)})\big)\\
    \nonumber
    \leq&~M_{\blambda}\cdot\text{dist} \big(\bg(\tilde{\bx}^{(t+1)}), \mathcal{N}_{\mathbb{R}_+^m}(\blambda^{(t+1)})\big) +M_{\blambda}\cdot \|\bg(\bx^{(t+1)})-\bg(\tilde{\bx}^{(t+1)})\|\\ \label{eq:epsilon_comple_slack_result_imela_exponent_d}
    \leq&~M_{\blambda}\cdot\text{dist} \big(\bg(\tilde{\bx}^{(t+1)}), \mathcal{N}_{\mathbb{R}_+^m}(\blambda^{(t+1)})\big) +\frac{M_{\blambda}B_g\epsilon_t}{p-L},
\end{align}
where the first inequality is by the Cauchy--Schwarz inequality, the second one is by Lemma~\ref{lmm:bound_lambda} and the fact that 
$g^2_i(\bx^{(t+1)})=\text{dist}(g_i(\bx^{(t+1)}),\mathcal{N}_{\mathbb{R}_+}(\lambda_i^{(t+1)}))^2$ for $i\in I_{g}^{(t+1)}$, the third one is by the triangle inequality, and the last one is by \eqref{eq:subprob_optimy_imela_corollary}, \eqref{eq:regularity_result_sc_exponent_d} with $\bx = \bx^{(t+1)}$, and Assumption~\ref{assume:general}E.

By summing up \eqref{eq:bounded_phi_exponent_d} for $t=0,1,\dots,T-1$, we have
\begin{align}
\nonumber
   &~\sum_{t=0}^{T-1}\left(\frac{p-L}{4}\|\bx^{(t)}-\tilde{\bx}^{(t+1)}\|^2+\frac{p}{12\theta}\|\bz^{(t)}-\bz^{(t+1)}\|^2
   +\frac{\tau}{2} \text{dist}\big( \bg(\tilde{\bx}^{(t+1)}), \mathcal{N}_{\mathbb{R}_+^m}(\blambda^{(t+1)})\big)^2 \right)\\
    \leq&~\phi^0-\phi^T+\frac{p\theta}{12(p-L)^2}\sum_{t=0}^{T-1}\epsilon_t^2 +T\cdot\Xi_d(\theta)
    \label{eq:constantC_d}
    \leq \phi^0-\underline{f}+\frac{cp\theta\pi^2}{72(p-L)^2} +T\cdot\Xi_d(\theta),
\end{align}
where the second inequality is by \eqref{eq:phi^t_lowerbound} and the fact that
\begin{align}
\label{eq:bounded_sum_epsilon_t_imela_exponent_d}
    \textstyle
    \sum_{t=0}^{T-1}\epsilon_t^2\leq\sum_{t=0}^{\infty}\epsilon_t^2=c\pi^2/6.
\end{align}
According to Definition~\ref{dfn:epsilon_KKT}, denote the squared KKT residual measured at $\bx^{(t+1)}$ by
$$
\mathcal R_{t+1}
:= \left(
\begin{array}{l}
\text{dist}\big(-\nabla f(\bx^{(t+1)})-\sum_{i=1}^m \lambda_i^{(t+1)}\nabla g_i(\bx^{(t+1)}), \N_\X(\bx^{(t+1)})\big)^2\\
+\|[\bg(\bx^{(t+1)})]_+\|^2
+\big(\sum_{i=1}^m |\lambda_i^{(t+1)}g_i(\bx^{(t+1)})|\big)^2
\end{array}
\right).
$$
Squaring and summing up both sides of \eqref{eq:epsilon_stat_result_imela_exponent_d}, the right inequality in \eqref{eq:epsilon_feas_xtilde_imela_exponent_d} and \eqref{eq:epsilon_comple_slack_result_imela_exponent_d} for $t=0,1,\dots,T-1$ leads to 
\begin{footnotesize}
\begin{align}
\nonumber
    \sum_{t=0}^{T-1}\mathcal R_{t+1}
    =&~\sum_{t=0}^{T-1}
    \left(
    \begin{array}{l}
    \text{dist}\big(-\nabla f(\bx^{(t+1)})-\sum_{i=1}^m
    \lambda_i^{(t+1)}\nabla g_i(\bx^{(t+1)}),\N_\X(\bx^{(t+1)})\big)^2\\
    +\|[\bg(\bx^{(t+1)})]_+\|^2
    +\big(\sum_{i=1}^m |\lambda_i^{(t+1)}g_i(\bx^{(t+1)})|\big)^2
    \end{array}
    \right)\\
\nonumber
    \leq&~
    \frac{2p^2}{\theta^2}\sum_{t=0}^{T-1}\|\bz^{(t+1)}-\bz^{(t)}\|^2
    +\left(2+\frac{2(M_{\blambda}^2+1)B_g^2}{(p-L)^2}\right)\sum_{t=0}^{T-1}\epsilon_t^2\\
\nonumber
    &+2(M_{\blambda}^2+1)\sum_{t=0}^{T-1}\text{dist}\big(\bg(\tilde{\bx}^{(t+1)}),\mathcal{N}_{\mathbb{R}_+^m}(\blambda^{(t+1)})\big)^2\\
    \nonumber
    \leq& \left(1+\frac{(M_{\blambda}^2+1)B_g^2}{(p-L)^2}\right)\frac{c\pi^2}{3}+\max\left\{\frac{24p}{\theta},\frac{4(M_{\blambda}^2+1)}{\tau}\right\}\left(\phi^0-\underline{f}+\frac{cp\theta\pi^2}{72(p-L)^2} +T\cdot\Xi_d(\theta) \right)\\
    \label{eq:main_result_constant_bound_exponent_d}
    =&~
    \mathcal C_0(\theta)+\mathcal C_1(\theta) T\cdot\Xi_d(\theta),
\end{align}
\end{footnotesize}where the second inequality is because of \eqref{eq:constantC_d} and \eqref{eq:bounded_sum_epsilon_t_imela_exponent_d}, and the equality is by the definitions of $\mathcal C_0(\theta)$ and $\mathcal C_1(\theta)$ in \eqref{eq:C0C1_def_exponent_d} as follows:
\begin{small}
\begin{align*}
\begin{aligned}
    \mathcal C_0(\theta)
    :=&~ \left(1+\frac{(M_{\blambda}^2+1)B_g^2}{(p-L)^2}\right)\frac{c\pi^2}{3}
    +\max\left\{\frac{24p}{\theta},\frac{4(M_{\blambda}^2+1)}{\tau}\right\}
    \left(\phi^0-\underline{f}+\frac{cp\theta\pi^2}{72(p-L)^2}\right),\\
    \mathcal C_1(\theta)
    :=&~ \max\left\{\frac{24p}{\theta},\frac{4(M_{\blambda}^2+1)}{\tau}\right\}.
\end{aligned}
\end{align*}
\end{small}

From now on, we specify the choice of $\theta$ and the complexity of $T$ separately for $d\in(1,2]$ and $d=1$. 

Firstly, we consider the case $d\in(1,2]$. Fix $\epsilon>0$. Recall from Proposition~\ref{thm:bounded_phi_exponent_d} and the definition of $\tilde{\theta}_d$ in \eqref{eq:tilde_theta_d} that, for any constant $\theta$ satisfying $0<\theta\leq \tilde{\theta}_d$, the sequence generated by Algorithm~\ref{alg:imela} satisfies
$$
\textstyle \sum_{t=0}^{T-1}\mathcal R_{t+1} \leq \mathcal C_0(\theta)+\mathcal C_1(\theta) T\cdot\Xi_d(\theta) 
$$
according to \eqref{eq:main_result_constant_bound_exponent_d}. We then specify 
\begin{align*}
    \theta=\theta(\epsilon)
    :=&~ \min\left\{ \tilde{\theta}_d, \frac{6p\tau}{M_{\blambda}^2+1}, \left(\frac{\epsilon^2}{2A_d}\right)^{d-1}
    \right\},\\
    \text{ where }
    A_d :=&~ 24p\left(1-\frac{1}{d}\right)
    \left(\frac{2}{\tau d}\right)^{\frac{1}{d-1}} (6pC_d)^{\frac{d}{d-1}}.
\end{align*}
Clearly, $\theta(\epsilon)\leq\tilde{\theta}_d$.
Moreover, since $\theta(\epsilon)\leq \frac{6p\tau}{M_{\blambda}^2+1}$, we have
\begin{align}
\label{eq:C1_theta_epsilon_exponent_d}
    \frac{24p}{\theta(\epsilon)}\geq \frac{4(M_{\blambda}^2+1)}{\tau},~
    \text{ and hence }~
    \mathcal C_1(\theta(\epsilon))
    = \max\left\{\frac{24p}{\theta(\epsilon)}, \frac{4(M_{\blambda}^2+1)}{\tau}\right\}
    = \frac{24p}{\theta(\epsilon)}.
\end{align}
Recalling the definition of $\Xi_d(\theta)$ in~\eqref{eq:Xi_d_def}, for $d\in(1,2]$ we have
$$
\Xi_d(\theta)
=
\left(1-\frac{1}{d}\right)\left(\frac{2}{\tau d}\right)^{\frac{1}{d-1}}
\big(6p\theta C_d\big)^{\frac{d}{d-1}},
$$
which implies
\begin{align*}
    \mathcal C_1(\theta(\epsilon)) \cdot\Xi_d(\theta(\epsilon))
    =&~ \frac{24p}{\theta(\epsilon)} \left( 1-\frac{1}{d} \right) \left( \frac{2}{\tau d} \right)^{\frac{1}{d-1}} \big( 6p\theta(\epsilon)\,C_d \big)^{\frac{d}{d-1}}\\
    =&~ 24p\left( 1-\frac{1}{d} \right) \left( \frac{2}{\tau d} \right)^{\frac{1}{d-1}} (6pC_d)^{\frac{d}{d-1}} \cdot \theta(\epsilon)^{\frac{1}{d-1}}
    = A_d\cdot \theta(\epsilon)^{\frac{1}{d-1}}.
\end{align*}
By the definition of $\theta(\epsilon)$, we also have
\begin{align}
\label{eq:theta_epsilon_exponent_d}
    \theta(\epsilon)\leq \left(\frac{\epsilon^2}{2A_d}\right)^{d-1},~
    \text{ hence }~
    \theta(\epsilon)^{\frac{1}{d-1}} \leq \frac{\epsilon^2}{2A_d}.
\end{align}
Substituting \eqref{eq:theta_epsilon_exponent_d} into the previous result gives
\begin{align}
\label{eq:C1_Xi_theta_epsilon_exponent_d}
    \mathcal C_1(\theta(\epsilon)) \cdot \Xi_d(\theta(\epsilon))
    = A_d\cdot \theta(\epsilon)^{\frac{1}{d-1}}
    \leq A_d\cdot \frac{\epsilon^2}{2A_d}
    = \frac{\epsilon^2}{2}.
\end{align}
Next, denote
\begin{align}
\label{eq:B0_exponent_d}
    B_0:= \left(1+\frac{(M_{\blambda}^2+1)B_g^2}{(p-L)^2}\right)\frac{c\pi^2}{3} +\frac{24cp^2\pi^2}{72(p-L)^2},
    \quad
    B_1:=B_0+24p \big( \phi^0-\underline{f} \big),
\end{align}
which are independent of $\epsilon$. Recalling from \eqref{eq:C0C1_def_exponent_d}, we obtain
\begin{align}
\nonumber
    \mathcal C_0(\theta(\epsilon))
    =&~ \left( 1+\frac{(M_{\blambda}^2+1)B_g^2}{(p-L)^2} \right) \frac{c\pi^2}{3}
    + \mathcal C_1(\theta(\epsilon))\left(\phi^0-\underline{f}+\frac{cp\theta(\epsilon)\pi^2}{72(p-L)^2}\right)\\
    \nonumber
    =&~ \left( 1+\frac{(M_{\blambda}^2+1)B_g^2}{(p-L)^2} \right) \frac{c\pi^2}{3}
    + \frac{24p}{\theta(\epsilon)} \left(\phi^0-\underline{f}+\frac{cp\theta(\epsilon)\pi^2}{72(p-L)^2}\right)\\
    \label{eq:C0_theta_epsilon_exponent_d}
    =&~B_0 + \frac{24p(\phi^0-\underline{f})}{\theta(\epsilon)}
    \leq \frac{B_1}{\theta(\epsilon)},
\end{align}
where the second equality is by \eqref{eq:C1_theta_epsilon_exponent_d}, the third equality is by \eqref{eq:B0_exponent_d}, and  the last inequality is by \eqref{eq:B0_exponent_d} and the fact that $0<\theta(\epsilon)\leq 1$.  We next bound $1/\theta(\epsilon)$ using the concrete choice of $\theta(\epsilon)$.
Denote
$$
\overline{\theta}_d:=\min\left\{\tilde{\theta}_d,\ \frac{6p\tau}{M_{\blambda}^2+1}\right\}>0.
$$
If $\big(\frac{\epsilon^2}{2A_d}\big)^{d-1}\leq \overline{\theta}_d$, then $\theta(\epsilon)=\big(\frac{\epsilon^2}{2A_d}\big)^{d-1}$ and hence
$$
\frac{1}{\theta(\epsilon)}=\left(\frac{2A_d}{\epsilon^2}\right)^{d-1}.
$$
Otherwise, $\theta(\epsilon)=\overline{\theta}_d$ and thus
$$
\frac{1}{\theta(\epsilon)}=\frac{1}{\overline{\theta}_d}
\leq \frac{1}{\overline{\theta}_d}\left(\frac{1}{\epsilon^2}\right)^{d-1},
$$
where the inequality holds since $\epsilon\in(0,1)$ and $d>1$ imply $(1/\epsilon^2)^{d-1}\geq 1$.
Thus, we conclude that
\begin{align}
\label{eq:theta_epsilon_inverse_exponent_d}
    \frac{1}{\theta(\epsilon)} \leq \frac{\max\{ 1/\overline{\theta}_d, (2A_d)^{d-1} \}}{\epsilon^{2(d-1)}},
    \quad 
    \mathcal C_0(\theta(\epsilon)) \leq \frac{B_1 \max\{ 1/\overline{\theta}_d, (2A_d)^{d-1} \}}{\epsilon^{2(d-1)}}.
\end{align}
Hence, we obtain
$$
T := \left\lceil\frac{2\mathcal C_0(\theta(\epsilon))}{\epsilon^2}\right\rceil \leq \left\lceil \frac{2B_1 \max\{ 1/\overline{\theta}_d, (2A_d)^{d-1} \}}{\epsilon^{2(d-1)}} \Big/\epsilon^2\right\rceil = O(\epsilon^{-2d}),
$$
where the inequality is by \eqref{eq:C0_theta_epsilon_exponent_d} and \eqref{eq:theta_epsilon_inverse_exponent_d}. Substituting $\theta=\theta(\epsilon)$ into \eqref{eq:main_result_constant_bound_exponent_d}, dividing both sides by $T$, and using the definition of $T$ together with \eqref{eq:C1_Xi_theta_epsilon_exponent_d} yield
$$
\frac{1}{T}\sum_{t=0}^{T-1}\mathcal R_{t+1}
\leq
\frac{\mathcal C_0(\theta(\epsilon))}{T}
+\mathcal C_1(\theta(\epsilon))\,\Xi_d(\theta(\epsilon))
\leq
\frac{\epsilon^2}{2}+\frac{\epsilon^2}{2}
=
\epsilon^2.
$$
Therefore, there exists $s\in\{0,1,\dots,T-1\}$ such that $\mathcal R_{s+1}\leq \epsilon^2$, namely,
\begin{equation}
\label{eq:main_result_epsilon_KKT_exponent_d}
\begin{aligned}
    &~\text{dist}\big(-\nabla f(\bx^{(s+1)}) -\textstyle\sum_{i=1}^m \lambda_i^{(s+1)}\cdot\nabla g_i(\bx^{(s+1)}), \N_\X(\bx^{(s+1)})\big)^2\\
    &+\|[\bg(\bx^{(s+1)})]_+\|^2 +\big(\textstyle\sum_{i=1}^m |\lambda_i^{(s+1)}g_i(\bx^{(s+1)})|\big)^2~\leq~\epsilon^2.
\end{aligned}
\end{equation}
According to Definition~\ref{dfn:epsilon_KKT}, $\bx^{(s+1)}$ is an $\epsilon$-KKT point of \eqref{eq:gco_general}.

Secondly, we consider the case $d=1$. Fix $\epsilon>0$. By setting $d=1$ in Proposition~\ref{thm:bounded_phi_exponent_d}, we have $\Xi_{1}(\theta)=0$ by \eqref{eq:Xi_d_def}. Moreover, $C_1=\frac{\gamma^{1}K}{p-L}=\frac{\gamma K}{p-L}$, and the quantity $\hat{\theta}_d$ defined in \eqref{eq:hat_theta_d} simplifies to
$$
\hat{\theta}_1 = \frac{\tau}{12p C_1}\left(\frac{\tau}{12p}\right)^{1-1} = \frac{\tau}{12p C_1} = \frac{(p-L)\tau}{12p\gamma K}.
$$
Consequently, the admissible constant stepsize upper bound $\tilde{\theta}_d$ in \eqref{eq:tilde_theta_d} specializes to
$$
\tilde{\theta}_1
= \min\left\{\frac{p-L}{18 p}, 6p\tau, \frac{R^2(\delta) (p-L)^2\tau}{2\cdot 12^3 p^3 M_{\blambda}^2}, \frac{(p-L)\tau}{12p\gamma K}\right\},~
\text{ with }~
\tau=\frac{p-L}{4B_g^2}.
$$
Fix any constant $\theta$ satisfying $0<\theta\leq \tilde{\theta}_1$, and set $\theta_t=\theta$ for all $t\geq 0$.
Then the bound \eqref{eq:main_result_constant_bound_exponent_d} reduces to
$$
\sum_{t=0}^{T-1}\mathcal R_{t+1}
\leq
\mathcal C_0(\theta),
$$
where $\mathcal C_0(\theta)$ is independent of $\epsilon$. Analogous to the case $d\in(1,2]$, choosing
$$
T:=\left\lceil \frac{2\mathcal C_0(\theta)}{\epsilon^2}\right\rceil = O(\epsilon^{-2})
$$
ensures that
$$
\frac{1}{T}\sum_{t=0}^{T-1}\mathcal R_{t+1}
\leq \frac{\mathcal C_0(\theta)}{T}
\leq \frac{\epsilon^2}{2} \leq \epsilon^2.
$$
Therefore, there exists $s\in\{ 0,1,\dots,T-1 \}$ such that $\mathcal R_{s+1}\leq \epsilon^2$, and hence $\bx^{(s+1)}$ is an $\epsilon$-KKT point of \eqref{eq:gco_general}. In this case, the total number of outer iterations satisfies $T= O(\epsilon^{-2})=O(\epsilon^{-2d})$ with $d=1$.

By \eqref{eq:subprob_complexity}, the complexity of the $t$-th iteration of Algorithm~\ref{alg:imela} is $O\left(\sqrt{\frac{K}{p-L}}\ln(\epsilon_t^{-1})\right) = O(\ln(\epsilon^{-1}))$ for $t=0,1,\dots,T-1$. Therefore, the total oracle complexity of Algorithm~\ref{alg:imela} for finding an $\epsilon$-KKT point is $O(\epsilon^{-2d}\ln(\epsilon^{-1}))$. 
\end{proof}

\section{Additional Details of Numerical Experiments}
\label{sec:additional_exp}
In this section, we provide additional details of the numerical experiments presented in Section~\ref{sec:experiment}.

\subsection{Details of Datasets}
\label{sec:dataset}
We solve problem \eqref{eq:DPfairnessclassification_linear} on three datasets: \textit{a9a}~\cite{kohavi1996scaling}, \textit{bank}~\cite{moro2014data} and \textit{COMPAS}~\cite{angwin2016compas}.  Details about these datasets are given in 
Table~\ref{tbl:data}. Each dataset is splited into two subsets with a ratio of 2 : 1.
The larger subset serves as $\mathcal{D}$ in the constraint, while the smaller subset is further partitioned into $\mathcal{D}_p$ and $\mathcal{D}_u$ based on the binary group variable specified in Table~\ref{tbl:data}.

\begin{table}[h]
\caption{Dataset information. The binary group variable represents males or females in the a9a dataset, users with age within $[25,60]$ and outside $[25,60]$ in the bank dataset, and Caucasian versus non-Caucasian individuals in the COMPAS dataset.}
\begin{center}
\begin{sc}
\begin{tabular}{c|ccccc}
\hline
Datasets & $n$ & $d$ & Label & Groups \\
\hline
a9a & 48,842 & 123 & Income & Gender\\
Bank & 41,188 & 54 & Subscription & Age\\
COMPAS & 6,172 & 16 & Recidivism & Race\\
\hline
\end{tabular}
\end{sc}
\end{center}
\vskip 0.05in
\label{tbl:data}
\vskip -0.2in
\end{table}

\subsection{Local LICQ and EBC for the DP classification instance \eqref{eq:DPfairnessclassification_linear}}
\label{sec:DP_LICQ}
\begin{proposition}
\label{prop:DP_LICQ}
Consider problem~\eqref{eq:DPfairnessclassification_linear} with $\X=\{\bx\in\mathbb{R}^d:\|\bx\|_1\leq r\}$. Suppose 
$$
\left\{\bx\in\mathbb{R}^d: \mathcal{L}(\bx)\leq \mathcal{L}^*+\kappa\right\}\subset
\left\{\bx\in\mathbb{R}^d: \|\bx\|_1\leq r/2\right\}.
$$
Then~\eqref{eq:DPfairnessclassification_linear} satisfies the local LICQ condition~\eqref{eq:LICQ_sigmamin} at every KKT point. Consequently, Assumption~\ref{assume:local_EBC_exponent_d} holds with exponent $d=1$.
\end{proposition}
\begin{proof}
Let $\bx^*\in\X^*$ be any KKT point of~\eqref{eq:DPfairnessclassification_linear}. 
By the hypothesis of this proposition, the set constraint $\X$ is inactive in a neighborhood of $\bx^*$ and $J_A(\bx^*)=\emptyset$ in \eqref{eq:LICQ_sigmamin}.

If $\mathcal{L}(\bx^*)<\mathcal{L}^*+\kappa$, then no inequality constraint is active
at $\bx^*$, i.e., $J_g(\bx^*)=\emptyset$ in \eqref{eq:LICQ_sigmamin}, and the LICQ condition holds trivially.

Suppose next that $\mathcal{L}(\bx^*)=\mathcal{L}^*+\kappa$. Since the set constraint $\X$ is inactive, the only active constraint is $\mathcal{L}(\bx)\leq \mathcal{L}^*+\kappa$. In this case, LICQ \eqref{eq:LICQ_sigmamin} reduces to the requirement that there exists $\zeta>0$ such that $\|\nabla \mathcal{L}(\bx^*)\| \geq \zeta$. Because $\mathcal{L}$ is continuous and $\X$ is compact, the restricted level set
$$
\mathcal{C}_\kappa:=\{\bx\in\X:\mathcal{L}(\bx)=\mathcal{L}^*+\kappa\}
$$
is compact whenever it is nonempty.
Therefore, the quantity
$$
\xi := \min\big\{\|\nabla \mathcal{L}(\bx)\| :
\bx\in\X,\ \mathcal{L}(\bx)=\mathcal{L}^*+\kappa\big\}
$$
is attained.

Moreover, we assert that $\xi>0$, since otherwise there would exist $\bar{\bx}\in\mathcal{C}_\kappa$ satisfying
$$
\nabla \mathcal{L}(\bar{\bx})=\mathbf{0}
~\text{ and }~
\mathcal{L}(\bar{\bx})=\mathcal{L}^*+\kappa,
$$
which contradicts the definition of $\mathcal{L}^*$ as the minimal value of $\mathcal{L}$ over $\X$.

Therefore, the local LICQ condition~\eqref{eq:LICQ_sigmamin} holds at every KKT point.
The claim follows from Proposition~\ref{thm:LICQ_implies_local_EBC}.
\end{proof}

\subsection{Implementation Details of the Methods in Comparison}
\label{sec:alg_details}
Recall that, at its $t$-th iteration, the iMELa method obtains 
\begin{align}
\label{eq:imela_subprob}
    \bx^{(t+1)} \approx \argmin_{\bu\in\X}\left\{f(\bu)+{\textstyle\sum_{i=1}^m\lambda_i^{(t+1)}\cdot g_i(\bu)}+\frac{p}{2}\|\bu-\bz^{(t)}\|^2\right\}.
\end{align}
Similarly, at its $t$-th iteration, the iPPP method approximately solves the subproblem
\begin{align}
\label{eq:ippp_subprob}
    \bx^{(t+1)} \approx \argmin_{\bu\in\X}\left\{f(\bu)+\frac{\rho_t}{2}{\textstyle\sum_{i=1}^m\left(\big[g_i(\bu)]_+\right)^2}+\frac{p_t}{2}\|\bu-\bx^{(t)}\|^2\right\},
\end{align}
where $\rho_t>0$ is the penalty parameter and $p_t>0$ is the proximal parameter.
In addition, DPALM follows a proximal augmented Lagrangian (AL) framework with a possibly damped dual step-size.
At its $t$-th iteration, DPALM approximately solves the subproblem
\begin{align}
\label{eq:dpalm_subprob}
    \bx^{(t+1)}
    \approx
    \argmin_{\bu\in\X}
    \left\{ f(\bu) + \frac{p}{2}\|\bu-\bx^{(t)}\|^2 + \frac{1}{2\beta_t} {\textstyle\sum_{i=1}^m} \Big( \big[ \lambda_i^{(t)}+\beta_t g_i(\bu) \big]_+^2 - \big( \lambda_i^{(t)} \big)^2 \Big)
    \right\},
\end{align}
where $p>0$ is the proximal parameter and $\beta_t>0$ is the AL penalty parameter. The (possibly damped) dual step-size is chosen as
$$
\alpha_t := \min\left\{ \beta_t, \frac{v_t}{ \|[\bg(\bx^{(t+1)})]_+\| } \right\},
$$
and the Lagrangian multipliers are  updated by
$$
\lambda_i^{(t+1)} = \lambda_i^{(t)} + \alpha_t\max\left\{-\frac{\lambda_i^{(t)}}{\beta_t}, g_i(\bx^{(t+1)})\right\},\; \forall\,i\in[m].
$$

For a fair comparison, we apply the accelerated projected gradient (APG) method by~\citet[Eq.\,(2.2.63)]{nesterov2018lectures} to approximately solve \eqref{eq:imela_subprob}, \eqref{eq:ippp_subprob}, and \eqref{eq:dpalm_subprob}. Let $F(\bu)$ be the objective function in \eqref{eq:imela_subprob}, \eqref{eq:ippp_subprob}, or \eqref{eq:dpalm_subprob} and let $\bu^{(k)}$ be the main iterate generated by the APG method. For the task $\min_{\bu\in\X} F(\bu)$ where $F$ is $\mu_F$-strongly convex and $L_F$-smooth in $\bu$, we denote $q_F:=\mu_F/L_F$ and then present the APG method~\citep[Eq.\,(2.2.63)]{nesterov2018lectures} in Algorithm~\ref{alg:apg}. The APG method is terminated when 
\begin{align}
\label{eq:subprob_stopping}
    \|\mathcal{G}_{1/\eta_t}(\bu^{(k)})\| :=\big\|\eta_t^{-1}\cdot\big( \bu^{(k)}-T_{1/\eta_t}(\bu^{(k)})\big)\big\|\leq\epsilon_t',
\end{align}
where $T_{1/\eta_t}(\bu^{(k)}):=\text{proj}_{\X} \left(\bu^{(k)}-\eta_t\nabla F(\bu^{(k)})\right)$,
$\eta_t$ is a step-size and $\epsilon_t'$ is a targeted tolerance. Then we take $\bx^{(t+1)}=T_{1/\eta_t}(\bu^{(k)})$ as 
the approximate solution to \eqref{eq:imela_subprob}, \eqref{eq:ippp_subprob}, or \eqref{eq:dpalm_subprob}. Here, $\mathcal{G}_{1/\eta_t}(\bu^{(k)})$ is known as the gradient mapping of $F(\cdot)$ at $\bu^{(k)}$. By~\citet[Eq.\,(2.15)]{nesterov2013gradient}, \eqref{eq:subprob_stopping} ensures 
$\text{dist}\left( -\nabla F(\bx^{(t+1)}),\N_\X(\bx^{(t+1)})\right)\leq\epsilon_t$ with some $\epsilon_t=O(\epsilon_t')$, i.e., condition \eqref{eq:subprob_optim_imela} holds if $F$ is from \eqref{eq:imela_subprob}.

\begin{algorithm}[t]
\caption{Accelerated Projected Gradient (APG) Method~\citep[Eq.\,(2.2.63)]{nesterov2018lectures} for $\min_{\bu\in\X} F(\bu)$}
\label{alg:apg}
\begin{algorithmic}[1]
    \STATE {\bfseries Input:} targeted tolerance $\epsilon_t'>0$, step-size $\eta_t>0$.
    \STATE {\bfseries Initialization:} $\bu^{(0)} =\bv^{(0)} =\bx^{(t)}$.
    \FOR{iteration $k=0,1,2,\dots$}
        \STATE $\tilde{\bu}^{(k)} = T_{1/\eta_t}(\bu^{(k)})$.
        \IF{ $\|\mathcal{G}_{1/\eta_t}(\bu^{(k)})\|=\eta_t^{-1}\|\bu^{(k)}-\tilde{\bu}^{(k)}\|\leq\epsilon_t'$ }
            \STATE {\bfseries Output:} $\bx^{(t+1)}=\tilde{\bu}^{(k)}$.
            \STATE {\bfseries Stop.}
        \ENDIF
        \STATE
        $\bu^{(k+1)} =\text{proj}_{\X} (\bv^{(k)}-\eta_t\nabla F(\bv^{(k)}))$.
        \STATE
        $\bv^{(k+1)}= \bu^{(k+1)} + \frac{1-\sqrt{q_F}}{1+\sqrt{q_F}} (\bu^{(k+1)} - \bu^{(k)})$.
    \ENDFOR
\end{algorithmic}
\end{algorithm}

In addition, at its $t$-th iteration, SP-LM updates $\bx$ by
\begin{align}
\label{eq:sp_lm_x_update}
    \bx^{(t+1)} =\text{proj}_{\X} \left(\bx^{(t)}-\eta_t\nabla_{\bx}\left( f(\bx^{(t)})+{\textstyle\sum_{i=1}^m\lambda_i^{(t+1)}\cdot g_i(\bx^{(t)})}+\frac{p}{2}\|\bx^{(t)}-\bz^{(t)}\|^2\right)\right),
\end{align}
and, at its $t$-th iteration, the SSG method updates $\bx$ by
\begin{small}
\begin{align}
\label{eq:ssg_x_update}
    \bx^{(t+1)}=
    \begin{cases}
        \text{proj}_{\X} (\bx^{(t)}-\eta_t\nabla f(\bx^{(t)}))\text{ if }\max_{i\in[m]}g_i(\bx^{(t)})\leq\epsilon_t, \\[0.5ex]
        \text{proj}_{\X} (\bx^{(t)}-\eta_t\bzeta_g^{(t)})\text{ for some }\bzeta_g^{(t)}\in\mathrm{conv}(\{\nabla g_j(\bx^{(t)}):j\in I(\bx^{(t)})\})\text{ otherwise,}
    \end{cases}
\end{align}
\end{small}where $I(\bx):=\{j\in[m]:g_j(\bx)=\max_{i\in[m]}g_i(\bx)\}$.

In the experiments, we first solve \eqref{eq:ermL_linear} by the projected gradient method up to near optimality to obtain a close approximation of $\mathcal{L}^*$ and a solution $\bx_{\text{feas}}$. Specifically, we use a constant step-size of $0.1$ in the projected gradient method and terminate it when it finds a solution $\bx_{\text{feas}}$ satisfying $\text{dist}(-\nabla\mathcal{L}(\bx_{\text{feas}}),\N_\X(\bx_{\text{feas}}))\leq0.001$. Then we set $\kappa=0.001\mathcal{L}^*$ in \eqref{eq:DPfairnessclassification_linear}.

The methods in comparison are initialized at $\bx^{(0)}=\bx_{\text{feas}}$ obtained above. Based on the definitions of $\mathcal{R}$ and $\mathcal{L}$, we can calculate an upper bound of the smooth parameter of the objective and constraint functions in \eqref{eq:DPfairnessclassification_linear} using the data. 
For the objective function in \eqref{eq:DPfairnessclassification_linear}, we notice that $\mathcal{R}(\bx)\in[0,1]$ for any $\bx\in\mathbb{R}^d$ and, by a derivation similar to~\citet[Eq.\,(75--77)]{huang2023oracle}, $\mathcal{R}(\bx)$ and $\nabla \mathcal{R}(\bx)$ are Lipschitz continuous with constants of
\begin{align*}
    \alpha=\frac{1}{4n_p}\sum_{i=1}^{n_p} \|\ba_i^p\|+\frac{1}{4n_u}\sum_{i=1}^{n_u} \|\ba_i^u\|\text{ and }
    \beta=\frac{1}{4n_p}\sum_{i=1}^{n_p} \|\ba_i^p\|^2+\frac{1}{4n_u}\sum_{i=1}^{n_u} \|\ba_i^u\|^2,
\end{align*}
respectively. 
Additionally, it is easy to show that $\frac{1}{2}(\mathcal{R}(\bx))^2$ is $\alpha$-Lipschitz continuous and $(\beta+\alpha^2)$-smooth in $\bx$. For the constraint function in \eqref{eq:DPfairnessclassification_linear}, by the fact that $\ell(z)$ is 1-Lipschitz continuous and $\sigma(z)$ is $\frac{1}{4}$-Lipschitz continuous,  we can show that $\mathcal{L}(\bx)$ and $\nabla\mathcal{L}(\bx)$ are Lipschitz continuous with constants of
\begin{align*}
    \gamma=\frac{1}{n}\sum_{i=1}^n\|\ba_i\|
    \text{ and }
    \gamma'=\frac{1}{4n}\sum_{i=1}^n\|\ba_i\|^2,
\end{align*}
respectively. 
Therefore, \eqref{eq:DPfairnessclassification_linear} satisfies Assumption~\ref{assume:general} with any $L\geq\max\{\beta+\alpha^2,\gamma'\}$. After estimating $\max\{\beta+\alpha^2,\gamma'\}$ using data, we set $L=10$, 10, and 2.5 on \textit{a9a}, \textit{bank}, and \textit{COMPAS}, respectively, and set $p=2L$. In the APG method applied to \eqref{eq:imela_subprob}, \eqref{eq:ippp_subprob}, and \eqref{eq:dpalm_subprob}, we use a constant step-size $\eta$. On \textit{a9a} and \textit{bank}, we select $\eta$ from $\{0.005,0.01,0.02,0.05\}$; on \textit{COMPAS}, we select $\eta$ from $\{0.02,0.05,0.1,0.2\}$.
The choices of other parameters are uniform across all datasets. For the iMELa method, we select $\tau_t=\tau$ from $\{5,10,20,50\}$ and $\theta_t=\theta$ from $\{0.5,0.75,1\}$, and set $\epsilon_t'=\frac{c}{t+1}$ with $c$ selected from $\{1,2,5,10\}$. For the iPPP method, following~\citet[Eq.\,(5.14)]{lin2022complexity}, we set $p_t\equiv p$, $\rho_t=\rho\sqrt{t+1}$ and $\epsilon_t'=\frac{1}{\rho_t(t+1)}$ with $\rho$ selected from $\{200,500,1000,1500\}$. For DPALM, following~\citet[Eq.\,(21)]{dahal2023damped}, we set $\beta_t=\beta_0\sqrt{t+1}$ and $v_k=\frac{v_0}{\sqrt{t+1}(\ln(t+1))^2}$ with $\beta_0$ selected from $\{0.02, 0.05, 0.1, 0.2\}$ on \textit{a9a} and \textit{bank} and $\{10^{-4},2\times10^{-4}, 5\times10^{-4},10^{-3}\}$ on \textit{COMPAS} and $v_0$ selected from $\{50, 100, 150, 200\}$; following \citet[Eq.\,(28) \& Sec. 4.5]{dahal2023damped}, we set $\epsilon_t'=\min\{\frac{\epsilon'}{8},\frac{1}{2}\sqrt{p/\beta_t},1\}$ with $\epsilon':=10^{-2}$ according to \citet[Sec. 4.5]{dahal2023damped}. For SP-LM, we select $\eta$ from $\{0.005,0.01,0.02,0.05\}$ on \textit{a9a} and \textit{bank} and from $\{0.02,0.05,0.1,0.2\}$ on \textit{COMPAS}. We then select $\tau_t=\tau$ from $\{5,10,20,50\}$ and $\theta_t=\theta$ from $\{0.5,0.75,1\}$. For the SSG method, following~\citet{huang2023oracle}, we adopt both static and diminishing step-sizes. For the static step-size, we select $\epsilon_t'=\epsilon$ from $\{10^{-6}, 2\times10^{-6},5\times10^{-6}, 10^{-5}\}$ and $\eta_t=\eta$ from $\{2\times10^{-4}, 5\times10^{-4},10^{-3}, 2\times10^{-3}\}$. For the diminishing step-size, we select $\epsilon_t'=\frac{E_1}{\sqrt{t+1}}$ and $\eta_t=\frac{E_2}{\sqrt{t+1}}$ and select $E_1$ from $\{5\times10^{-5}, 10^{-4},2\times10^{-4}, 5\times10^{-4}\}$ and $E_2$ from $\{0.02,0.05,0.1,0.2\}$. Since the iMELa method, DPALM, and SP-LM are primal-dual based, and the primal-based iPPP method can employ $\blambda^{(t)}=\rho_{t-1}[\bg(\bx^{(t)})]_+$ as its auxiliary dual estimates according to Step 5 of Algorithm 1 in~\citet{lin2022complexity}, their best combination of parameters is chosen as the one that minimizes the smallest value of 
\begin{align*}
    &~\text{dist}\big(-\mathcal{R}(\bx^{(s)})\cdot\nabla \mathcal{R}(\bx^{(s)}) -\lambda^{(s)}\cdot\nabla \mathcal{L}(\bx^{(s)}), \N_\X(\bx^{(s)})\big)\\
    &+\big\|\big[\mathcal{L}(\bx^{(s)})-(\mathcal{L}^*+\kappa)\big]_+\big\| + \big|\lambda^{(s)}\big(\mathcal{L}(\bx^{(s)})-(\mathcal{L}^*+\kappa)\big)\big|,
\end{align*}
among iterations $s=0,1,\dots,T_{\text{tuning}}-1$. Here, $T_{\text{tuning}}$ is the number of iterations for parameter tuning. 
Since the SSG method does not generate $\blambda^{(t)}$, its best combination of parameters is identified as the one that minimizes the smallest value of $\frac{1}{2}(\mathcal{R}(\bx^{(s)}))^2$ when $\mathcal{L}(\bx^{(s)})\leq \mathcal{L}^*+\kappa+10^{-5}$ among iterations $s=0,1,\dots,T_{\text{tuning}}-1$. For SP-LM and the SSG method, we set $T_{\text{tuning}}=1000$, 5000 and 75000 for \textit{a9a}, \textit{bank} and \textit{COMPAS}, respectively. For a fair comparison, for the iMELa method, the iPPP method, and DPALM, we set $T_{\text{tuning}}$ to be the smallest integer satisfying $\sum_{t=0}^{T_{\text{tuning}}-1}k_t\geq1000\text{, }5000\text{ and }75000$ for \textit{a9a}, \textit{bank} and \textit{COMPAS}, respectively, where $k_t$ is the number of inner steps performed at the $t$-th outer iteration. Once all parameters are chosen, we run SP-LM and the SSG method for $T=4T_{\text{tuning}}$ iterations, and run the iMELa method, the iPPP method, and DPALM for $T$ iterations with $T$ being the smallest integer satisfying $\sum_{t=0}^{T-1} k_t\geq4\sum_{t=0}^{T_{\text{tuning}}-1}k_t$. 

\end{document}